\documentclass[12pt]{article}
\usepackage{amsmath,mathrsfs}
\usepackage{graphicx}
\usepackage{enumerate}
\usepackage{natbib} 
\usepackage{algorithm}
\usepackage{algorithmicx}
\usepackage[noend]{algpseudocode}
\usepackage{makecell}

\usepackage{tikz}
\usepackage{xparse}

\NewDocumentCommand{\circnum}{O{-0.5pt} m}{%
  \tikz[baseline=(char.base)]{%
    \node[shape=circle,draw,inner sep=#1] (char) {#2};%
  }%
}

\newcommand{\blind}{1}

\usepackage{amsmath,amssymb,amsfonts,amsthm,bbm}

\usepackage{mathtools}   
\usepackage{stmaryrd} 

\usepackage{enumitem}   
\usepackage{romannum}

\usepackage[titletoc,title]{appendix}  

\usepackage[table,xcdraw,dvipsnames]{xcolor}   
\usepackage{hyperref}
\newcommand\myshade{85}
\colorlet{mylinkcolor}{YellowOrange}
\colorlet{mycitecolor}{Aquamarine}
\colorlet{myurlcolor}{violet}

\hypersetup{
  linkcolor  = mylinkcolor!\myshade!black,
  citecolor  = mycitecolor!\myshade!black,
  urlcolor   = myurlcolor!\myshade!black,
  colorlinks = true,
}

\usepackage{caption}
\usepackage{subcaption}

\usepackage[normalem]{ulem}
\usepackage{setspace}

\usepackage{graphicx}
\usepackage{comment}
\graphicspath{{figures/}}

\renewcommand{\hat}{\widehat}
\renewcommand{\tilde}{\widetilde}

\def\one {\mathbbm{1}}

\providecommand{\abs}[1]{\left\lvert#1\right\rvert}

\usepackage{mathtools}
\DeclarePairedDelimiterX{\infdivx}[2]{(}{)}{%
  #1 \; \delimsize\| \; #2%
}

\DeclareMathOperator*{\argmin}{argmin}

\newcommand*\xbar[1]{%
  \hbox{%
    \vbox{%
      \hrule height 0.4pt 
      \kern0.5ex
      \hbox{%
        \kern-0em
        \ensuremath{#1}%
        \kern-0em
      }%
    }%
  }%
} 
\newtheorem{definition}{Definition}
\newtheorem{assumption}[definition]{Assumption}
\newtheorem{lemma}[definition]{Lemma}

\newtheorem{theorem}[definition]{Theorem}

\newtheorem{example}[definition]{Example}

\theoremstyle{definition}

\newcommand{\E}{\mathbb{E}}
\newcommand{\Prob}{\mathbb{P}}

\definecolor{royalpurple}{rgb}{0.47, 0.32, 0.66}
\definecolor{greenfresh}{HTML}{00897B}
\definecolor{bluefresh}{HTML}{1E88E5}
\definecolor{redfresh}{HTML}{E53935}

\definecolor{royalpurple}{rgb}{0.47, 0.32, 0.66}

\providecommand\Halmos{\rule{1ex}{1.2ex}}

\def\beq{\begin{equation}}
\def\eeq{\end{equation}}

\def\bet{\begin{theorem}}
\def\eet{\end{theorem}}

\def\bel{\begin{lemma}}
\def\eel{\end{lemma}}

\begin{document}
\pagenumbering{arabic}

\def\spacingset#1{\renewcommand{\baselinestretch}%
{#1}\small\normalsize} \spacingset{1}

%
%
%

\def\TITLE{Optimizing Server Locations in Spatial Queues: \\ Parametric and Nonparametric Bayesian Optimization Approaches}

\if1\blind
{
  \title{\Large\bf \TITLE}
  \author{
  Cheng Hua$^{\flat *}$ \hspace{1.5ex}
  Arthur J. Swersey$^{\natural *}$ \hspace{1.5ex}
  Wenqian Xing$^{\sharp *}$ \hspace{1.5ex}
  Yi Zhang$^{\dag}$\footnote{Author names are listed in alphabetical order.} \\ \normalsize
  \medskip
  {\footnotesize  
  $^{\flat}$Shanghai Jiao Tong University \hspace{2ex}
  $^\natural$ Yale University \hspace{2ex}
  $^\sharp$Stanford University \hspace{2ex}
  $^{\dag}$Columbia University
  }
  }
  \maketitle
} \fi

\if0\blind
{
  \bigskip
  \bigskip
  \bigskip
  \begin{center}
    {\LARGE\bf \TITLE}
\end{center}
  \medskip
} \fi

\bigskip

\begin{abstract}
\spacingset{1.68}

This paper presents a new model for solving the optimal server location problem in a spatial hypercube queueing model. Unlike deterministic location models, our approach accounts for server availability, varying utilization levels, and dependencies across servers. We prove that the problem is NP-hard and establish lower and upper bounds, as well as asymptotic results, by relating it to special cases of the classical $p$-Median problem. To address the computational challenge, we propose two Bayesian optimization approaches: (i) a parametric approach based on a sparse Bayesian linear model with second-order interactions, and (ii) a nonparametric approach using a Gaussian process surrogate with the $p$-Median objective as the prior mean function. We prove that both methods achieve sublinear regret and converge to the optimal solution. 
Numerical experiments and a case study using real-world data from the St. Paul, Minnesota, emergency response system show that our approaches consistently identify optimal solutions and outperform all baseline methods. 

\end{abstract}

\noindent%


\newpage
\spacingset{1.9} 

\addtolength{\textheight}{.1in}%

\section{Introduction} 
Over the past 50 years, a large literature has appeared on operations research models for emergency services deployment~\citep{green2004anniversary}. One of the most often cited and widely applied models is the hypercube queueing model for spatial queues~\citep{larson74hypercube, larson75approx}. For discussions of the early literature, including the hypercube model, see~\cite{kolesar1986deployment},~\cite{swersey1994deployment}, and~\cite{iannoni2023review}. The hypercube framework models \textit{server-to-customer} systems in which units travel to incident locations, typical in fire, police, and emergency medical service (EMS) operations. Calls arrive at atoms within geographic regions, units are stationed at designated locations, and once dispatched, a unit becomes unavailable for a probabilistic service time. In this descriptive model, the goal is to evaluate performance measures such as the average response time.


We develop the first optimization model for determining server locations in the hypercube setting. Unlike the descriptive hypercube model, which evaluates performance for fixed locations, our approach optimizes unit placement to minimize mean response time. Prior work provides only heuristic approaches. 
We define \textit{response time} as the sum of \textit{turnout time}, which is the time from call receipt until the response unit begins traveling, and the \textit{travel time} to the incident location. 
In EMS and fire systems, turnout time reflects the delay before a response unit departs the station. In police systems, units are on patrol, so turnout time is essentially zero. The average response time is a key performance measure in emergency service systems. In EMS, shorter response times reduce patient morbidity and increase survival, while in police systems, they increase the likelihood of apprehension.  In fire systems, shorter response time reduce property damage, injury, and death. 

The main contributions of this paper are as follows: 

\begin{itemize} 
    \item[(i).] We provide the first prescriptive model for the server location problem in spatial hypercube queueing models. A key feature is its ability to capture stochastic server availability, heterogeneous utilization levels, and dependencies in availability among servers. Classical approaches, such as the $p$-Median model~\citep{hakimi1964optimum} and covering-based models~\citep{toregas1971location, church1974maximal, daskin1983maximum}, typically assume full unit availability, identical utilization rates, and independence across servers. 

    \item[(ii).] We derive lower and upper bounds for the optimal solution by relating it to a weighted $p$-Median problem and analyze the model’s asymptotic behavior under different traffic regimes. These results serve to benchmark the achievable performance and provide theoretical guarantees on the optimal objective value. In the vanishing-load regime, where demand is light, the objective reduces to that of a weighted $p$-Median problem~\citep{hakimi1964optimum}. In the heavy-load regime, where the offered load grows without bound, the optimal strategy is to locate all servers at the $1$-Median site. 

    \item[(iii).] We develop a parametric model based on a sparse Bayesian linear framework that captures interaction effects between server locations. These interaction effects represent how the placement of one unit influences the effectiveness of others. Building on this parametric model, our method introduces two key innovations: (i) a sparse surrogate model for location decisions (\S\ref{ssec:blm}), and (ii) a submodular relaxation approach that reformulates the binary quadratic program acquisition function into a min-cut problem, enabling efficient optimization~(\S\ref{ssec:acquisition}). We establish a regret bound of $\mathcal{O}(s_0 \sqrt{T \log D})$, where $s_0$ denotes the sparsity level, $T$ the number of iterations, and $D$ the input dimension.
    
    \item[(iv).] We also develop a nonparametric Gaussian process (GP) model that uses the $p$-Median objective as its prior mean (\S\ref{sec_GP}). Conventional GP-based methods~\citep{frazier2018tutorial}, designed for continuous domains, are not well-suited for discrete location problems. To adapt the framework, we introduce several specialized techniques for solving our Bayesian optimization that are described in \S\ref{sec_BO}. We show that this approach also achieves sublinear regret and converges to the optimal solution. 

    \item[(v).] Numerical experiments demonstrate that both our parametric and nonparametric models consistently deliver strong performance across problem sizes, from small-scale to large-scale instances. We apply the methods to an ambulance location problem in St. Paul, Minnesota, using real-world data, and show that both algorithms converge rapidly to the global optimum while outperforming all baseline approaches, including the $p$-Median model, metaheuristic methods, and other Bayesian optimization techniques. 
\end{itemize}


The remainder of this paper is organized as follows. In \S \ref{sec_literature_review}, we review the relevant literature. In \S \ref{sec_problem}, we define the problem and show its relationship to the $p$-Median problem. In \S \ref{sec:HS_linear}, we introduce the parametric sparse Bayesian linear model. In \S \ref{sec_BO}, we develop the nonparametric GP-based model. In \S \ref{sec_numerical}, we investigate the performance of the proposed algorithms through numerical experiments. In \S \ref{sec_realdata}, we apply our solution method to the real data of the St. Paul, Minnesota, emergency service system and discuss managerial insights. Finally, in \S \ref{sec_conclusion}, we present conclusions. 

\section{Literature Review}\label{sec_literature_review}








\subsection{Location Models in Emergency Service Systems}
This section reviews literature on location problems in emergency service systems.
\citet{ahmadi2017survey} provide an extensive overview of healthcare facility location models, while \citet{ingolfsson2013ems} focuses on empirical studies and stochastic modeling for emergency medical services. Location problems are typically classified as either discrete or continuous. Discrete models assign units to predefined sites~\citep{berman2007facility}, while continuous models allow placement anywhere within a region~\citep{baron2008facility}. Our work focuses on the discrete setting. \citet{daskin2008you} offers a detailed survey of such problems. 


Among the most important discrete models are the $p$-Median and covering models. The $p$-Median model \citep{hakimi1964optimum} chooses $p$ facility sites to minimize total distance or travel time to demand points. Covering models, including the location set covering model (LSCM) of \cite{toregas1971location} and the maximal covering location problem (MCLP) of \cite{church1974maximal}, aim to ensure that as many demand locations as possible are covered within a specified response time or distance. These classic models provided a foundation for emergency facility deployment by focusing on geographic coverage and average travel metrics. However, they do not capture the stochastic realities of emergency service systems, where both demand and service times are random, and units are often unavailable due to ongoing calls. 

To better reflect real emergency systems, researchers have extended basic models to account for unit availability. Probabilistic covering models, such as MEXCLP \citep{daskin1983maximum}, estimate expected coverage using fixed unit busy fractions, while MALP \citep{revelle1989maximum} specifies coverage with a stated probability. Although more realistic than deterministic models, these formulations often treat unit busy probabilities as fixed input parameters that do not depend on unit interactions, and assume that unit availabilities are independent of one another. An important contribution by \citet{berman2007facility} incorporates reliability into the 
$p$-Median model by allowing facilities to fail with specified probabilities. Even under the assumption of independent failures, the problem is challenging, with exact solutions only feasible for special cases; heuristics are typically required in general settings.
In contrast, our work considers unit-specific busy probabilities that are not fixed and may be dependent. When one unit is busy, the likelihood of other units being busy increases. This dependence further complicates the problem.

\subsection{Spatial Hypercube Queueing Model}
Researchers have developed models that account for the stochastic nature of emergency service systems. The problem we address is the same as that addressed by the spatial hypercube queueing model introduced by \citet{larson74hypercube, larson75approx}. However, unlike the descriptive hypercube model, which evaluates average system-wide response time for a given set of unit locations, our approach is prescriptive and identifies the optimal set of locations. The hypercube model has been widely cited and applied to police patrol planning~\citep{chaiken1978transfer} and ambulance siting~\citep{brandeau1986extending}, and has also served as input to models such as the maximal expected covering location problem~\citep{batta1989maximal} and network or discrete location formulations~\citep{daskin2011network}. A major computational challenge of the exact hypercube is that it requires solving $2^p$ equations, where $p$ is the number of servers. For larger systems,~\cite{larson75approx} developed an approximation algorithm that requires only solving $p$ linear equations. 

\cite{ghobadi2021integration} surveyed papers that integrate location models with the hypercube framework and identified several works aiming to minimize average response time in the hypercube setting. Unlike our approach, however, these studies are heuristic and do not provide performance guarantees. Examples include \citet{geroliminis2004districting}, using a heuristic method; \citet{geroliminis2009spatial}, a simulation-based approach; \citet{geroliminis2011hybrid}, a heuristic genetic algorithm; \citet{chanta2011minimum}, using a tabu search; and \citet{toro2013joint}, a genetic algorithm. 

The spatial hypercube model has been incorporated into other location models to address stochasticity. For example, \cite{berman1987stochastic} developed two heuristics combining median problems with the hypercube model to solve the stochastic $p$-Median problem, while \cite{ingolfsson2008optimal} proposed an iterative method to determine ambulance locations that minimize fleet size subject to service-level constraints. Unlike these heuristic approaches, we propose a solution method with theoretical performance guarantees, ensuring convergence to the optimal solution with provable regret bounds.

\subsection{Bayesian Optimization}

Bayesian optimization (BO)~\citep{frazier2018bayesian, frazier2018tutorial} is a sample-efficient strategy for optimizing objective functions that are typically expensive to evaluate. The approach has been successfully applied to a wide range of real-world problems, including food safety control~\citep{hu2010contamination} and drug discovery~\citep{negoescu2011knowledge}. 
At its core, BO models the unknown objective function using a probabilistic \textit{surrogate} model that captures both function estimates and uncertainty, and selects the next query point by optimizing an \textit{acquisition function}~(\S\ref{sec_acq}) that balances exploration and exploitation. The process begins by placing a \emph{prior} distribution over the objective function. As function evaluations are collected, Bayes' rule is used to update the belief, yielding a \emph{posterior} distribution.
A common and powerful choice of surrogate is the Gaussian Process (GP), which is defined by a \emph{kernel} function~(\S\ref{sec_GP}), encoding prior assumptions about the function's regularity. Although GPs are widely used, other \emph{parametric} surrogate models can be employed when appropriate. A detailed illustration of BO can be found in~\S\ref{ec:GP_illustration}. 

While many Bayesian optimization techniques have been developed for continuous spaces~\citep{frazier2018tutorial}, applying Bayesian optimization to discrete optimization problems is more challenging because it was originally designed for continuous functions, which are generally easier to optimize due to their smooth and predictable nature~\citep{frazier2018bayesian}. Discrete optimization problems, such as combinatorial optimization problems, are often non-smooth and have a larger number of discrete solutions, which makes the optimization process more difficult. Additionally, Gaussian processes (discussed in \S \ref{sec_GP}), which are commonly used in Bayesian optimization as surrogate models, are not well-suited for discrete spaces, making it difficult to develop accurate models~\citep{williams2006gaussian}. Therefore, special methods need to be developed to handle optimization in discrete spaces, as we do in this paper. 


Several recent studies have extended Bayesian optimization (BO) to discrete or categorical domains. For example, \citet{baptista2018bayesian} proposed a second-order monomial representation for the objective function. \citet{ru2020mix} combined ideas from multi-armed bandits with GP-based BO improving exploration in combinatorial settings, while \citet{wan2021think} introduced trust regions to better navigate structured search spaces. Other lines of work include diffusion kernels for discrete objects \citep{deshwal2020mercer}, variational  optimization methods \citep{wu2020variationalCombinatorial}, and GP surrogates defined over graphs \citep{oh2019combinatorial}. 
Despite these advances, existing methods still face three challenges: (i) lack of theoretical convergence guarantees, (ii) limited scalability in high-dimensional discrete spaces, and (iii) limited applicability to location problems. Our work addresses these gaps by providing regret guarantees, scalable algorithms, and a modeling approach for location models in spatial queues.

\section{Model}\label{sec_problem}
This section formulates the \emph{approximate} spatial hypercube location model. 
Section~\ref{ssec:problem} introduces the problem setting and the $p$-Hypercube model. 
Section~\ref{sec:bounds} establishes lower and upper bounds by relating our model to the deterministic $p$-Median problem, linking stochastic queueing-based models with classical location theory. 
Section~\ref{ssec:asymptotic} analyzes the model’s asymptotic behavior, showing that the optimal configuration converges to the $p$-Median and $1$-Median solutions under vanishing-load and heavy-load regimes, respectively. 

\subsection{Problem Formulation}\label{ssec:problem}
We select $p$ locations from $N$ candidate sites, denoted by $\mathcal{I}=\{1,\cdots,N\}$. The service region is partitioned into $M$ subregions, denoted by $\mathcal{J}=\{1,\cdots,M\}$, with calls in each subregion $j \in \mathcal{J}$ arriving as a Poisson process with rate $\lambda_j$. The service time for a unit at location $i$ serving subregion $j$ follows a general distribution with mean $1/\mu_{ij}$, covering the time from dispatch to unit return. The response time consists of the turnout time $\tau_i$ and the travel time $t_{ij}$ from location $i$ to subregion $j$. 

When a call arrives, it is assigned to the highest-ranked available unit according to a predetermined preference list $(\gamma_{jl})$, where $\gamma_{jl}$ denotes the index of the $l$-th preferred unit for subregion $j$ (e.g., ordered by proximity or travel time). For simplicity, we assume that if no units are available, the call is lost, which is appropriate when the probability that all units are simultaneously busy is negligible. In more congested settings, our model can incorporate queues~\citep{larson75approx}; all theoretical results here remain valid. 

We define binary variables $x_i$, where $x_i = 1$ indicates a unit at location $i$ and $x_i = 0$ indicates no unit. The vector $\boldsymbol{x} = \{x_1,\cdots, x_{N}\} \in  \mathcal{D}$ represents unit locations, where $\mathcal{D}= \{0,1\}^N$. Our objective is to determine the optimal unit locations $\boldsymbol{x}^*$ that minimize the system-wide mean response time. To achieve this, we solve an integer program that we call \emph{$p$-Hypercube}.
\begin{subequations}\label{eq_all_approx}
\setlength{\abovedisplayskip}{3pt}
\setlength{\belowdisplayskip}{3pt}
\begin{align}
\text{OPT}_{H} = \min_{\boldsymbol{x}} & \sum_{i \in \mathcal{I}} \sum_{j \in \mathcal{J}}(\tau_i + t_{i j}) \frac{\lambda_j}{\sum_{j\prime\in \mathcal{J}} \lambda_{j\prime}}\frac{q_{i j}(\boldsymbol{x})}{\sum_{i' \in \mathcal{I}} q_{i' j}(\boldsymbol{x})} \qquad \text{\emph{(\textbf{$p$-Hypercube})}} \label{eq_obj}\\
\text { s.t. } & \sum_{i \in \mathcal{I}} x_{i}=p, \label{eq_feasibility_approx}\\
& q_{i j}(\boldsymbol{x})=  Q(N,\bar{\rho}, \eta_{ij}-1)\left(\prod_{l=1}^{\eta_{ij}-1} \rho_{\gamma_{j l}}(\boldsymbol{x})\right)\left(1-\rho_{i}(\boldsymbol{x})\right),  \quad \forall i \in \mathcal{I}, \forall j \in \mathcal{J}, \label{eq_qij_approx}\\
& x_{i} \in\{0,1\}, \qquad \forall i \in \mathcal{I}.
\end{align}
\end{subequations}

In this formulation, $q_{ij}(\boldsymbol{x})$ represents the proportion of calls from subregion $j$ that are handled by the unit located at $i$, given the unit locations $\boldsymbol{x}$, and $\frac{q_{i j}(\boldsymbol{x})}{\sum_{i' \in \mathcal{I}} q_{i' j}(\boldsymbol{x})}$ represents the normalized fraction of demand from subregion $j$ that is served by unit $i$.
The term $\rho_i(\boldsymbol{x})$ denotes the utilization of unit $i$ under configuration $\boldsymbol{x}$, $Q$ is a correction factor that accounts for dependencies among unit availabilities, and $\bar{\rho}$ is the average server utilization. Finally, $\eta_{ij}$ specifies the rank of unit $i$ in the preference list for subregion $j$.

Equation~\eqref{eq_qij_approx} uses the approximate hypercube model from \citet{larson75approx}, since computing the exact values of $q_{ij}(\boldsymbol{x})$ would require solving the full spatial hypercube model of \citet{larson74hypercube}, which becomes computationally infeasible for large-scale systems. 
This approximation method is derived by treating server workloads as being independent and then adjusting for this obvious error through a correction factor $Q$, which is a function of the number of servers, traffic intensity, and server preference. We also provide the exact formulation in \S~\ref{ssec:exact_Hypercube}, where the lower and upper bounds, as well as the asymptotic results established in this section, continue to hold.

Define $P(k)$ as the probability that exactly $k$ servers are busy. The correction factor is given by
\begin{equation}
    Q(N,\bar{\rho},r) = \sum_{k=r}^{N-1} \frac{{k \choose r}}{{N \choose r}}\frac{N-k}{N-r}\frac{P(k)}{\bar{\rho}^{r}(1-\bar{\rho})}, \label{eq:Q_factor}
\end{equation}
where probabilities $P(k)$ are obtained by the Erlang loss formula~\citep{erlang1917loss}. The unit utilizations $\rho_i$ are obtained by solving the following $p$ non-linear equations,
\begin{equation}
    \rho_i = 1 - \left( 1 + \sum_{k=1}^N \sum_{j \in G_{ik}} \frac{\lambda_j Q(N, \bar{\rho}, k-1) \prod_{l=1}^{k-1} \rho_{\gamma_{jl}}}{\mu_{ij}} \right)^{-1},\label{Eq:rhoi}
\end{equation}
where $G_{ik}$ denotes the set of nodes for which unit $i$ is the $k$-th preferred option. Equivalently, $G_{ik}$ can be viewed as the inverse mapping of the preference list $\gamma_{j k}$. The average server utilization $\bar{\rho} = \sum_{j=1}^{M}\frac{\lambda_j}{N}\sum_{i=1}^{N}\frac{q_{ij}(\boldsymbol{x})}{\mu_{ij}}$ depends only on the mean service rates $\mu_{ij}$ and not its distribution \citep{sevastyanov_ergodic_1957}. This utilization measure is then used to update the correction factors in the next iteration. In contrast, the exact hypercube formulation requires solving $2^p$ detailed balance equations to determine the steady-state distribution across all possible server states, making it computationally intractable for larger systems. 


\begin{theorem}
\label{thm:NP-hard_approx}
The $p$-Hypercube problem is NP-hard.
\end{theorem}

We prove that the $p$-Hypercube problem is NP-hard by a polynomial reduction of the $p$-Median problem, which is known to be NP-hard \citep{kariv1979algorithmic}, as detailed in~\S \ref{ssec:proof_NPhard_approx}.

\subsection{Lower and Upper Bounds}
\label{sec:bounds}

This section establishes lower and upper bounds for the optimal solution of the $p$-Hypercube problem by establishing its connection to the deterministic $p$-Median problem. 
The $p$-Median problem \citep{hakimi1964optimum} determines the minimum average response time to all calls while assuming units are always available.

We first derive a lower bound. In the deterministic $p$-Median problem, binary variables $x_i\in \{0,1\}$ indicate whether a unit is located at position $i$, similar to the variables in $p$-Hypercube. Additionally, we introduce binary variables $y_{ij}\in \{0,1\}$ to indicate the assignment of calls from subregion $j$ to the unit located at location $i$. We denote the optimal solution to the $p$-Median problem as $(\boldsymbol{\hat{x}}^{*}, \boldsymbol{\hat{y}}^{*})$. The objective of the $p$-Median model is to minimize the total call-weighted response time, where $w_j$ is the weight for subregion $j$. We have
\begin{subequations}\label{eq_all_pmedian}
\setlength{\abovedisplayskip}{3pt}
\setlength{\belowdisplayskip}{3pt}
\begin{align}
\text{OPT}_{M} &= \min_{\boldsymbol{x}, \boldsymbol{y}} \sum_{i \in \mathcal{I}} \sum_{j \in \mathcal{J}} w_j (\tau_i + t_{i j}) y_{i j} \qquad \text{\emph{(\textbf{$p$-Median})}} \label{eq:objective} \\
\text{s.t.} \quad &\sum_{i \in \mathcal{I}} x_{i}=p, \label{eq:x_sum} \\
               &\sum_{i \in \mathcal{I}} y_{i j}=1, \quad \forall j \in \mathcal{J}, \label{eq:y_sum} \\
               &y_{i j} \leq x_{i},  \quad \forall i \in \mathcal{I}, \forall j \in \mathcal{J}, \label{eq:y_valid} \\
               &y_{i j} \in\{0,1\},  \quad\forall i \in \mathcal{I}, \forall j \in \mathcal{J}, \label{eq:y_binary} \\
               &x_{i} \in\{0,1\},  \quad \forall i \in \mathcal{I}.
\end{align}
\end{subequations}

Constraint \eqref{eq:y_sum} requires that each subregion is assigned to a unit, and Constraint \eqref{eq:y_valid} requires that calls in subregion $j$ are only assigned to an occupied location. 
We next show that a special case of the $p$-Median problem provides a lower bound for the optimal value of $p$-Hypercube. The proof is detailed in~\S\ref{ssec_lower_bound_proof_approx}.

\begin{theorem}
\label{thm_low_bound_approx}
The optimal value to the $p$-Median problem is a lower bound for the optimal value of the  $p$-Hypercube problem when $w_j = \lambda_j/\sum_{j'\in \mathcal{J}} \lambda_{j'}$, i.e.,
\begin{equation*}
\text{OPT}_{H} \;\geq\; \text{OPT}_M.
\end{equation*}

\end{theorem}

Next, we derive an upper bound, as presented in Lemma \ref{lemma_upper_bound}.
\begin{lemma}
\label{lemma_upper_bound}
By setting $w_j = \lambda_j/\sum_{j'\in \mathcal{J}} \lambda_{j'}$, applying the optimal solution $\boldsymbol{\hat{x}}^{*}$ obtained from $p$-Median to $p$-Hypercube provides an upper bound on the optimal value of $p$-Hypercube.
\end{lemma}
Given that the optimal solution $\boldsymbol{\hat{x}}^{*}$ is always a feasible solution for the $p$-Hypercube problem, we leverage this to establish an upper bound on the optimal value of the $p$-Hypercube problem, which is a minimization problem.

\subsection{Asymptotic Results}
\label{ssec:asymptotic}
We next analyze the asymptotic behavior of the proposed model across different traffic regimes. By varying the offered load, we establish connections to deterministic location problems and derive structural insights into the limiting behavior of optimal solutions. 

Formally, define the scaled arrival rates $
\lambda_j(\theta)=\theta\lambda_j,
$ 
where $\theta>0$ is a scaling factor that adjusts the demand intensity in each subregion. The corresponding objective function at scale $\theta$ is 
\[
\text{OPT}_{H}(\theta,\boldsymbol{x})\;:=\;\sum_{i\in\mathcal I}\sum_{j\in\mathcal J}(\tau_i+t_{ij})\frac{\lambda_j}{\sum_{j\prime\in \mathcal{J}} \lambda_{j\prime}}\frac{q_{i j}(\theta;\boldsymbol{x})}{\sum_{i' \in \mathcal{I}}  q_{i' j}(\theta;\boldsymbol{x})},
\]
where $q_{ij}(\theta; \boldsymbol{x})$ is the conditional share among calls from $j$ served by unit $i$ under configuration $\boldsymbol{x}$ when the arrival rate is scaled by $\theta$. We identify two canonical limiting regimes: the \emph{vanishing-load regime}, where demand intensity is small, and the \emph{heavy-load regime}, where demand becomes arbitrarily large. Together, these regimes characterize the extremes of the system and link the optimal $p$-Hypercube solution to classical location problems.

\begin{theorem}[Vanishing load]\label{thm:asym-vanishing-lambda}
There exists $\theta_0>0$ such that for all $0<\theta<\theta_0$, any optimal solution $\boldsymbol{x}^*(\theta)\in \arg\min_{\boldsymbol{x}\in\mathcal D} \text{OPT}_{H}(\theta,\boldsymbol{x})$ to $p$-Hypercube is also an optimal solution to the weighted $p$-Median.
\end{theorem}

In the vanishing-load regime, the probability that multiple units are simultaneously busy becomes negligible. As a result, unit availability constraints effectively disappear, and the problem reduces to minimizing the weighted travel time (including turnout times). This is precisely the weighted $p$-Median objective~\citep{hakimi1964optimum}. Thus, the $p$-Hypercube model generalizes the $p$-Median in light-traffic systems, and the two coincide asymptotically. The formal proof is provided in \S\ref{ssec:asym-vanishing-lambda}.

\begin{theorem}[Heavy load]\label{thm:asym-heavy-lambda}
If colocating multiple servers at the same site is allowed, then as $\theta\to\infty$, the optimal solution to the $p$-Hypercube problem is to place all $p$ units at the $1$-Median location. 
\end{theorem}

In the heavy-load regime, units are almost always busy, and the system effectively operates under maximum congestion. The bottleneck lies not in balancing coverage across subregions but in ensuring that service capacity is concentrated at the single most central location. Hence, colocating all units at the $1$-Median minimizes the average response time when a unit does become available. The proof is provided in \S\ref{ssec:asym-heavy-lambda}.


\section{Parametric Surrogate: Sparse Bayesian Linear Model}\label{sec:HS_linear}

This section develops a parametric surrogate model using a sparse Bayesian linear formulation to address the $p$-Hypercube problem. The following subsection provides an overview of the modeling framework and its key components. 

\subsection{An Overview of the Parametric Approach}

The parametric approach approximates the objective function using a sparse Bayesian linear model. A feasible solution is a binary vector $\boldsymbol{x}\in\{0,1\}^N$ with $\|\boldsymbol{x}\|_0=p$, and the goal is to find $\boldsymbol{x}^*$ that minimizes the mean response time $f(\boldsymbol{x}^*)$.

We begin by randomly selecting an initial set of feasible configurations $\boldsymbol{x}_1, \ldots, \boldsymbol{x}_\mathcal{T}$ and evaluating their objective values $f(\boldsymbol{x}_1),\ldots, f(\boldsymbol{x}_\mathcal{T})$ using the $p$-Hypercube objective~\eqref{eq_obj}. These observations are then used to fit a Bayesian linear model with a hierarchical horseshoe prior (\S\ref{ssec:blm}). At each iteration, we sample $\alpha_t$ from the posterior, defining a new realization of the surrogate model $f(\boldsymbol{x})$. We choose $\boldsymbol{x}_{t+1}$ by minimizing the sampled surrogate via a binary quadratic program, solved efficiently using a submodular relaxation that reformulates the problem as a minimum-cut (\S\ref{ssec:acquisition}). The new evaluation updates the posterior, refining parameter estimates and uncertainty quantification. This iterative process continues until the evaluation budget $T$ is reached. Theoretical analysis (\S\ref{ssec:theory}) shows the algorithm achieves sublinear cumulative regret.

\subsection{Bayesian Linear Model with Horseshoe Prior}\label{ssec:blm}

The $p$-Hypercube location problem is NP-hard, making exact optimization infeasible. To address this, we develop surrogates that approximate the hypercube formulation while building upon both the hypercube framework and deterministic location models, thereby preserving the key stochastic characteristics of emergency service systems and enabling tractable optimization.

We introduce a parametric surrogate model inspired by the double coverage concept \citep{hogan1986Concepts, pirkul1988Siting} in emergency service systems. The idea of double coverage is that each demand region should be covered by at least two nearby response units, so that when one unit is busy, another is still available to respond. This enhances reliability and increases the chance of prompt service. Our model builds on this principle by explicitly considering the interaction between pairs of units when determining locations.


\begin{example}
    If a subregion is primarily covered by two nearby units, each with busy probability $b$ (assuming independence), the probability that at least one is available is $1 - b^2$. With $b = 0.3$, which is a relatively high utilization level in emergency services, the probability that both are busy drops to $b^2 = 9\%$, giving a 91\% likelihood that a unit is free to respond. Thus, a third unit is rarely needed, consistent with~\citet{iannoni2007multiple} for Brazilian-highway EMS dispatch prioritizing the nearest and second-nearest units.
\end{example}


Motivated by this, we propose a linear surrogate model that incorporates both first-order and second-order interaction terms to capture the effects of individual unit placements and the pairwise interactions between units. The model is
{
\setlength{\abovedisplayskip}{3pt}
\setlength{\belowdisplayskip}{3pt}
\begin{equation}\label{eqn:bayesian-linear-surrogate}
f_\alpha(\boldsymbol{x}) = \alpha_0 + \sum_j \alpha_j x_j + \sum_{i,j>i} \alpha_{ij} x_i x_j.
\end{equation}
}


Here, the first-order terms ($\alpha_j x_j$) represent the marginal benefit of placing a unit at site $j$, while the second-order terms ($\alpha_{ij} x_i x_j$) capture how the placement of one unit influences the contribution of another, such as when their service areas overlap or when they jointly cover the same high-demand region. These models are also known as regression metamodels in simulation. The theoretical foundation is that any smooth response function can be locally approximated by a low-order Taylor expansion \citep{box1992experimental}. Thus second-order terms capture how stations influence each other’s performance. When stations are located near one another, their coverage areas overlap and additional units provide limited improvement to overall performance. In contrast, when stations are placed in distinct regions with minimal overlap, they collectively enhance system coverage by serving nearby but separate demand areas. 




However, model~\eqref{eqn:bayesian-linear-surrogate} is linear in $\boldsymbol{\alpha} = (\alpha_i, \alpha_{ij}) \in \mathbb{R}^D$, with $D = 1 + N + \binom{N}{2}$ and is dense, including a large number of interaction terms that are often insignificant in practice. In particular, interactions between certain distant or weakly related location pairs contribute little to the objective. To encourage sparsity and focus on influential parameters, we adopt a Bayesian framework and, following \citet{baptista2018bayesian}, impose a horseshoe prior on $\boldsymbol{\alpha}$.
The horseshoe prior is a continuous shrinkage prior that pulls small, uninformative coefficients strongly toward zero while allowing a few important ones to remain large. Under this prior, we have
\begin{equation}\label{eqn:horseshoe-prior}
    \begin{aligned}
\boldsymbol{y} \mid \boldsymbol{X}, \boldsymbol{\alpha}, \sigma &\sim \mathcal{N}(\boldsymbol{X}\boldsymbol{\alpha}, \sigma^2 \boldsymbol{I}_T)\\
\alpha_k \mid \beta_k, \tau, \sigma &\sim \mathcal{N}(0, \beta_k^2 \tau^2 \sigma^2) \quad k = 1, \ldots, D \\
\tau, \beta_k &\sim C^+(0,1) \quad k = 1, \ldots, D \\
\sigma^2 &\sim \sigma^{-2} \text{d} \sigma^{2}.
\end{aligned}
\end{equation}
Specifically, each coefficient $\alpha_k$ is normally distributed with variance scaled by a global parameter $\tau^2$ and a local parameter $\beta_k^2$, both of which follow independent half-Cauchy distributions denoted by $C^+(0,1)$ with probability density function
$p(z) = \frac{2}{\pi(1 + z^2)}$ for $z > 0$. Note that the index 
$k$ runs over all elements of $\boldsymbol{\alpha}$, including both single-location terms ($\alpha_i$) and interaction terms ($\alpha_{ij}$). We are using an improper prior\footnote{An improper prior is a prior distribution that does not integrate to one over the parameter space and therefore is not a valid probability distribution on its own.  When applied to the noise variance $\sigma^2$, an improper prior reflects the absence of strong prior beliefs and allows the data to fully determine the inference about the variance. It can be used in Bayesian inference provided that the resulting posterior distribution is proper, see e.g. \citet{gelman1995bayesian}.} for $\sigma^2$ with density proportional to $1/\sigma^2$ over $(0, \infty)$ and $\text{d}\sigma^2$ refers to the Lebesgue measure over the continuous domain of $\sigma^2$. 

The horseshoe prior induces sparsity by shrinking unimportant $\alpha_k$ toward zero while allowing significant ones to retain their magnitude~\citep{carvalho2010horseshoe}. Unlike traditional shrinkage priors such as lasso~\citep{tibshirani1996regression}, it preserves large coefficients rather than over-penalizing them, enabling both variable selection and accurate estimation. 
The global shrinkage parameter $\tau$ determines the overall level of sparsity by controlling how strongly all coefficients are pulled toward zero. Each local shrinkage parameter $\beta_k$ adjusts shrinkage on $\alpha_k$, allowing important parameters to remain large while suppressing less relevant ones. 
This hierarchical prior structure is particularly effective in high-dimensional settings, where most pairwise interaction terms, which describe how the presence of one unit influences the performance of another (for example, through overlapping coverage or shared demand regions), are expected to have minimal impact. 
Details on deriving a tractable posterior distribution using the reparameterization method are provided in~\S\ref{ec:hs_posteiror}. 


\subsection{Acquisition Function}\label{ssec:acquisition}

The acquisition function guides the search by selecting the point to evaluate. We adopt Thompson sampling~\citep{thompson1935TS}. 
At each iteration $t$, we sample $\boldsymbol{\alpha}_t \sim \mathcal{P}(\boldsymbol{\alpha} \mid \boldsymbol{X}, \boldsymbol{y})$ 
from the current posterior. This sampled parameter defines a specific surrogate realization $f_{\boldsymbol{\alpha}_t}(\boldsymbol{x})$. Using this sampled model, we then identify the next configuration $\boldsymbol{x}_{t+1}$ that minimizes  $f_{\boldsymbol{\alpha}_t}(\boldsymbol{x})$.
\begin{equation}
\argmin_{\boldsymbol{x} \in \mathcal{D},\|\boldsymbol{x}\|_0=p} f_{\boldsymbol{\alpha}_t}(\boldsymbol{x})= \argmin_{\boldsymbol{x} \in \mathcal{D},\|\boldsymbol{x}\|_0=p} \sum_j \alpha_{jt} x_j + \sum_{i<j} \alpha_{ijt} x_i x_j,
\end{equation}
The above minimization problem can be written as a binary quadratic programming (BQP) problem~\citep{Kochenberger2014UBQP} 
\begin{equation}\label{eqn:BQP}
\arg\min_{\boldsymbol{x} \in \mathcal{D},\|\boldsymbol{x}\|_0=p} \boldsymbol{x}^T \boldsymbol{A} \boldsymbol{x} + \boldsymbol{b}^T \boldsymbol{x},
\end{equation}
where the matrix $\boldsymbol{A} \in \mathbb{R}^{N \times N}$ captures the interaction coefficients $\left(\alpha_{i j}\right)$, and the vector $\boldsymbol{b} \in \mathbb{R}^N$ contains the linear terms $\left(\alpha_j\right)$. 

One could solve ~\eqref{eqn:BQP} via a semidefinite programming (SDP) relaxation~\citep{VandenbergheBoyd1996}, followed by randomized rounding to obtain a binary solution in $\mathcal{D}$. However, this approach is unsuitable for our problem due to two key limitations: (i) rounding may produce solutions far from optimal, as integer solutions cannot be reliably obtained through LP rounding; (ii) our problem requires a strict cardinality constraint ($\|\boldsymbol{x}\|_0 = p$), which the SDP relaxation cannot enforce.


\paragraph{Submodular Relaxation.}
We propose a fast and scalable approach for solving this BQP problem with a cardinality constraint, building upon recent advances in submodular relaxation~\citep{ito2016Sub}. The objective function~\eqref{eqn:BQP} is \textit{submodular} when $A_{ij} \leq 0$ for all $i,j$. Such functions admit exact minimization via graph-cut algorithms~\citep{fujishige2005submodular}. However, for general cases where the objective may not be submodular, we employ a relaxation technique developed by~\cite{ito2016Sub}. This approach constructs a submodular lower bound by decomposing the matrix $A$ into its positive and negative components:
\begin{equation*}
    \boldsymbol{A} = \boldsymbol{A}^+ + \boldsymbol{A}^-, \quad A^+_{ij} = \max(A_{ij}, 0), \quad A^-_{ij} = \min(A_{ij}, 0),
\end{equation*}
and introduces a matrix of relaxation parameters $\Gamma$ to define the following submodular lower bound:
\begin{equation} \label{eq:sub_lower}
\boldsymbol{x}^T (\boldsymbol{A}^+ \circ \Gamma) \boldsymbol{1} + \boldsymbol{1}^T (\boldsymbol{A}^+ \circ \Gamma) \boldsymbol{x} - \boldsymbol{1}^T (\boldsymbol{A}^+ \circ \Gamma) \boldsymbol{1} \leq \boldsymbol{x}^T \boldsymbol{A}^+ \boldsymbol{x}, 
\end{equation}
where $\circ$ denotes element-wise multiplication and $\boldsymbol{1}$ is the all-ones vector. The matrix $\Gamma$ controls the trade-off between relaxation tightness and computational tractability. The lower bound \eqref{eq:sub_lower} is linear in $\boldsymbol{x}$ for the $\boldsymbol{A}^+$ part, ensuring submodularity, and the full relaxed submodular function can be written as 
\begin{equation}
    f_{\text{sub}}(\boldsymbol{x}) = \boldsymbol{x}^T \boldsymbol{A}^- \boldsymbol{x} + \tilde{\boldsymbol{b}}^T \boldsymbol{x} - \boldsymbol{1}^T (\boldsymbol{A}^+ \circ \Gamma) \boldsymbol{1}, 
\end{equation}
where $\tilde{\boldsymbol{b}}^T \boldsymbol{x} = \boldsymbol{b}^T \boldsymbol{x} + \boldsymbol{x}^T (\boldsymbol{A}^+ \circ \Gamma) \boldsymbol{1} + \boldsymbol{1}^T (\boldsymbol{A}^+ \circ \Gamma) \boldsymbol{x}$. 

To solve the relaxed problem, we build a directed graph whose edge capacities represent the submodular objective. Each interaction term $A_{i j}^{-}$ is represented by an edge between nodes $i$ and $j$ with capacity $\left|A_{i j}^{-}\right|$, while the adjusted linear term $\tilde{b}_i$ determines the capacities of edges connecting node $i$ to the source or target. The cut separates nodes into source-side (corresponding to $x_i = 0$) and target-side (corresponding to $x_i = 1$). The following example illustrates the directed-graph construction.

\begin{figure}[btph]
\vspace{-0.5cm}
\centering
  \includegraphics[width=0.5\textwidth]{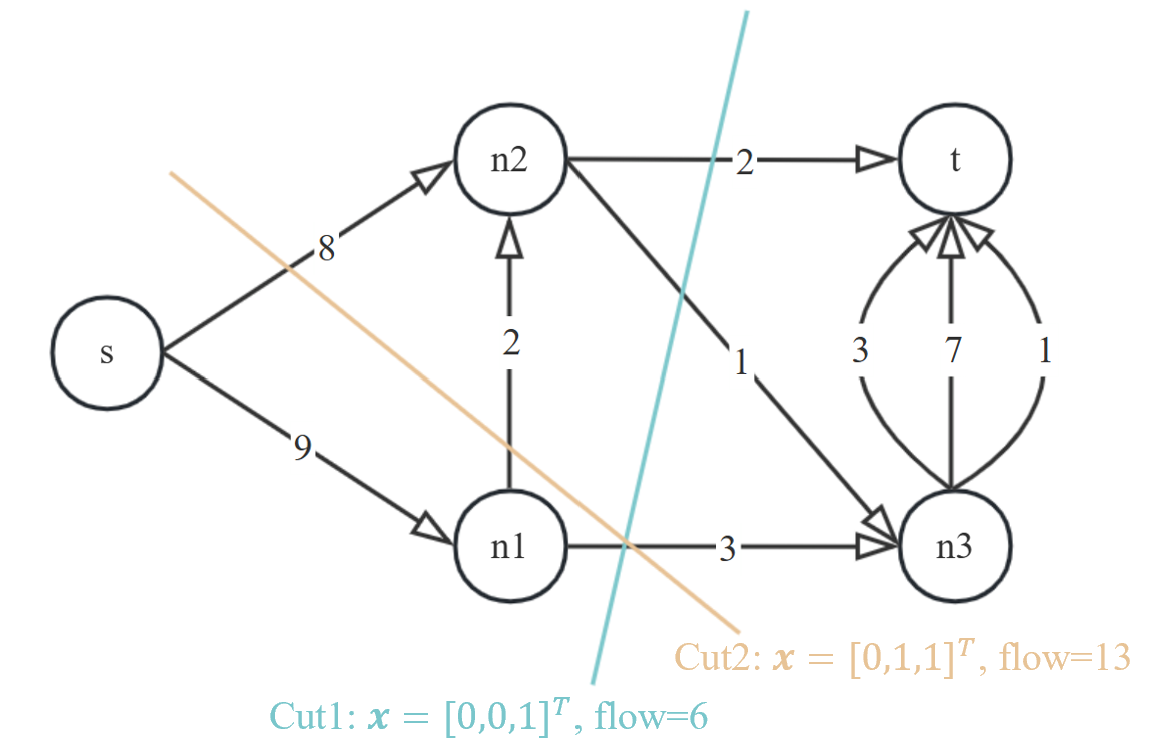}
\caption{An example of graph construction.}
\raggedright
\label{fig:graph}
\vspace{-0.2cm}
\end{figure}

\begin{example}
    Consider the quadratic function 
$f(\boldsymbol{x}) = \boldsymbol{x}^{T}
\begin{bmatrix}
0 & -2 & -3 \\
-2 & 0 & -1 \\
-3 & -1 & 0
\end{bmatrix}
\boldsymbol{x} +
\begin{bmatrix}
9 \\
8 \\
-7
\end{bmatrix}^{T} \boldsymbol{x}$,
which attains $\min f(\boldsymbol{x}) = -7$ at $\boldsymbol{x} = [0,0,1]^T$. We transform this quadratic function into the directed graph shown in Figure \ref{fig:graph}. The graph has a source node $s$, a sink node $t$, and intermediate nodes $ n_1, n_2, n_3$ corresponding to the decision variables. Edge weights encode the quadratic and linear terms of the function. The figure shows two possible cuts. Among all feasible cuts, Cut 1 is the minimum-cut and matches the true minimizer of the quadratic function.
\end{example}

\begin{algorithm}[h]
\caption{Subroutine: Solving BQP with Submodular Relaxation} 
\label{algo:submodular_relaxation}

\begin{algorithmic}[1]
\State \textbf{Input:} $\boldsymbol{A}, \boldsymbol{b}, p, \epsilon$
\State Decompose $\boldsymbol{A} = \boldsymbol{A}^+ + \boldsymbol{A}^-$
\State Initialize $\Gamma \leftarrow \boldsymbol{1}\boldsymbol{1}^T$
\Repeat
   \State Construct submodular lower bound:
   \vspace{-10pt}
   \begin{equation*}
       f_{\text{sub}}(\boldsymbol{x})
       :=\boldsymbol{x}^T\boldsymbol{A}^-\boldsymbol{x}+\tilde{\boldsymbol{b}}^T\boldsymbol{x} - \boldsymbol{1}^T(\boldsymbol{A}^+ \circ \Gamma)\boldsymbol{1}
   \end{equation*}
   \vspace{-20pt}
   \State Build graph $G$ with edge capacities:
   \State \quad - Interaction edges: $c_{ij} = |A^-_{ij}|$
   \State \quad - Source/sink edges: $c_{\mathcal{S}i} = \tilde{b_i}^+$, $c_{i\mathcal{T}} = |\tilde{b}_i^-|$
   
   \State Solve min-cut on $G$ → $\boldsymbol{x}^*$
   \State Update $\Gamma$ via projected gradient descent:
   \vspace{-10pt}
   \begin{equation*}
       \Gamma \leftarrow \text{Proj}_{[0,1]}\left(\Gamma - \eta \nabla_\Gamma \mathcal{L}(\Gamma)\right)
   \end{equation*}
   \vspace{-20pt}
\Until{$\|\Gamma^{(k)} - \Gamma^{(k-1)}\|_F < \epsilon$}
\State \textbf{Output:} $\boldsymbol{x}^*$
\end{algorithmic}
\end{algorithm}

\paragraph{Cardinality Constraint.}
In addition to the $s$–$t$ cut, we impose a cardinality constraint $\|\boldsymbol{x}\|_0 = p$ to ensure exactly $p$ units are selected. We adopt the LP-based polynomial-time algorithm of \citet{Chen2016MincutSize} and update $\Gamma$ using projected gradient descent~\citep{rosen1960gradient}:
\begin{equation*}
    \Gamma \leftarrow \operatorname{Proj}_{[0,1]}\left(\Gamma-\eta \nabla_{\Gamma} \mathcal{L}(\Gamma)\right),
\end{equation*}
where the loss function $\mathcal{L}(\Gamma) = \|\boldsymbol{x^*}^{T} \boldsymbol{A} \boldsymbol{x^*}+\boldsymbol{b}^\top\boldsymbol{x}^* - f_{\text{sub}}(\boldsymbol{x^*})\|^2$ measures the deviation between the true objective value and its submodular approximation. This iterative process continues until convergence, i.e., until the change in $\Gamma$ is below a threshold $\epsilon$. The final solution $\boldsymbol{x^*}$ is then returned as the solution to the original BQP problem.



Algorithm~\ref{algo:submodular_relaxation} solves the BQP via submodular relaxation; Algorithm~\ref{algo:BO_with_submodular} gives the full procedure.
At each iteration, the algorithm samples model parameters, builds a quadratic acquisition function, and approximately solves a binary quadratic program under a cardinality constraint to select the next point. It then evaluates the true function at this point, updates the dataset, and ultimately outputs the best solution obtained.

\begin{algorithm}[h]
\caption{Sparse Bayesian Linear Model (SparBL)}
\label{algo:BO_with_submodular}
{\fontsize{10}{16}\selectfont
\begin{algorithmic}[1]

\State \textbf{Input:} Search space $\mathcal{D} = \{0,1\}^N$, cardinality $p$, initial points $\boldsymbol{X}_0$, observations $\boldsymbol{y}_0$, iterations $T$
\State Initialize dataset $\mathcal{D}_0 \leftarrow (\boldsymbol{X}_0, \boldsymbol{y}_0)$

\For{$t = 1$ \textbf{to} $T$}
    \State \textbf{Surrogate Model:} Fit SparBL with Horseshoe Prior
    \State \quad Sample from posterior via Gibbs sampling (Eq.~\ref{eqn:horseshoe-posterior}):
    \vspace{-10pt}
    \[
    \boldsymbol{\alpha}_t \sim p(\boldsymbol{\alpha} \mid \mathcal{D}_{t-1}), \quad \sigma^2_t \sim p(\sigma^2 \mid \mathcal{D}_{t-1})
    \]
    \vspace{-20pt}
    \State \textbf{Acquisition:} Solve BQP with Submodular Relaxation
    \State \quad Construct quadratic form from $\boldsymbol{\alpha}_t$:
    \vspace{-10pt}
    \[
    f_{\boldsymbol{\alpha}_t}(\boldsymbol{x}) = \boldsymbol{x}^T \boldsymbol{A}_t \boldsymbol{x} + \boldsymbol{b}_t^T \boldsymbol{x}
    \]
    \vspace{-20pt}
    \State \quad Run Algorithm~\ref{algo:submodular_relaxation} with matrix $\boldsymbol{A}_t$, vector $\boldsymbol{b}_t$, cardinality constraint $p$
    \State \quad Obtain next point: $\boldsymbol{x}_t \leftarrow \argmin_{\|\boldsymbol{x}\|_0=p} f_{\boldsymbol{\alpha}_t}(\boldsymbol{x})$
    
    \State \textbf{Evaluation:} Observe $y_t \leftarrow f(\boldsymbol{x}_t) + \epsilon_t$
    \State \textbf{Update:} $\mathcal{D}_t \leftarrow \mathcal{D}_{t-1} \cup \{(\boldsymbol{x}_t, y_t)\}$
\EndFor

\State \textbf{Output:} Best solution $\boldsymbol{x}^* \leftarrow \argmin_{\boldsymbol{x} \in \{\boldsymbol{x}_1,...,\boldsymbol{x}_T\}} f(\boldsymbol{x})$
\end{algorithmic}}
\end{algorithm}

\subsection{Theoretical Results}\label{ssec:theory}
To establish the theoretical results, we first outline the technical assumptions on which our analysis is based. 

\begin{assumption}[Boundedness Assumption]\label{assump:true-param}
There exist positive constants $s_0 \in \mathbb{N}$ and $W \in \mathbb{R}^+$ such that $\|\boldsymbol{\alpha}\|_0 = s_0$ and $\|\boldsymbol{\alpha}\|_1 \leq W$.
\end{assumption}
Following a standard assumption commonly adopted in bandit and high-dimensional regression studies, we require boundedness of the true parameter vector $\boldsymbol{\alpha}$ to ensure that the regret bound remains scale-free~\citep{abbasi-yadkori2011improved, bastani2020online}. In addition, sparsity is promoted in our model through the horseshoe prior. 
The next assumption is about the design matrix. We use $\boldsymbol{X}_t$ to denote the matrix $(\boldsymbol{x}^{(1)}, \ldots, \boldsymbol{x}^{(t)})^\top \in \mathbb{R}^{t \times D}$. The corresponding empirical covariance matrix is defined as $\hat{\Sigma}_t = \boldsymbol{X}_t^\top \boldsymbol{X}_t / t$. 

\begin{assumption}[Compatibility Assumption]\label{assump:compatibility}

Writing $\Phi_t=[\phi(x_1),\ldots,\phi(x_t)]^\top$,
there exist $t_0\in\mathbb{N}$ and $\kappa>0$ such that, with probability at least $1-2e^{-c t}$ for some $c>0$ and all $t\ge t_0$,
\begin{equation}
\frac{\|\Phi_t v\|_2}{\sqrt{t}} \ \ge\ \kappa\,\|v\|_2
\quad\text{for all }v\in\mathcal{C}(s_0):=\{v:\ \|v_{(s_0)^c}\|_1\le 3\|v_{s_0}\|_1\}.
\label{eq:compatibility}
\end{equation}


\end{assumption}

Assumption~\ref{assump:compatibility} is the standard compatibility condition along the Thompson sampling trajectory. Its validity in high-dimensional sparse linear contextual bandits has been established using the transfer lemma and anti-concentration arguments (see Appendix~B in~\cite{chakraborty2023thompsonsamplinghighdimensionalsparse}). In deriving the regret upper bound, we condition on the high-probability event where \eqref{eq:compatibility} holds, and subsequently relax this conditioning at the conclusion of the analysis. 




\begin{theorem}[Regret Bound of SparBL]\label{thm:linear-regret}
Let $r_T$ denote the cumulative regret up to iteration $T$. Under Assumptions~\ref{assump:true-param} and ~\ref{assump:compatibility}, SparBL admits the following regret bound:
\begin{equation*}
r_T \le C_0s_0 \sqrt{T\,\log D},
\end{equation*}
where $D = 1 + N + \binom{N}{2}$, and $s_0$ is the constant specified in Assumption~\ref{assump:true-param}.
\end{theorem}

The proof of Theorem~\ref{thm:linear-regret} is provided in \S\ref{ssec_linear_regret}. This result shows that the cumulative regret scales logarithmically with dimension $D$ and sublinearly with the time horizon $T$. This theoretical result, together with the acquisition function solution introduced in the previous section, is new to the literature of optimization with sparse Bayesian linear models.

\section{Nonparametric Surrogate: Gaussian Process with $p$-Median Prior}\label{sec_BO}

The parametric surrogate offers interpretability and guarantees, but performance may depend on its functional form and sparsity. To address this, we also consider a nonparametric Gaussian process (GP) surrogate. Unlike the parametric model, a GP surrogate makes minimal assumptions and flexibly captures complex nonlinear dependencies.


\subsection{An Overview of the Nonparametric Approach}\label{ssec:overview}



The nonparametric approach follows the same iterative framework but differs in key components. We begin by randomly selecting $\mathcal{T}$ initial solutions $\boldsymbol{x}_1, \cdots, \boldsymbol{x}_{\mathcal{T}}$, and evaluating $f(\boldsymbol{x}_1),\cdots, f(\boldsymbol{x}_\mathcal{T})$ using \eqref{eq_obj}. We then fit a Gaussian-process surrogate (\S\ref{sec_GP}). Next, an acquisition function identifies a feasible trust region (FTR), a subset satisfying feasibility and proximity constraints.Within the FTR, an adaptive swapping search identifies the next candidate $\boldsymbol{x}_{t+1}$ (\S\ref{sec_acq}). We then evaluate $f(\boldsymbol{x}_{t+1})$, and update the surrogate model. When the FTR is exhausted or improvement stalls, a restart shifts the search to a new region (\S\ref{sec_restart}). This process repeats until the evaluation budget $T$ is reached. The final output is the configuration with the lowest observed objective value among all evaluated solutions. In \S\ref{sec_example}, we illustrate the approach with an example problem, and in \S\ref{ssec_theory}, we present the theoretical results. 



\subsection{Gaussian Process}
\label{sec_GP}
A \textit{Gaussian process} $GP(m, k)$ is specified by its mean $m(\boldsymbol{x})=E[f(\boldsymbol{x})]$ and kernel function $k\left(\boldsymbol{x}, \boldsymbol{x}^{\prime}\right)=E \left[(f(\boldsymbol{x})-m(\boldsymbol{x}))(f(\boldsymbol{x}^{\prime})-m(\boldsymbol{x}^{\prime}))\right]$, which represents the covariance between the function values at solutions $\boldsymbol{x}$ and $\boldsymbol{x}^\prime$. Let $\boldsymbol{X}$ represent the set of evaluated solutions, where each solution corresponds to the locations of all units. Similarly, let $\boldsymbol{Y}$ be the set of the corresponding actual objective function values (mean response times), i.e., $\boldsymbol{Y} = \{f(\boldsymbol{x})| \boldsymbol{x}\in \boldsymbol{X}\}$. The Gaussian process distribution captures the joint distribution of all the evaluated solutions, and any subset of these solutions follows a multivariate normal distribution. Under $GP(m, k)$, for a new solution $\boldsymbol{\tilde{x}}$, the joint distribution of $\boldsymbol{Y}$ and the mean response time value $\tilde{y}$ of $\boldsymbol{\tilde{x}}$ is
\begin{equation}
\begin{bmatrix}
\boldsymbol{Y} \\
\tilde{y}
\end{bmatrix} \sim \mathcal{N}\left(m\left(\begin{bmatrix}
\boldsymbol{X} \\
\boldsymbol{\tilde{x}}
\end{bmatrix}\right), \begin{bmatrix}
K(\boldsymbol{X}, \boldsymbol{X})+\sigma^2 \boldsymbol{I} & K(\boldsymbol{X}, \boldsymbol{\tilde{x}}) \\
K(\boldsymbol{X}, \boldsymbol{\tilde{x}})^\top & k(\boldsymbol{\tilde{x}}, \boldsymbol{\tilde{x}})
\end{bmatrix}\right),
\end{equation}
where $K(\boldsymbol{X}, \boldsymbol{X})=\left[k\left(\boldsymbol{x}, \boldsymbol{x}^{\prime}\right)\right]_{\boldsymbol{x}, \boldsymbol{x}^{\prime} \in \boldsymbol{X}}$ is the covariance matrix between the previously evaluated solutions, $K(\boldsymbol{X}, \boldsymbol{\tilde{x}})=\left[k\left(\boldsymbol{x}, \boldsymbol{\tilde{x}}\right)\right]_{\boldsymbol{x} \in \boldsymbol{X}}$ is the covariance vector between the previously evaluated solutions and the new solution, and $k(\boldsymbol{\tilde{x}}, \boldsymbol{\tilde{x}})$ is the covariance between the new solution and itself. In addition, $\sigma^2$ represents the variance of the observed noise, $\boldsymbol{I}$ is the identity matrix, and superscript $\top$ denotes matrix transpose. 

\paragraph{$\boldsymbol{p}$-Median Mean Function.} We use the $p$-Median objective as the mean function $m(x)$, which serves as a prior incorporating structural knowledge from classical facility location theory. We refer to this model as GP-$p$M. This design grounds the Gaussian process in the deterministic $p$-Median framework while allowing the kernel to model additional spatial patterns and dependencies that the $p$-Median formulation alone cannot capture. 

\paragraph{Zero Mean Function.} As a baseline, we also use a zero-mean prior $m(\boldsymbol{x}) = 0$, referred to as GP-zero. A zero mean is commonly adopted when no prior structural knowledge about the problem is available. This baseline allows us to assess the benefit of incorporating problem-specific initialization in Bayesian optimization for facility location problems. In particular, it helps distinguish the performance gains arising from an informed mean function from those contributed by the kernel’s ability to learn residual spatial relationships from data.

\paragraph{Kernel Function.} The kernel function $k(\cdot, \cdot)$ measures solution similarity. In our problem, we design a specific kernel function to measure how similar the new solution is to the existing ones. We let 
\begin{equation}
    k\left(\boldsymbol{x}, \boldsymbol{x}^{\prime}\right)= 
    \underbrace{\exp \left( \sum_{i=1}^{N} \ell_{i}\delta\left(x_{i}, x_{i}^{\prime}\right)/N \right)}_{\text{location-level similarity}} + \underbrace{(\tanh{\gamma})^{\frac{H(\boldsymbol{x}, \boldsymbol{x}^{\prime})}{2}}}_{_{\text{configuration-level similarity}}}, \label{eq:kernel}
\end{equation}
where $\ell_{i}$ and $\gamma$ are trainable parameters. The kernel includes the Kronecker delta function $\delta(x_i, x_i^\prime)$, which equals 1 if $x_i=x_i^\prime$ and 0 otherwise, and the Hamming distance $H(\boldsymbol{x}, \boldsymbol{x}^\prime)$ as defined in \eqref{eq_hamming}.


Kernel~\eqref{eq:kernel} is designed to measure similarity between two solutions at two levels. At the \textit{location-level}, it computes the weighted sum of the differences between the locations of units, where each weight is a learned parameter $\ell_{i}$. At the \textit{configuration-level}, it accounts for the number of differing unit locations between the two solutions. For instance, given two solutions $\boldsymbol{x} = \{1,1,0\}$ and $\boldsymbol{x}^\prime = \{1,0,1\}$, our kernel value is $k\left(\boldsymbol{x}, \boldsymbol{x}^{\prime}\right) = e^{(l_2+l_3)/3}+\tanh \gamma$, because $\delta\left(x_{i}, x_{i}^{\prime}\right)=1$ for $i=2,3$, and $H(\boldsymbol{x}, \boldsymbol{x}^{\prime})=2$.

\begin{lemma}
The kernel $k\left(\boldsymbol{x}, \boldsymbol{x}^{\prime}\right)$ is Hermitian and positive semi-definite.
\end{lemma}

As both components in the kernel $k\left(\boldsymbol{x}, \boldsymbol{x}^{\prime}\right)$ are positive semi-definite, adding two positive semi-definite kernels is also positive semi-definite. 

\paragraph{Posterior Distribution.}
Given solutions $\boldsymbol{X}$ and values $\boldsymbol{Y}$, the \textit{posterior distribution} of a new solution $\boldsymbol{\tilde{x}}$ follows the normal distribution $p(\tilde{y}|\boldsymbol{\tilde{x}},\boldsymbol{X},\boldsymbol{Y}) \sim \mathcal{N}\left(\mu(\boldsymbol{\tilde{x}}; \boldsymbol{X}, \boldsymbol{Y}), \sigma^2(\boldsymbol{\tilde{x}}; \boldsymbol{X}, \boldsymbol{Y})\right)$, where $\mu(\boldsymbol{\tilde{x}}; \boldsymbol{X}, \boldsymbol{Y})$ and $\sigma^2(\boldsymbol{\tilde{x}}; \boldsymbol{X}, \boldsymbol{Y})$ are the mean and variance that are given by
{
\setlength{\abovedisplayskip}{3pt}
\setlength{\belowdisplayskip}{3pt}
\begin{align}
    \mu(\boldsymbol{\tilde{x}}; \boldsymbol{X}, \boldsymbol{Y})&=m(\boldsymbol{\tilde{x}})+K(\boldsymbol{\tilde{x}}, \boldsymbol{X}) [K(\boldsymbol{X}, \boldsymbol{X})+\sigma^2 \boldsymbol{I}]^{-1}(\boldsymbol{Y}-m\left(\boldsymbol{X})\right), \label{GP_mean}\\
    \sigma^2(\boldsymbol{\tilde{x}}; \boldsymbol{X}, \boldsymbol{Y})&=k(\boldsymbol{\tilde{x}}, \boldsymbol{\tilde{x}})-K(\boldsymbol{\tilde{x}}, \boldsymbol{X}) [K(\boldsymbol{X}, \boldsymbol{X})+\sigma^2 \boldsymbol{I}]^{-1} K(\boldsymbol{X}, \boldsymbol{\tilde{x}}). \label{GP_var}
\end{align}
}

The training of GPs directly follows the Bayesian updating rules, where we find the optimal parameters in the kernel by maximizing the likelihood of the training data.

\subsection{Acquisition Function}
\label{sec_acq}
We implement a two-level Gaussian Process (GP) framework to explore the solution space. $GP_{global}$ identifies promising regions, while $GP_{local}$ performs fine-grained searches within them. This structure balances global exploration and local exploitation, improving efficiency in high-dimensional discrete optimization.

Specifically, $GP_{global}$ first identifies a \textit{center solution} that defines the FTR, a constrained subset with high-potential configurations. Within this region, $GP_{local}$, trained only on solutions from $\mathcal{X}_{FTR}$, conducts fine-grained optimization to select the next candidate solution. 

\paragraph{Center Solution.}
$GP_{global}$ selects a center solution $\boldsymbol{x}^{c}$ via a \textit{lower confidence bound} (LCB) acquisition strategy~\citep{srinivas2009gaussian}: 
\begin{equation}\label{eq_center}
    \boldsymbol{x}^{c}=\argmin_{\boldsymbol{x}} \mu(\boldsymbol{x};\boldsymbol{X}, \boldsymbol{Y})-\beta^{1 / 2} \sigma(\boldsymbol{x};\boldsymbol{X}, \boldsymbol{Y}),
\end{equation}
where $\mu(\boldsymbol{x};\boldsymbol{X}, \boldsymbol{Y})$ and $\sigma(\boldsymbol{x};\boldsymbol{X}, \boldsymbol{Y})$ represent the mean and standard deviation, respectively, of the Gaussian process $GP_{global}$, as shown in \eqref{GP_mean} and \eqref{GP_var}, where $\beta$ is a trade-off parameter. The center $\boldsymbol{x}^{c}$ is feasible and promising. We leverage this solution as the starting point for the search within the feasible trust region.

\paragraph{Feasible Trust Region.} The FTR is a constrained search neighborhood around $\boldsymbol{x}^{c}$ whose candidates satisfy distance and cardinality constraints \eqref{eq_feasibility_approx}. The distance constraint restricts how far a new solution can deviate from $\boldsymbol{x}^{c}$, while the cardinality constraint ensures that exactly $p$ units are selected, maintaining feasibility within the location problem. 
We define the feasible trust region $\mathrm{FTR}_{d}\left(\boldsymbol{x}^{c}\right)$ as
\begin{equation}
\mathrm{FTR}_{d}\left(\boldsymbol{x}^{c}\right)=\left\{\boldsymbol{x}\in \{0,1\}^N \mid H(\boldsymbol{x}, \boldsymbol{x}^{c}) \leq d \text{  and  } \sum_{i=1}^N x_i = p \right\} \label{eq_feasible_TR},
\end{equation}
where $\boldsymbol{x}^{c}$ is the center of the FTR, $H(\boldsymbol{x}, \boldsymbol{x}^{c})$ is the Hamming distance between a solution $\boldsymbol{x}$ and the center solution $\boldsymbol{x}^{c}$, and $d$ is the edge-length of the FTR. The edge-length $d$ controls exploration by limiting the Hamming distance from $\boldsymbol{x}^{c}$. 

In the context of optimizing unit locations, an edge-length refers to the maximum number of unit locations in which a new solution can differ from the center solution $\boldsymbol{x}^c$. This means that the FTR defines a region in the search space that is restricted to solutions that are close to the center solution in terms of the number of differing unit locations. The Hamming distance $H\left(\boldsymbol{x}, \boldsymbol{x}^c\right)$ is a measure of the number of positions at which two solutions $\boldsymbol{x}$ and $\boldsymbol{x}^c$ differ, and the FTR constrains this distance to be less than or equal to the edge-length $d$. By limiting the search to this small region of feasible solutions, the optimization algorithm can more efficiently find the next solution that improves the objective function.

\paragraph{Adaptive Swapping in FTR.} The \textit{adaptive swapping search} method conducts an efficient local exploration within the FTR centered at $\boldsymbol{x}^{c}$. This approach iteratively applies a \textit{random swapping operation} $\xi: \{0,1\}^N \mapsto \{0,1\}^N$, where
{
\setlength{\abovedisplayskip}{3pt}
\setlength{\belowdisplayskip}{3pt}
\begin{equation*}
\xi(\boldsymbol{x}) = \{x_1,\cdots, x_j, \cdots, x_i, \cdots, x_{N}\}\qquad \text{s.t.}\qquad x_i\not=x_j,
\end{equation*}
}
which generates new candidate solutions by exchanging values between distinct unit locations in the current solution. Each execution of $\xi$ produces a modified solution where exactly two location indicators have been swapped. Through repeated executions, this mechanism enables controlled movement through the solution space while maintaining feasibility.

The following lemma demonstrates the property of repeated executions of the random swapping function, which forms the basis for our adaptive swapping method design. The proof is presented in \S \ref{ssec_swap_lemma}.

\begin{lemma}
\label{lemma_rand_swap}
Let $\xi^m(\boldsymbol{x})$ be the solution after performing $m$ random swaps for $\boldsymbol{x}$. The Hamming distance $H(\boldsymbol{x}, \xi^m(\boldsymbol{x}))$ satisfies the following properties: 
\begin{enumerate}
	\item For $m=1$, $H(\boldsymbol{x}, \xi^m(\boldsymbol{x})) = 2$.
	\item For $m\geq 2$, $H(\boldsymbol{x}, \xi^m(\boldsymbol{x})) \in \llbracket 0, 2\min\{m, p, N-p\} \rrbracket \cap \mathbb{N}_{2k}$, where $\llbracket a, b \rrbracket$ indicates the interval of all integers between $a$ and $b$ included and $\mathbb{N}_{2k}$ is the set of all even natural numbers. 
\end{enumerate}
\end{lemma}
The lemma shows $H(\boldsymbol{x}, \xi^m(\boldsymbol{x}))$ is even and bounded by $2\min\{m, p, N-p\}$. This property aids us in determining the number of random swaps needed to adequately explore the feasible trust region.

By Lemma \ref{lemma_rand_swap}, we set the number of repetitions of applying the random swapping function to $s(d)=\lfloor\min\{d/2, p, N-p\}\rfloor$ to ensure sufficient exploration. For large $p$ and $N-p$, the FTR exploration is sufficient to reach every solution. For small $p$ or $N-p$, the number of feasible solutions is bounded by problem size, so $p$ or $N-p$ swaps suffice.

We start by adaptively swapping from $\boldsymbol{x}^{c}$ using $GP_{local}$ trained on all observations $(\boldsymbol{X}, \boldsymbol{Y})$ and evaluate candidates with \textit{expected improvement} (EI) function. EI quantifies improvement over the current best $g = \min (\boldsymbol{Y})$ by integrating over the GP posterior:
\begin{equation}
\begin{aligned}
\mathrm{EI}(\boldsymbol{x}) &= \mathrm{E}\left[\max \left(0, g-f\left(\boldsymbol{x}\right)\right) \mid\boldsymbol{X}, \boldsymbol{Y}\right]= \int_{-\infty}^{\infty} \max \left(0, g-f\left(\boldsymbol{x}\right)\right) \varphi(z) \mathrm{d} z\\
   &=\left(g - \mu(\boldsymbol{x};\boldsymbol{X}, \boldsymbol{Y})\right)\Phi\left(\frac{g - \mu(\boldsymbol{x};\boldsymbol{X}, \boldsymbol{Y})}{\sigma(\boldsymbol{x};\boldsymbol{X}, \boldsymbol{Y})}\right)+ \sigma(\boldsymbol{x};\boldsymbol{X}, \boldsymbol{Y})  \varphi\left(\frac{g - \mu(\boldsymbol{x};\boldsymbol{X}, \boldsymbol{Y})}{\sigma(\boldsymbol{x};\boldsymbol{X}, \boldsymbol{Y})}\right),
\end{aligned}
\end{equation}
where $\varphi(\cdot)$ and $\Phi(\cdot)$ denote the standard normal PDF and CDF respectively. This closed-form expression results from a reparameterization of the GP posterior.

Each iteration replaces the current FTR best with any higher-EI solution. The algorithm executes $K$ such iterations before selecting the highest-EI candidate for actual evaluation. The complete procedure is formalized in Algorithm \ref{ASS_algo}.


\begin{algorithm}
\caption{Subroutine: Adaptive Swapping Search}
\label{ASS_algo}
{\fontsize{10}{16}\selectfont
\begin{algorithmic}[1]

\State \textbf{Input:} Feasible trust region with center $\boldsymbol{x^c}$; edge-length $d$; total number of iterations $K$.
\State $k \leftarrow 0$; $\boldsymbol{x_0} \leftarrow \boldsymbol{x^c}$.
\While {$k < K$}
\State Conduct $s(d)$ random swaps to obtain a new candidate solution $\xi^{s(d)}(\boldsymbol{x_k})$.

\If{EI$(\xi^{s(d)}(\boldsymbol{x_k})) > $ EI$(\boldsymbol{x_k})$}

\State $\boldsymbol{x_{k+1}} \leftarrow \xi^{s(d)}(\boldsymbol{x_k})$.

\Else
\State $\boldsymbol{x_{k+1}} \leftarrow \boldsymbol{x_{k}} $.

\EndIf

\State $k \leftarrow k + 1$.
\EndWhile
\Statex \textbf{Output:} New candidate solution $\boldsymbol{x_K}$.
\end{algorithmic}}
\end{algorithm}

\begin{algorithm}
\caption{Gaussian Process with $p$-Median Mean Prior (GP-$p$M)}
\label{BO_algo}
{\fontsize{10}{16}\selectfont
\begin{algorithmic}[1]
\State \textbf{Input:} Initial sample size $\mathcal{T}$; initial edge-length $d_0$; evaluation budget $T$; $\alpha_s, \alpha_f, n_s, n_f, \beta$. 
\State Randomly select $\mathcal{T}$ initial samples $\boldsymbol{X}_{\mathcal{T}} = \{\boldsymbol{x}_1, ..., \boldsymbol{x}_{\mathcal{T}}\}$ and evaluate $\boldsymbol{Y}_{\mathcal{T}} = \{f(\boldsymbol{x}_1),f(\boldsymbol{x}_2),\cdots, f(\boldsymbol{x}_{\mathcal{T}})\}$.
\State Set $t \leftarrow \mathcal{T}+1$, restart $\leftarrow$ True, $\boldsymbol{X}^\prime_{\mathcal{T}}\leftarrow \boldsymbol{X}_{\mathcal{T}}$, $\boldsymbol{Y}^\prime_{\mathcal{T}}\leftarrow \boldsymbol{Y}_{\mathcal{T}}$.
\While {$t < T$}
\If{restart is True}
\State Reset edge-length $d \leftarrow d_0$, restart $\leftarrow$ False, $y^\prime \leftarrow \infty$.
\State Update Gaussian process $GP^*$ with $(\boldsymbol{X}^\prime_{t-1}, \boldsymbol{Y}^\prime_{t-1})$
\State Determine a new center of the feasible trust region $\boldsymbol{x^c}$ using (\ref{eq_center}) with $GP^*$.
\Else


\State Update Gaussian process $GP$ with $(\boldsymbol{X}_{t-1}, \boldsymbol{Y}_{t-1})$.

\State Perform \textit{Adaptive Swapping Search} with $GP$ to obtain the next candidate solution $\boldsymbol{x_t}$.

\State Evaluate the objective function $y_t=f(\boldsymbol{x_t})$ and update $\boldsymbol{X}_t \leftarrow \boldsymbol{X}_{t-1}\cup \{\boldsymbol{x_t}\}$, $\boldsymbol{Y}_t \leftarrow \boldsymbol{Y}_{t-1}\cup \{y_t\}$.

\If{$y_t < y^\prime$} 
\State Update the best solution in this FTR: $y^\prime \leftarrow y_t, \boldsymbol{x}^\prime \leftarrow \boldsymbol{x}_t$.
\EndIf

\If{$y^\prime$ has improved $n_s$ times} 
\State Increase the edge-length $d \leftarrow d \alpha_s$.
\EndIf

\If{$y^\prime$ remains unchanged for $n_f$ consecutive iterations}
\State Reduce the edge-length $d \leftarrow d \alpha_f$.
\EndIf

\EndIf
\If{$\lfloor d \rfloor < 2$}
\State Set restart $\leftarrow$ True, update $\boldsymbol{X}^\prime_t \leftarrow \boldsymbol{X}^\prime_{t-1}\cup \{\boldsymbol{x}^\prime\}$, $\boldsymbol{Y}^\prime_t \leftarrow \boldsymbol{Y}^\prime_{t-1}\cup \{y^\prime\}$.
\Else
\State Keep $\boldsymbol{X}^\prime_t \leftarrow \boldsymbol{X}^\prime_{t-1}$, $\boldsymbol{Y}^\prime_t \leftarrow \boldsymbol{Y}^\prime_{t-1}$.
\EndIf
\State Move to the next iteration $t \leftarrow t + 1$.
\EndWhile
\Statex \textbf{Output:}  Best solution $\boldsymbol{x}^*$ with the minimum mean response time $\boldsymbol{y}^*$ observed in  $(\boldsymbol{X}_{T}, \boldsymbol{Y}_{T})$.
\end{algorithmic}}
\end{algorithm}

\subsection{Restart Mechanism}
\label{sec_restart}

Our algorithm incorporates a restart mechanism to enhance global search capability and avoid local optima stagnation, drawing inspiration from restart strategies in non-stationary bandit problems \citep{besbes2014stochastic,wan2021think}. This mechanism dynamically adjusts the FTR's edge-length $d$ based on recent search performance: each solution evaluation is classified as either a success (improving upon the current best solution) or a failure (no improvement). The edge-length undergoes contraction when encountering $n_f$ consecutive failures ($d \leftarrow d \times \alpha_f$ where $\alpha_f \leq 1$) or expansion after $n_s$ successes ($d \leftarrow d \times \alpha_s$ where $\alpha_s \geq 1$). A critical restart condition triggers when $\lfloor d \rfloor < 2$, ensuring the algorithm escapes local search regions and maintains complete solution space exploration.

When a restart condition is met, the algorithm performs three key operations: (i) resetting the edge-length to its initial value $d$, (ii) updating $GP_{global}$'s training set with the best solution from the current FTR, and (iii) establishing a new FTR centered on the most promising region identified by the updated $GP_{global}$. This process continues until exhausting the evaluation budget $T$, with the global minimum solution retained as the final output. The complete procedure is formally presented in Algorithm~\ref{BO_algo}.

\subsection{An Illustrative Example} 
\label{sec_example}

This section demonstrates our GP-$p$M approach through a small working example. 
\begin{example}
Consider a service region partitioned into $M=5$ subregions where we aim to optimally locate $p=2$ service units among $N=5$ candidate locations to minimize mean response time. The initial trust region edge-length is set to $d_0=2$. In addition, we set $n_s=n_f=3, \alpha_s=1.5, \alpha_f=0.75$, and $\beta=25$. Applying the approach yields $\boldsymbol{x}^*=\{0,1,0,0,1\}$ with $f\left(\boldsymbol{x}^*\right)=4.15$, which we confirm as optimal by enumeration. The solution space is shown in Figure~\ref{hypercube_5}.

\begin{figure}[htbp]
\centering
  \includegraphics[width=10cm]{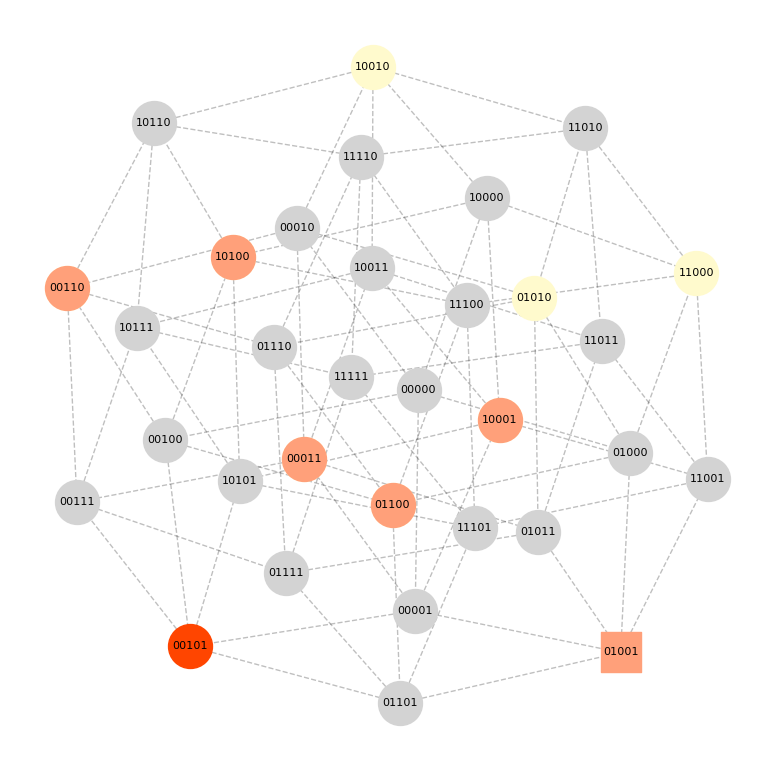}
  \caption{Illustration of the solution space.}
  \raggedright
    {\scriptsize \textit{Note.} The current FTR center, denoted by \{0,0,1,0,1\}, is highlighted in red. Solutions within the FTR are shown in orange, while other feasible solutions are marked in yellow. Infeasible solutions are indicated in gray. The optimal solution \{0,1,0,0,1\} is marked with a square. Dashed edges represent a Hamming distance of 1 between solutions.}
  \label{hypercube_5}
  \vspace{-0.2cm}
\end{figure}

At the start, the algorithm initializes with three random configurations: $\boldsymbol{x}_1 = \{1, 0, 1, 0, 0\}$, $\boldsymbol{x}_2 = \{0, 0, 1, 0, 1\}$ and $\boldsymbol{x}_3 = \{0, 0, 1, 1, 0\}$. Evaluating these solutions using the spatial hypercube model yields $f(\boldsymbol{x}_1) = 6.31$, $f(\boldsymbol{x}_2) = 5.28$, $f(\boldsymbol{x}_3) = 7.45$. Next, we train a global Gaussian process $GP_{global}$ with the initial observations, which identifies $\boldsymbol{x}^{c} =\{0, 0, 1, 0, 1\}$ as the FTR center based on the lowest LCB value. We highlight the FTR and the center in Figure \ref{hypercube_5}.

Next, we train a Gaussian process $GP$ with all the observations $X_3 = \{\boldsymbol{x}_1, \boldsymbol{x}_2, \boldsymbol{x}_3\}, \boldsymbol{Y}_3 = \{f(\boldsymbol{x}_1), f(\boldsymbol{x}_2), f(\boldsymbol{x}_3)\}$.
We use the adaptive swapping method with $GP$ to search for solutions within the FTR. In the first iteration, the solution $\{1, 0, 0, 0, 1\}$ has the highest EI value of 0.08. We let $\boldsymbol{x}_4= \{1, 0, 0, 0, 1\}$ and use the approximate hypercube model to compute its value as $f(\boldsymbol{x}_4) = 4.71$, which is then used to update $GP$. In the next iteration, the algorithm selects solution $\{0, 1, 0, 0, 1\}$ with the highest EI value of 0.19 as the next candidate solution $\boldsymbol{x}_5$. Upon evaluation, we obtain a value of $f(\boldsymbol{x}_5) = 4.15$.

The algorithm obtained the optimal solution $\boldsymbol{x^*} = \boldsymbol{x}_5 = \{0, 1, 0, 0, 1\}$ without a restart. However, for more complex problems that require more iterations, the edge-length will shrink after $n_f$ failures, and the algorithm will trigger the restart mechanism to update $GP_{global}$ and explore a new FTR. 
\end{example}

\subsection{Theoretical Results}\label{ssec_theory}

The next theorem shows the method converges to OPT$_{H}$ in finite iterations under the stated assumptions.

\begin{theorem}[Optimality]
\label{thm:convergence}
Let $\mathcal{S}$ be the set of feasible solutions to the $p$-Hypercube, and let $f: \mathcal{S} \rightarrow \mathbb{R}$ be the corresponding objective function. Let $\{\boldsymbol{x}_t\}$ be a sequence of solutions generated by our algorithm, and define $g_{t} = \min_{k\leq t} f(\boldsymbol{x}_{k})$. Then, for all $t$, we have $\boldsymbol{x}_{t} \in \mathcal{S}$, and $\lim_{t \rightarrow \infty} g_t = \text{OPT}_H$. Additionally, the algorithm converges in a finite number of iterations.
\end{theorem}

We provide a detailed proof of Theorem \ref{thm:convergence} in \S \ref{ssec_optimality_proof}. In summary, we establish the convergence of Theorem \ref{thm:convergence} by demonstrating that the edge-length of the FTR will always be reduced by the design of our algorithm, leading to a restart for every FTR explored. This guarantees that the algorithm will not get trapped in any local optima. Then, we prove optimality using the monotone convergence theorem~\citep{bibby1974axiomatisations} because $g_t$ is monotonically non-increasing and bounded. 

It is also important to show that our algorithm converges rapidly as the search space is large. 
Let $r_{v} = f(\boldsymbol{x}_{v}^{*}) - \text{OPT}_H$ denote the regret for the $v$-th restart, where $\boldsymbol{x}_{v}^{*}$ is the optimum in the FTR at the $v$-th restart. 

\color{black}
\begin{theorem}[Regret Bound of GP-$p$M]
\label{thm:bound}
    Define $\bar{r}_V = \frac{1}{V}\sum_v^V r_v$ as the average regret up to the $V$-th restart. Assume that a sample function from the Gaussian process model defined by the kernel function $k$ passes through all local optima of $f$. Then, for any $\delta \in(0,1)$, we have 
\begin{equation}
    \operatorname{Pr}\left\{\bar{r}_V \leq \sqrt{8 C\beta_{V}  \kappa_V / V}\right\} \geq 1-\delta, 
\end{equation}
where $\beta_V = 2\log(\pi^2 V^2 {N \choose p}/6\delta )$, $\kappa_V = \mathcal{O}(2^N \log V)$, and $C=1 / \log \left(1+\sigma^{-2}\right)$.
\end{theorem}

We present the proof of Theorem~\ref{thm:bound} in \S \ref{ssec_bound}. The regret bound states that, with high probability, the average regret is bounded by a sublinear function of the number of restarts.

\section{Numerical Experiments}\label{sec_numerical}

\begin{figure}[btph]
\vspace{-20pt}
    \centering
    \begin{minipage}{\linewidth}
        \centering
        \includegraphics[width=\textwidth]{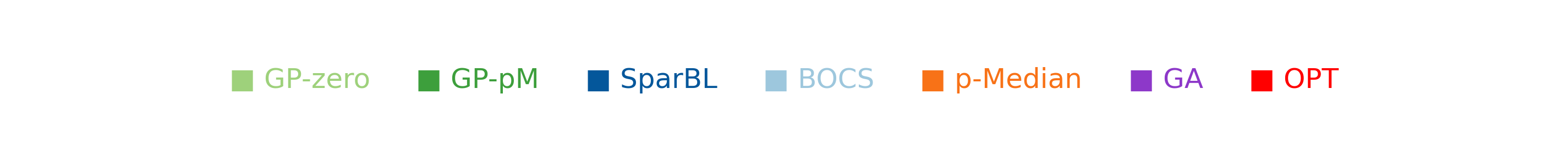}
    \end{minipage}
    
    \vspace{-0.5em}
    \begin{minipage}[t]{0.32\linewidth}
        \centering
        \includegraphics[width=\linewidth]{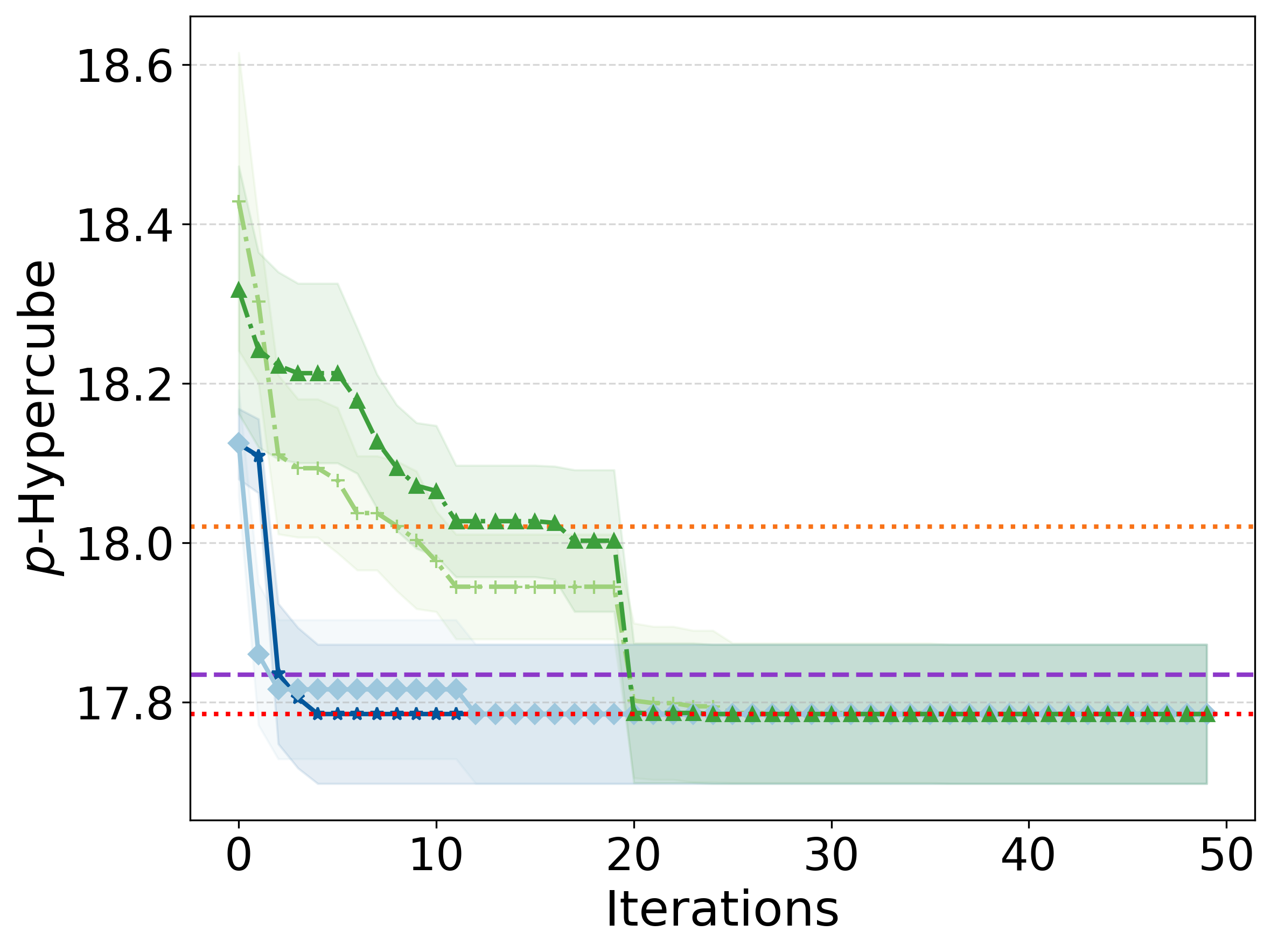}
        \vspace{0.1em}
        \small (a) $N = 10,\ p = 5$
    \end{minipage}%
    \hfill
    \begin{minipage}[t]{0.32\linewidth}
        \centering
        \includegraphics[width=\linewidth]{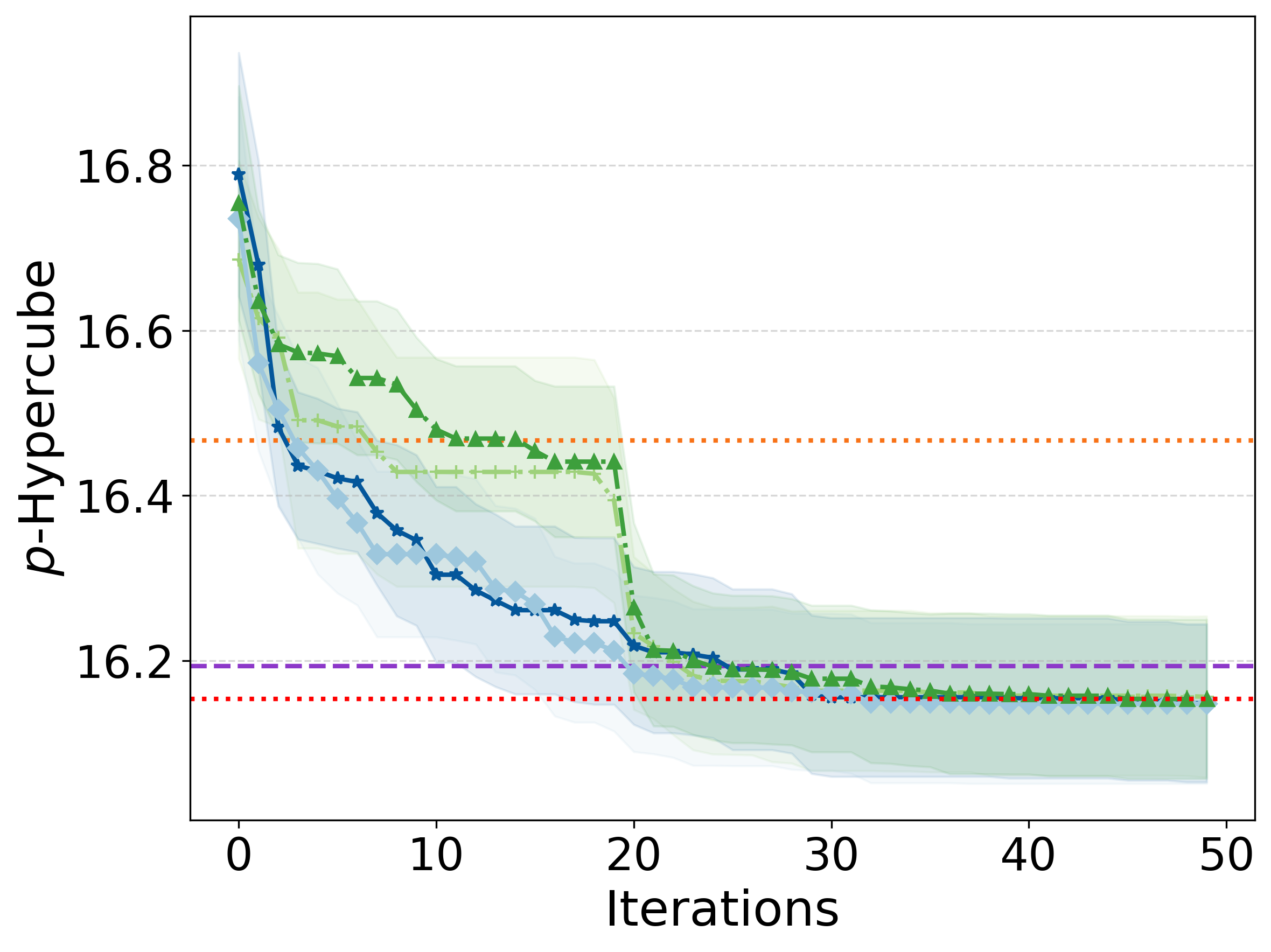}
        \vspace{0.1em}
        \small (b) $N = 20,\ p = 10$
    \end{minipage}%
    \hfill
    \begin{minipage}[t]{0.32\linewidth}
        \centering
        \includegraphics[width=\linewidth]{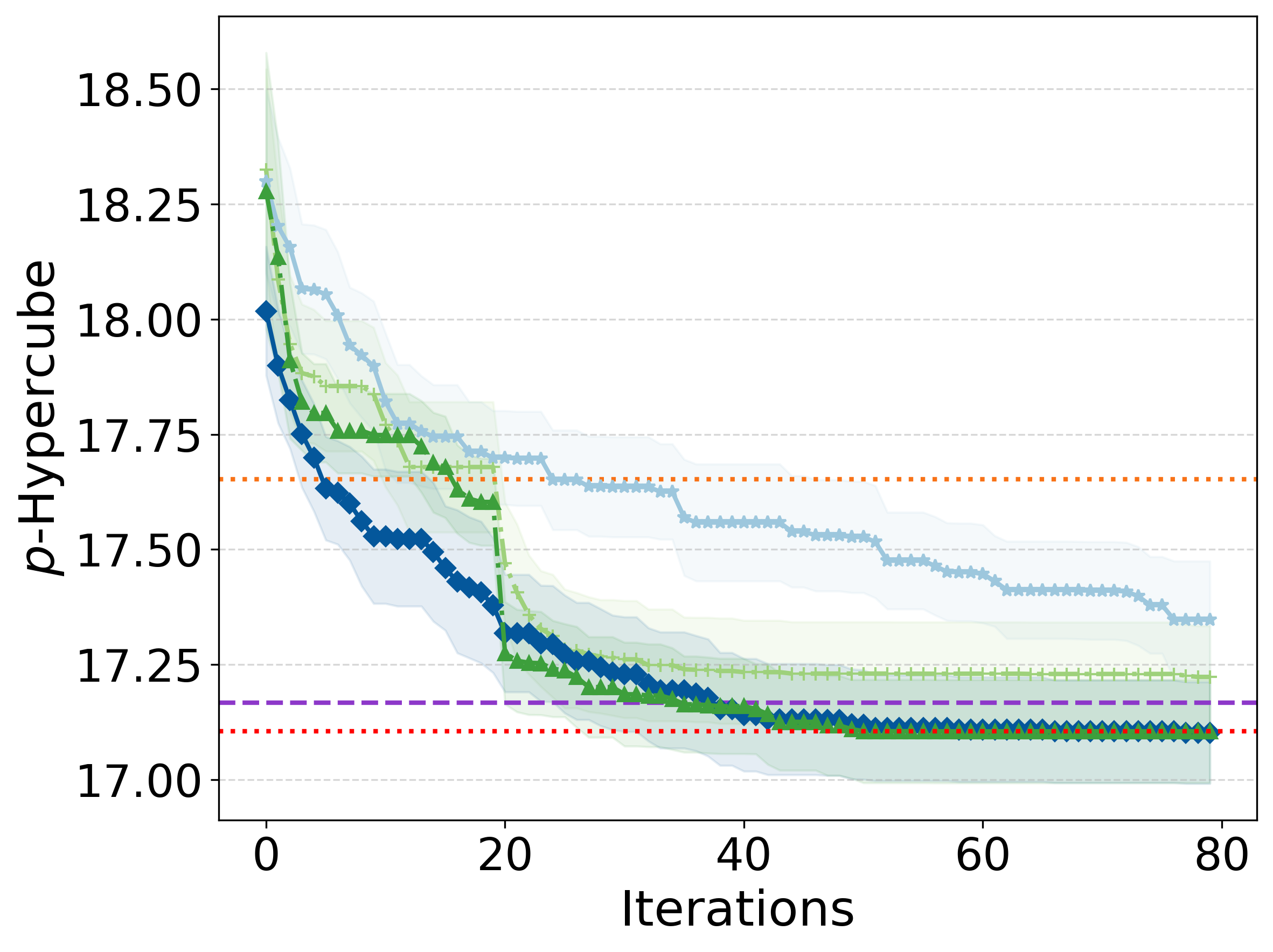}
        \vspace{0.1em}
        \small (c) $N = 30,\ p = 5$
    \end{minipage}
    
    \begin{minipage}[t]{0.32\linewidth}
        \centering
        \includegraphics[width=\linewidth]{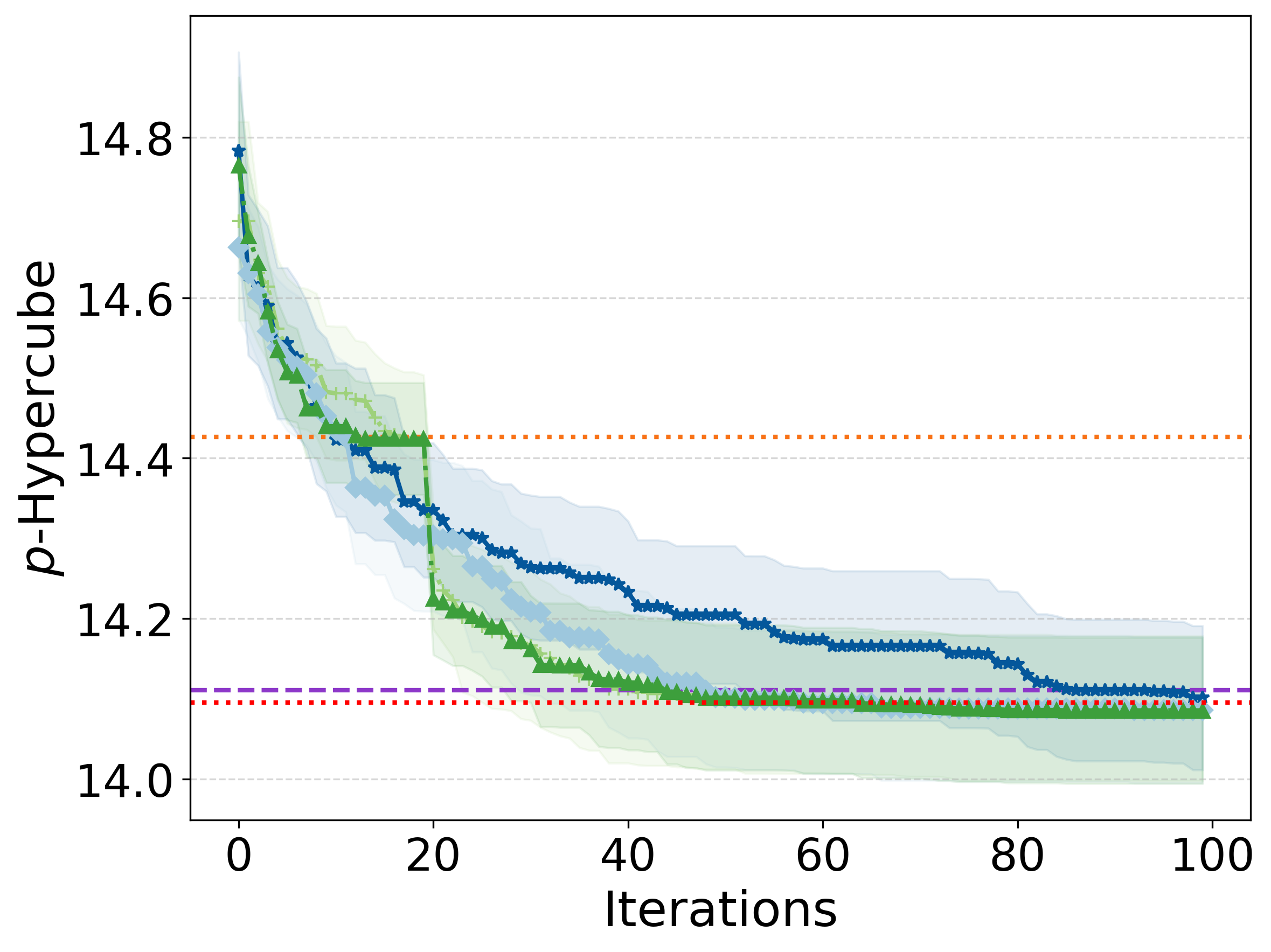}
        \vspace{0.1em}
        \small (d) $N = 30,\ p = 15$
    \end{minipage}%
    \hfill
    \begin{minipage}[t]{0.32\linewidth}
        \centering
        \includegraphics[width=\linewidth]{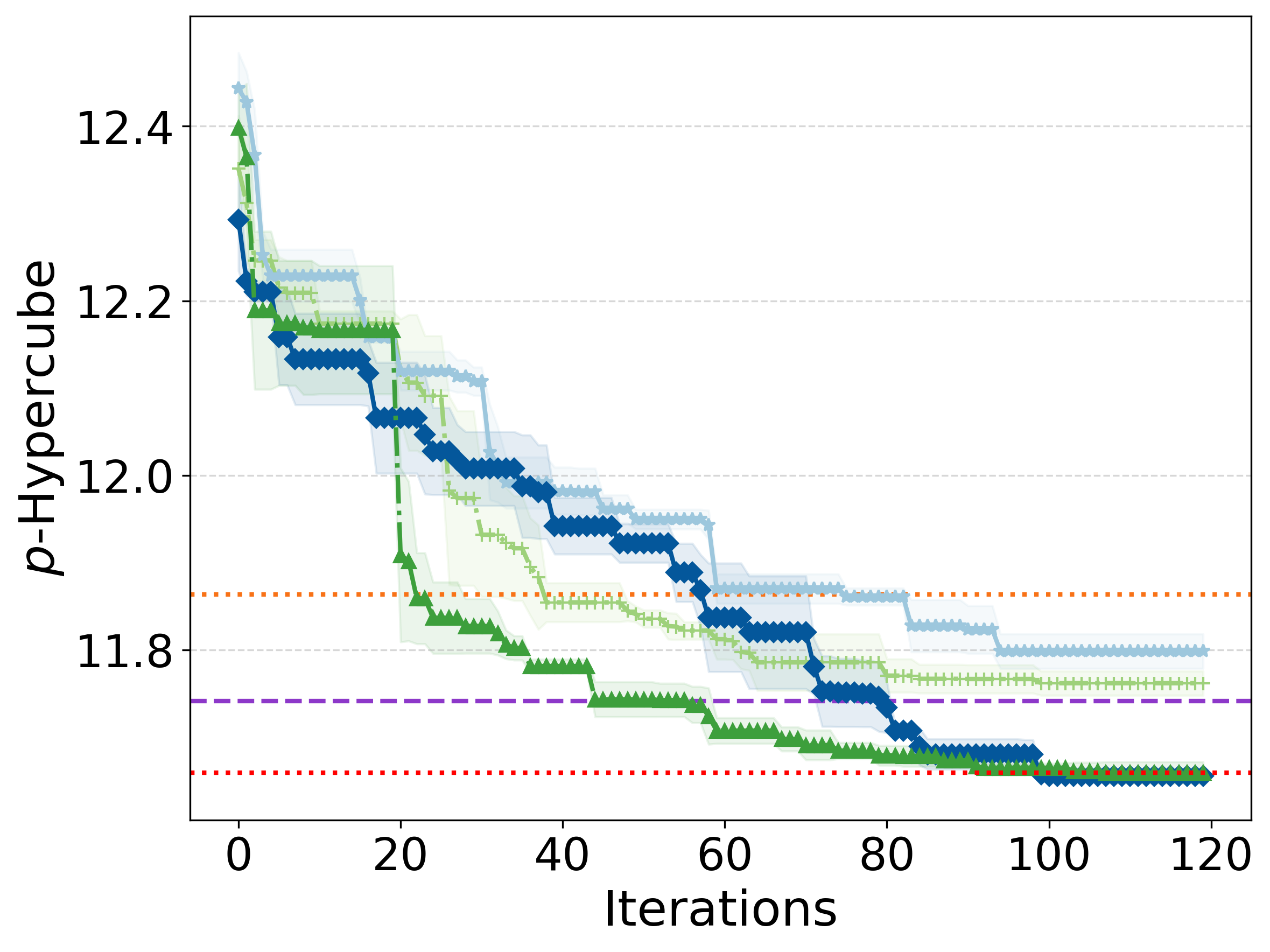}
        \vspace{0.1em}
        \small (e) $N = 50,\ p = 20$
    \end{minipage}%
    \hfill
    \begin{minipage}[t]{0.32\linewidth}
        \centering
        \includegraphics[width=\linewidth]{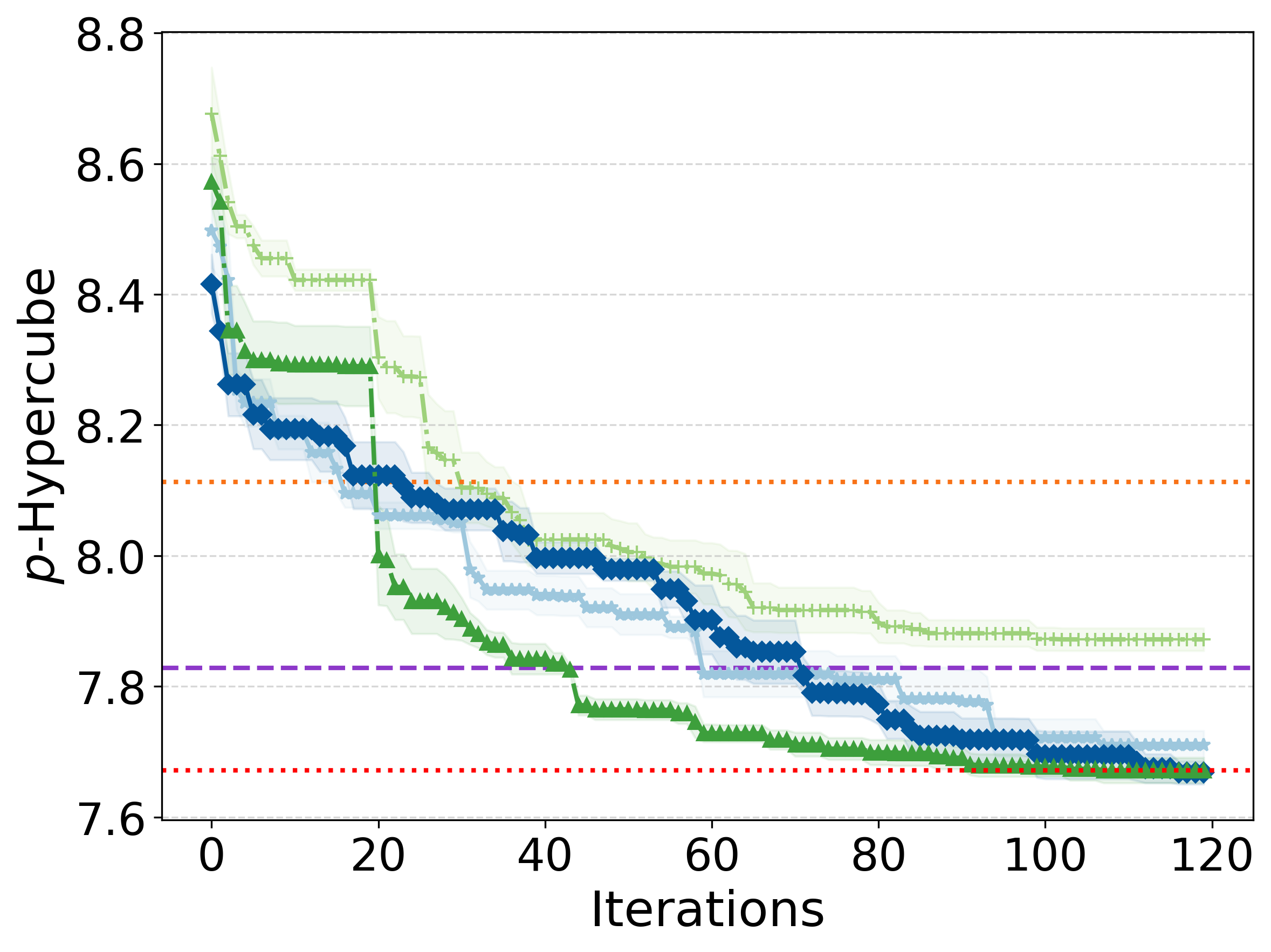}
        \vspace{0.1em}
        \small (f) $N = 50,\ p = 30$
    \end{minipage}
    \vspace{0.1em}
    \caption{Results for simulation experiments.}
    \raggedright
    {\scriptsize \textit{Note.} $N$ denotes the number of candidate stations, and $p$ the number of servers to be allocated. Results show consistent convergence across problem sizes, with the mean and 95\% confidence intervals computed from 10 independent runs. For the $p$-Median and GA benchmarks, only the final outcomes are reported.}
    \label{fig_numerical}
\end{figure}

This section presents systematic numerical evaluations of our approaches, including: (i) SparBL: Implementation of Algorithm~\ref{algo:BO_with_submodular}; (ii) GP-$p$M: Implementation of Algorithm~\ref{BO_algo} with the $p$-Median prior mean function; (iii) GP-zero: Implementation of Algorithm~\ref{BO_algo} with zero mean function, as a special case of GP-$p$M. 

We compare our solutions against the following baselines: (i) classical $p$-Median solution; (ii) BOCS (Bayesian Optimization of Combinatorial Structures)~\citep{baptista2018bayesian}; (iii) Genetic algorithm (GA) for emergency deployment~\citep{geroliminis2011hybrid}. The simulation environment consists of a $10\times10$ grid ($M=100$ subregions) with randomly sampled facility locations.

As demonstrated in Figure~\ref{fig_numerical}, our SparBL and GP-$p$M algorithms consistently achieve optimal solutions across all configurations. For small-scale problems (Fig~\ref{fig_numerical}. a--b), both methods perfectly match the globally optimal solutions obtained through exhaustive enumeration. In medium and large-scale problems (Fig~\ref{fig_numerical}. c--f), our approaches maintain their optimal performance while significantly outperforming the baselines. 
The consistent optimality of SparBL and GP-$p$M across all test cases demonstrates their robustness and scalability, particularly in high-dimensional optimization scenarios where benchmarks fail to explore the solution space effectively.

\section{An Ambulance Location Problem Using St. Paul, MN, Data}\label{sec_realdata}

This section applies the proposed methods to a real-world emergency medical service system in St.~Paul, Minnesota.
Section~\ref{ssec:setup} introduces the operational background of the St.~Paul system.
Section~\ref{ssec:performance} evaluates the performance of our approaches against benchmark methods.
Section~\ref{ssec:asymptotic_case} examines how optimal unit placements evolve under different arrival rates and illustrates the asymptotic properties established in Section~\ref{ssec:asymptotic}.
We also discuss the practical implications of our findings for real-world deployment in this section. 

\subsection{Background}\label{ssec:setup}
Based on 2020 census data, St. Paul has a population of 311,527. Our study utilizes 2014 emergency medical call records containing ambulance locations and incident reports (total of 30,911 cases), with each call's arrival time and nearest street intersection recorded. All data were made anonymous to comply with HIPAA regulations. The ambulance location optimization problem positions 9 ambulances among 17 candidate stations to minimize average response time. 

The St. Paul Department reports show an average turnout time of 1.75 minutes across all stations for medical calls, though our methodology in \S\ref{sec_problem} accommodates location-dependent variations. We aggregated demand by partitioning the city into 71 census tracts, calculating travel times from each tract centroid to candidate stations via the Google Maps API. Census tracts are small, standardized geographic regions defined by the U.S. Census Bureau that approximate neighborhoods.

\begin{figure}[hbtp]
    \centering
    \begin{subfigure}[t]{0.48\textwidth}
        \centering
        \includegraphics[width=\linewidth]{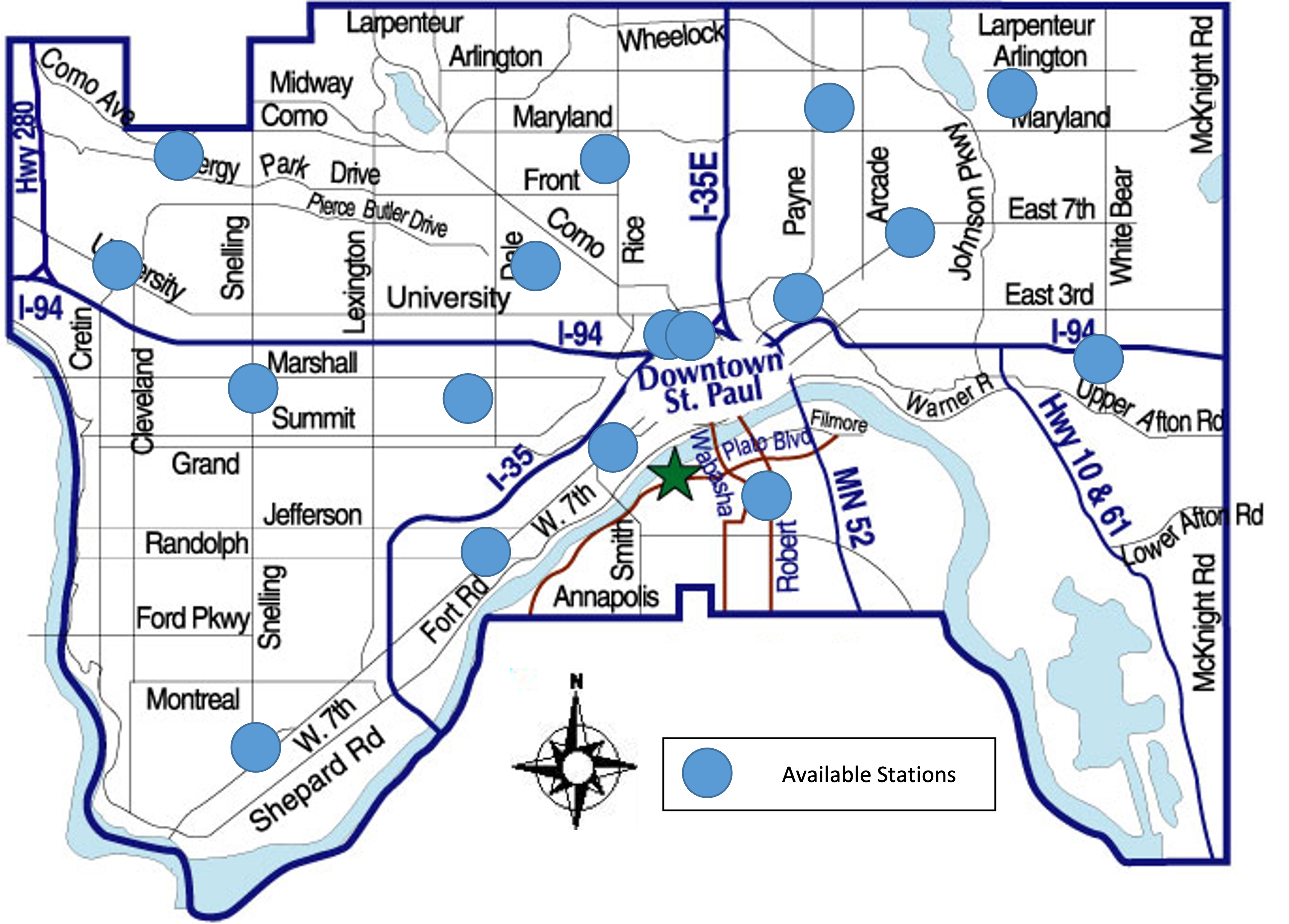}
\caption{\scriptsize St. Paul, MN map showing the 17 available station locations.}
\label{fig_StPaul}
    \end{subfigure}%
    \hfill
    \begin{subfigure}[t]{0.48\textwidth}
        \centering
        \includegraphics[width=\linewidth]{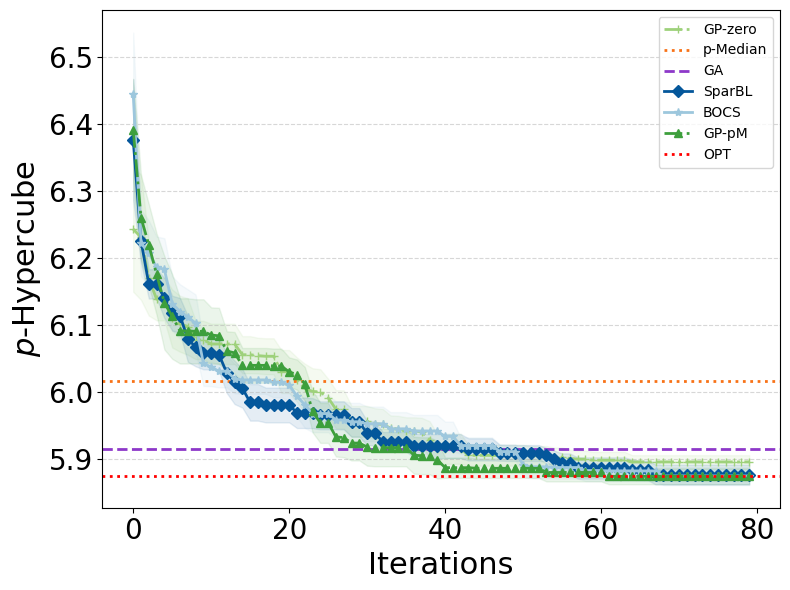}
        \caption{\scriptsize Comparison of our algorithms versus baselines for St. Paul.}
        \label{fig_st_paul}
    \end{subfigure}
    
    \begin{subfigure}[t]{0.48\textwidth}
        \centering
        \includegraphics[width=\linewidth]{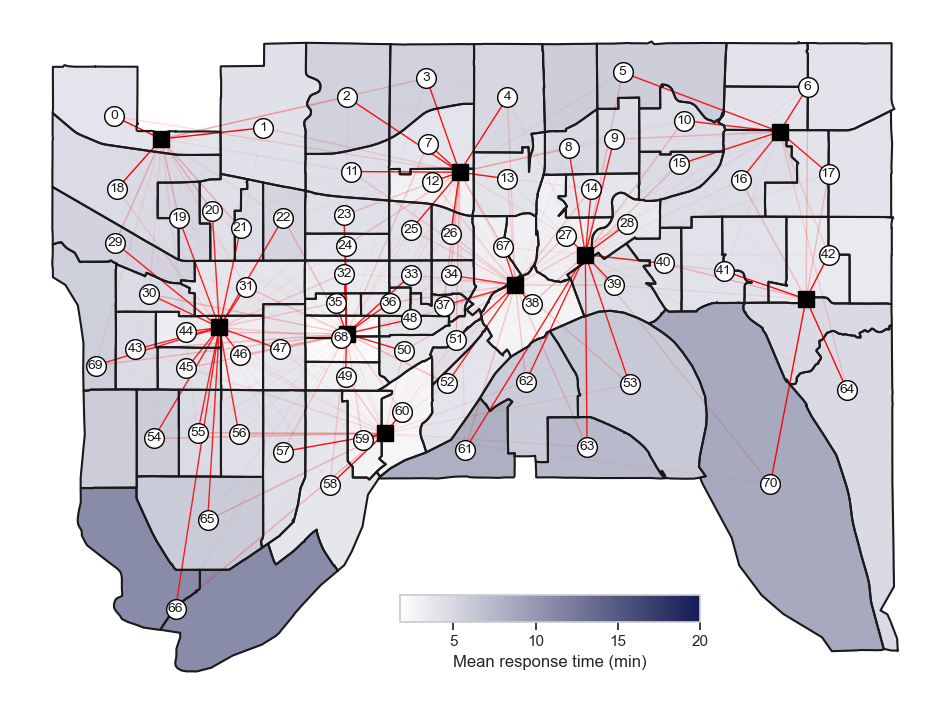}
        \caption{\scriptsize Unit locations of our solutions (coinciding with OPT).}
        \label{fig_StPaul_0}
    \end{subfigure}%
    \hfill
     \begin{subfigure}[t]{0.48\textwidth}
        \centering
        \includegraphics[width=\linewidth]{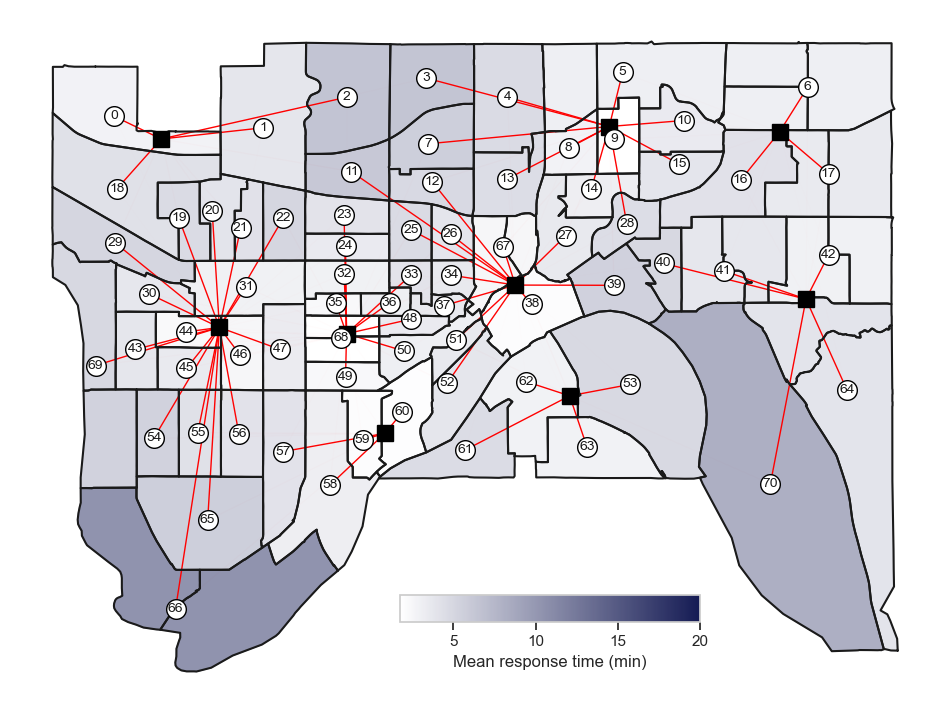}
        \caption{\scriptsize  Unit locations for $p$-Median solution.}
        \label{fig_pmedian_location}
    \end{subfigure}%
    \caption{Results for St. Paul case study.}
    \raggedright
    {\scriptsize \textit{Note.} (a) The map of St. Paul shows the 17 available stations, where nine units are to be placed. (b) St. Paul is a relatively small system, where both GP-$p$M and SparBL converge rapidly to the global optimum. (c) Color intensity within each subregion reflects its average response time, with darker shades indicating longer times. Solid lines connect subregions to stations, with thicker lines denoting higher probabilities of being served by that station. (d) In the $p$-Median formulation, each subregion is assumed to be served by exactly one station.}
    \label{fig:stpaul_overview}
\end{figure}

\subsection{Performance of Our Approaches}\label{ssec:performance}

We present the census tracts and candidate station locations in Figure~\ref{fig:stpaul_overview}, together with the performance of our methods and their comparison to benchmarks. 

Figure~\ref{fig_StPaul} presents the map of St. Paul with its 17 available stations, from which nine units are to be deployed. Figure~\ref{fig_st_paul} compares our methods against the baselines described in~\S\ref{sec_numerical}. By exhaustively enumerating all $\binom{17}{9} = 24,310$ possible configurations, we confirm that our approach identifies the globally optimal deployment, achieving the minimal average response time of 5.86 minutes within about 60 iterations. This rapid convergence is noteworthy given the vast number of possible configurations. In contrast, the $p$-Median average response time remains above the optimal objective, and other baselines either converge more slowly or fail to reach the optimum. Since St. Paul represents a relatively small-scale system, certain baseline methods still perform reasonably; however, as demonstrated in~\S\ref{sec_numerical}, their performance deteriorates markedly in  settings with more units.

Figure~\ref{fig_StPaul_0} shows the ambulance deployment identified by both of our approaches, which result in the same solution, the global optimum. Census tracts are depicted as circles numbered 0–70, and unit locations are marked by solid squares. The shading of each subregion reflects its average response time, with darker colors indicating longer delays. The spatial distribution of response times follows expected patterns: centrally located tracts and those near ambulance stations experience faster responses, while peripheral areas face longer delays. We also report the percentage of calls served by each station, which measures unit workload. Solid lines connect subregions to their assigned stations, with line thicknesses  proportional to the probability of being served. Most calls in each subregion are handled by the first-preferred station, with only a small fraction diverted to others.

We also compared ambulance locations under the $p$-Median solution with those in the optimal configuration. Figure~\ref{fig_pmedian_location} displays the $p$-Median deployment. In this formulation, each subregion is assigned to exactly one station, so every subregion is connected by a single line. The $p$-Median solution yields an average response time of 6.01 minutes, worse than the optimal 5.86 minutes. Relative to the optimal configuration, two units are shifted from central locations (tracts \circnum{27}/\circnum{28}/\circnum{39} and \circnum[0.6 pt]{7}/\circnum{12}) to peripheral areas (tracts \circnum[0.6 pt]{9} and \circnum{62}). This relocation reflects the $p$-Median assumption that all units are always available, which places heavy reliance on a central unit. Consequently, response times increase in several central tracts 
(\circnum[0.6 pt]{2}–\circnum[0.6 pt]{4}, \circnum[0.6 pt]{7}, \circnum{11}–\circnum{13}), even though peripheral regions improve. Such central delays would be amplified in larger cities with denser cores. Even in St. Paul, where the population is more dispersed, the central area remains the busiest, making the $p$-Median solution suboptimal.




\subsection{Optimal Unit Locations Under Varying Arrival Rates}\label{ssec:asymptotic_case}
\subsubsection{Asymptotic Results}

\begin{figure}[hbtp]
    \centering
    \begin{subfigure}[t]{0.48\textwidth}
        \centering
        \includegraphics[width=\linewidth]{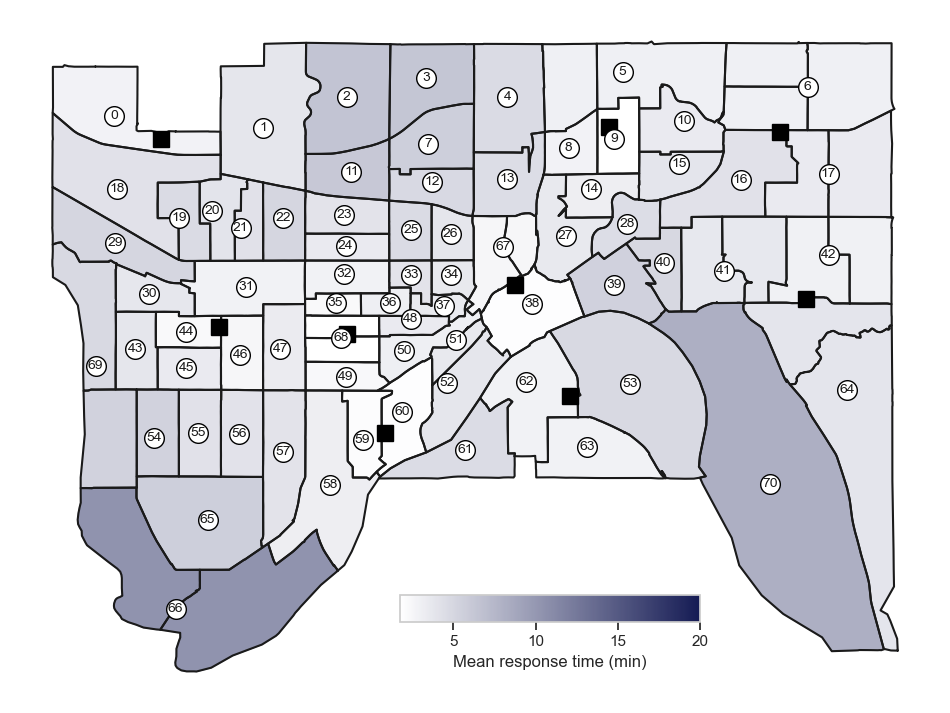}
\caption{\scriptsize $\theta= 10^{-3}$}
\label{fig_lambda_0}
    \end{subfigure}%
    \hfill
    \begin{subfigure}[t]{0.48\textwidth}
        \centering
        \includegraphics[width=\linewidth]{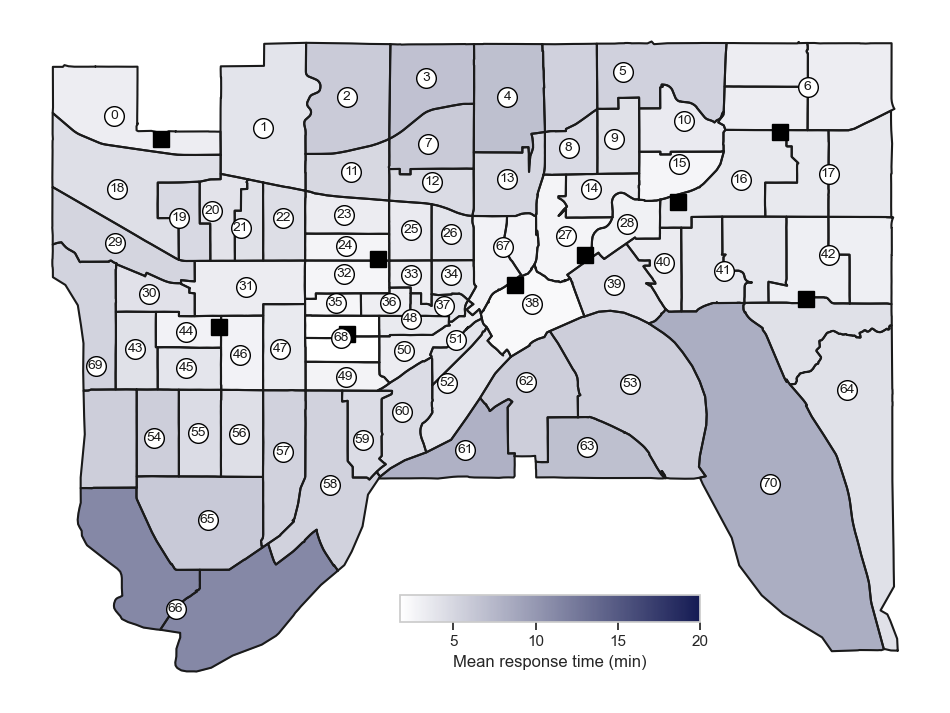}
        \caption{\scriptsize $\theta=2$}
        \label{fig_lambda_2}
    \end{subfigure}
    
    \begin{subfigure}[t]{0.48\textwidth}
        \centering
        \includegraphics[width=\linewidth]{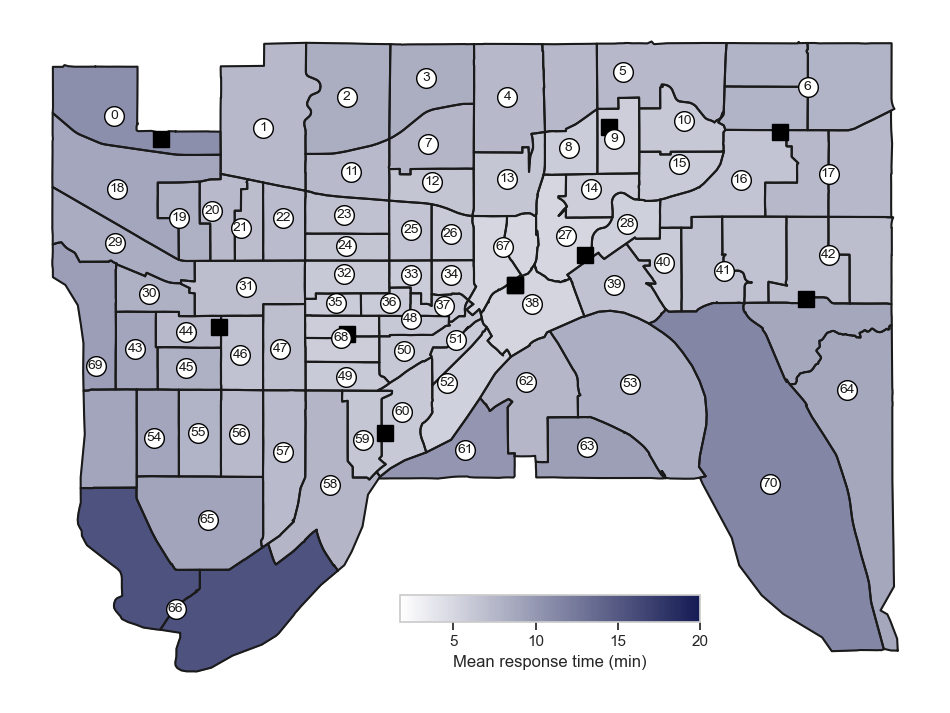}
        \caption{\scriptsize $\theta= 6$}
        \label{fig_lambda_6}
    \end{subfigure}%
    \hfill
     \begin{subfigure}[t]{0.48\textwidth}
        \centering
        \includegraphics[width=\linewidth]{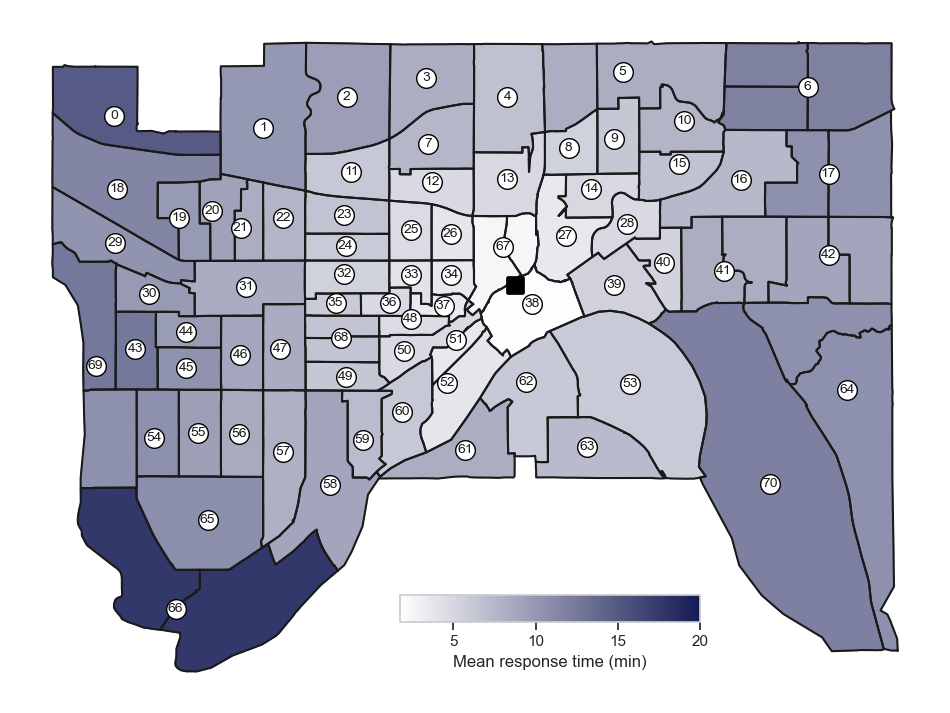}
        \caption{\scriptsize $\theta=10^3$}
        \label{fig_lambda_inf}
    \end{subfigure}
    \caption{Results for St. Paul case study under varying arrival rates.}
    \label{fig:vary_lambda}
\end{figure}

We next examine how the optimal deployment changes as the call arrival rate increases. Figure~\ref{fig:vary_lambda} reports the optimal solutions for the St. Paul case study under different scaling factors $\theta$, which illustrate the asymptotic results in Section~\ref{ssec:asymptotic}.

    

In the vanishing-load regime (Figure~\ref{fig_lambda_0}), the optimal configuration coincides with the weighted $p$-Median benchmark established in Theorem~\ref{thm:asym-vanishing-lambda}. In fact, it is identical to the solution shown in Figure~\ref{fig_pmedian_location}. Units are distributed broadly across the city, with coverage directed toward high-demand neighborhoods. Since congestion is negligible, performance is driven almost entirely by travel distances, leading the solution to emphasize geographic dispersion and unit availability. In this regime, the problem reduces to a deterministic facility location model, unaffected by stochastic server availability.

As demand intensifies to moderate levels (Figures~\ref{fig_lambda_2}–\ref{fig_lambda_6}), the system dynamics begin to shift. Congestion becomes more frequent, and the likelihood of multiple units being busy at the same time increases. As a result, the optimal deployment gradually reallocates units toward central neighborhoods where call intensities are higher. Peripheral tracts lose some dedicated coverage, while central regions gain unit density, striking a balance between travel efficiency and availability. This adjustment unfolds gradually as $\theta$ increases from 2 to 6: average response times rise, spatial coverage becomes less balanced as more units concentrate in high-demand regions, and the resulting configuration reveals how spatial demand variation interacts with stochastic unit availability in shaping placement. 

In the heavy-load regime (Figure~\ref{fig_lambda_inf}), when arrival rates become extremely high and colocation is permitted, the optimal configuration converges to the $1$-Median solution, as stated in Theorem~\ref{thm:asym-heavy-lambda}. In this setting, all units concentrate at the city center. Colocating units effectively pools capacity, allowing the system to minimize the weighted average response time by reducing the expected distance from the central location under severe congestion.

Viewed together, these results reveal how the optimal configuration evolves as system load increases: (i) under light traffic, units are widely dispersed, resembling the $p$-Median solution; (ii) under moderate traffic, units shift toward high-demand centers to balance travel efficiency and congestion risk; and (iii) under heavy traffic, all units converge to a single central hub, reflecting capacity pooling under extreme load. This progression highlights how the model bridges classical deterministic location model and the stochastic dynamics of congested service systems.

\subsubsection{Performance under Increasing Arrival Rate}

In this section, we evaluate performance under varying arrival rates, comparing our approaches and the classical $p$-Median model against the optimal solutions obtained through exhaustive enumeration. 

Table~\ref{diff_lambda} shows that as offered load increases, both SparBL and GP-$p$M consistently obtain the optimal solution across all tested scenarios. In contrast, the $p$-Median does not find the optimal solution as traffic intensifies. At very low loads (e.g., $\lambda=0.1$), the $p$-Median solution is virtually identical to the optimum, consistent with our asymptotic result in Theorem~\ref{thm:asym-vanishing-lambda}. As arrival rates grow, however, the gap between the $p$-Median and the optimal solution widens steadily, from $0.16$ minutes at $\lambda=0.225$ (the original load level) to $0.67$ minutes at $\lambda=1.0$. An average gap of $0.16$ minutes translates into a total annual delay of more than 4,800 minutes, given over 30,000 calls per year. This growing discrepancy highlights that the $p$-Median model underestimates congestion effects in higher-traffic settings. By contrast, our methods remain effective across all load levels, consistently maintaining optimality even as the $p$-Median benchmark deteriorates. 















\begin{table}[!ht]
\caption{Performance comparison at varying arrival rates $\lambda$. The number with an asterisk represents the original arrival rate.}
\centering
\footnotesize
        {\def\arraystretch{1}  
\begin{tabular*}{\textwidth}{@{\extracolsep{\fill}}ccccccc}
\hline
\hline
Offered Load
            & Opt (min) 
            & SparBL 
            & GP-$p$M
            & $p$-Median 
            & \makecell{Gap\\{\scriptsize ($p$-Median\,vs.\,Opt)}} 
            & \makecell{Gap\\{\scriptsize (SparBL/GP-$p$M\,vs.\,Opt)}} \\
\hline
{0.1}       & {5.54} & {5.54} & {5.54} & {5.55} & 0.01 & 0.00 \\
{$0.225^*$} & {5.86} & {5.86} & {5.86} & {6.01} & 0.15 & 0.00 \\
{0.3}       & {6.11} & {6.11} & {6.11} & {6.27} & 0.16 & 0.00 \\
{0.4}       & {6.51} & {6.51} & {6.51} & {6.70} & 0.19 & 0.00 \\
{0.5}       & {6.93} & {6.93} & {6.93} & {7.17} & 0.24 & 0.00 \\
{0.6}       & {7.27} & {7.27} & {7.27} & {7.58} & 0.31 & 0.00 \\
{0.7}       & {7.61} & {7.61} & {7.61} & {8.03} & 0.42 & 0.00 \\
{0.8}       & {7.90} & {7.90} & {7.90} & {8.41} & 0.51 & 0.00 \\
{0.9}       & {8.15} & {8.15} & {8.15} & {8.74} & 0.59 & 0.00 \\
{1.0}       & {8.34} & {8.34} & {8.34} & {9.01} & 0.67 & 0.00 \\
\hline
\hline
\end{tabular*}
\label{diff_lambda}
}
\end{table}












Many cities operate with substantially higher ambulance utilization rates than St. Paul. In such environments, our approaches are expected to deliver even greater improvements. This makes them particularly valuable in high-demand systems, where traditional deterministic models tend to greatly underestimate congestion effects and therefore struggle to produce near-optimal solutions.

Table~\ref{tab_utilization} reports estimated ambulance utilizations across major U.S. cities, based on data from the 2019 \textit{Firehouse Magazine} survey (pre-pandemic, to avoid distortions introduced by COVID-19). 
Cities without survey responses or without publicly available ambulance counts were excluded. To ensure comparability, utilization rates were approximated using St. Paul’s average service rate as a baseline. Actual service times in larger, busier cities are likely even longer, suggesting that true utilizations may be higher than those reported here. 
Taken together, these estimates underscore a key insight: as system load intensifies, our methods maintain optimal performance while the traditional approach diverges, making our solutions especially well-suited for metropolitan EMS systems with high utilization and heavy congestion. 






   
   


               
                
               
               
               
                
               

\begin{table}[!ht]
\caption{Utilization of EMS units in different regions and cities.}
\centering
\footnotesize
        {\def\arraystretch{1}  
\begin{tabular*}{\textwidth}{@{\extracolsep{\fill}}ccccc}
\hline
\hline
City & Yearly Med Calls & Calls per Minute & \# of Ambulances & Utilization\\
			\hline 
			{New York}&{2,128,560}&{4.05}&{450}&{30.6\%}\\

			{Chicago}&{596,807}&{1.14}&{80}&{48.3\%}\\

			{Baltimore}&{183,306}&{0.35}&{27}&{43.9\%}\\

			{Philadelphia}&{272,772}&{0.52}&{57}&{30.9\%}\\

			{Los Angeles}&{414,375}&{0.79}&{148}&{18.1\%}\\

			{Phoenix}&{194,406}&{0.37}&{36}&{34.9\%}\\
   
			{Washington DC}&{173,004}&{0.33}&{38}&{29.5\%}\\
   
               {San Antonio}&{164,458}&{0.31}&{40}&{26.9\%}\\

               {Las Vegas}&{89,728}&{0.17}&{24}&{24.5\%}\\

               {Boston}&{143,189}&{0.27}&{26}&{36.1\%}\\
               
                {Nashville}&{82,221}&{0.16}&{28}&{19.3\%}\\
                
               {Tampa}&{74,634}&{0.14}&{18}&{27.2\%}\\
               
               {St. Louis}&{71,439}&{0.14}&{10}&{46.8\%}\\
               
               {Milwaukee}&{70,461}&{0.13}&{12}&{38.5\%}\\
               
                {Memphis}&{125,144}&{0.24}&{32}&{25.6\%}\\
                
               {Albuquerque}&{96,421}&{0.18}&{20}&{31.6\%}\\
               
   		{St. Paul} &{30,911}&{0.07}&{9}&{22.5\%}\\
\hline
\hline
\end{tabular*}
\label{tab_utilization}
}
\end{table}






   
   


               
                
               
               
               
                
               

London, UK, provides another example where our approach would be particularly valuable due to its exceptionally high ambulance utilization. \citet{bavafa2023distributional} and \citet{bavafa2021variance} studied the London Ambulance Service using a comprehensive dataset of calls spanning more than ten years, focusing on how worker fatigue and experience affect service times. From discussions with the authors, we learned that the average ambulance utilization in their dataset was approximately 84\%, a level far higher than that observed in most U.S. cities. 
In such a heavily loaded system, our model could be applied to identify station locations that minimize average response times and to evaluate the impact of increasing the fleet size on system-wide performance. 
This illustrates the value of our framework not only for optimizing existing deployments but also for supporting strategic planning in high-utilization urban settings.

\section{Conclusions}\label{sec_conclusion}

In this paper, we studied the problem of optimizing server locations in spatial hypercube queueing models, and proposed two Bayesian optimization approaches: (i) a parametric method based on sparse Bayesian linear regression with horseshoe priors (SparBL), and (ii) a nonparametric method using Gaussian processes with a $p$-Median prior (GP-$p$M). We established theoretical lower and upper bounds for the optimal solution, analyzed asymptotic behavior by connecting our model to the $p$-Median and 1-Median problems, and proved that both algorithms achieve sublinear regret.

Our numerical experiments and a problem using real-world data from the St. Paul Fire Department show that both methods consistently outperform the baseline methods, particularly under high unit utilization. While the $p$-Median approach performs well in low-utilization settings, it assumes all units are available, and its performance deteriorates as call volumes increase. In contrast, our framework adapts to these dynamics, with its relative advantage becoming greater as utilization rises, a key feature for deployment in busy urban systems. 

To guide practical implementation, we identified high-utilization cities where our approach offers the most significant benefits and identified London as a city where our work would be especially useful. In our current model, we assume units are dispatched based on a pre-specified preference list. In future work we plan to extend this framework to jointly optimize dispatching policies and location decisions.

\spacingset{1.18}
\bibliographystyle{agsm}
\bibliography{main}

\newpage
\begin{appendices}

\begin{center}
{\large\bf SUPPLEMENTARY MATERIAL of \\
    ``Optimizing Server Locations in Spatial Queues: \\ 
    Parametric and Nonparametric Bayesian Optimization Approaches''}
\end{center}

\section{Illustration of BO in a Simple Continuous Example}
\label{ec:GP_illustration}

Bayesian optimization is an approach that is particularly useful for problems where the objective function is expensive or time-consuming to evaluate, or for which the gradient is not available~\citep{frazier2018tutorial}. It is based on the idea of building a \textit{surrogate model} of the objective function, and using a so-called \textit{acquisition function} based on the surrogate model to determine the next solution to be evaluated. 

\begin{figure}[htbp]
\centering
  \includegraphics[width=15cm]{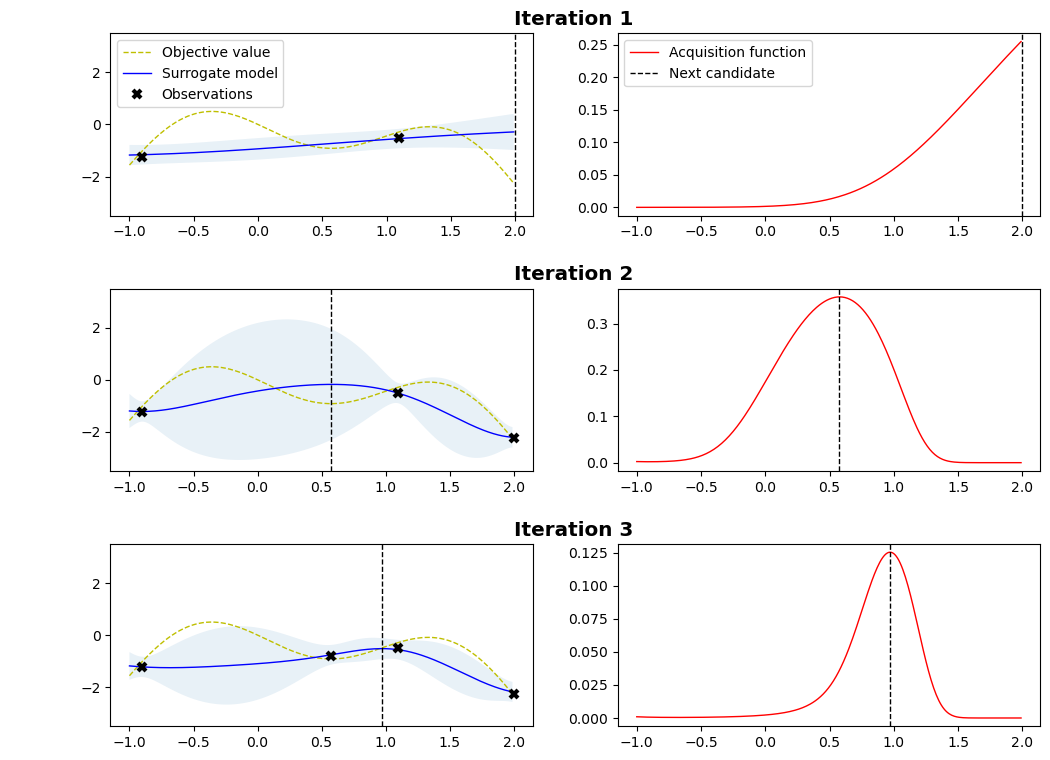}
\caption{A simple illustration of the Bayesian optimization framework in a one-dimensional continuous space.}
\label{fig_BO_illustration}
\vspace{-0.4cm}
\end{figure}

Figure \ref{fig_BO_illustration} illustrates the basic idea of BO using a simple one-dimensional example, where solution points are evaluated and added sequentially at each iteration. This illustration optimizes in a continuous space, in contrast to optimizing in a combinatorial space, which is more challenging. 

The left-hand side of the figure shows three subplots that illustrate the evolution of the surrogate model as solutions are observed and added to the model. The locations of the observed solutions are indicated by the $\times$ symbol. In the first iteration, we randomly select two solutions that are added to the model. In each subsequent iteration, the algorithm selects the solution with the highest acquisition function value for evaluation, shown by the dashed vertical lines. The plots show how the surrogate model is updated based on the newly observed solutions. The right-hand side of the figure shows three subplots that demonstrate the role of the acquisition function. The acquisition function provides a  trade-off between exploration (sampling new solutions with high predicted variances) and exploitation (sampling solutions with high predicted means). Typically, the points with high mean values and high variances have higher values for the acquisition functions, as shown in the three subplots on the right.

\section{Hypercube queueing Model}\label{sec:hypercube}
\subsection{Exact Hypercube Formulation}
\label{ssec:exact_Hypercube}

Exact Hypercube Formulation is as follows:
\begin{subequations}\label{eq_all}
\setlength{\abovedisplayskip}{3pt}
\setlength{\belowdisplayskip}{3pt}
\begin{align}
\text{OPT}_{H} = \min_{\boldsymbol{x}} & \sum_{i \in I} \sum_{j \in J}(\tau_i + t_{i j}) q_{i j}(\boldsymbol{x}) \qquad \text{\emph{(\textbf{Exact $p$-Hypercube})}} \\
\text { s.t. } & \sum_{i \in I} x_{i}=p, \label{eq_feasibility}\\
& q_{i j}(\boldsymbol{x})=  \frac{\sum_{\boldsymbol{b} \in S_{ij}(\boldsymbol{x})} \lambda_j P(\boldsymbol{b})}{\sum_{j\in J} \lambda_j \left(1-\Pi\right)}, \forall i \in I, \forall j \in J, \label{eq_qij}\\
& x_{i} \in\{0,1\}, \forall i \in I.
\end{align}
\end{subequations} 

In the above formulation, $q_{i j}(\boldsymbol{x})$ is the conditional share among calls from $j$ that are served by the unit located at $i$, given a unit configuration $\boldsymbol{x}$.  The state of all $p$ units is represented by $\boldsymbol{b}=\{b_1, \cdots, b_p\}\in \{0,1\}^p$, where $b_k=0$ indicates unit $k$ is available, and $b_k=1$ means unit $k$ is busy. The probability of being in state $\boldsymbol{b}$ is $P(\boldsymbol{b})$, and the blocking probability 
is $\Pi$.
Furthermore, $S_{ij}(\boldsymbol{x})$ denotes the set of states where the unit located at location $i$ is assigned to calls from subregion $j$. 

We note that $q_{ij}(\boldsymbol{x})$ depends on $P(\boldsymbol{b})$ and $\Pi$, which requires solving the balance equations of the spatial hypercube model of \cite{larson74hypercube}. For all $\boldsymbol{b}$, the balance equations are:
\begin{equation}
P(\boldsymbol{b})\Bigg[\underbrace{\sum_{j\in J}\lambda_{j}}_{\text{arrival}}+\underbrace{\sum_{k=1}^p b_k \mu_k}_{\text{service completion}} \Bigg]=\underbrace{\sum_{ \boldsymbol{b}': H(\boldsymbol{b}, \boldsymbol{b}')^+ = 1} P(\boldsymbol{b}') r(\boldsymbol{b}, \boldsymbol{b}')}_{\text{upward transitions}} + \underbrace{\sum_{\boldsymbol{b}': H(\boldsymbol{b}, \boldsymbol{b}')^- = 1} P(\boldsymbol{b}')\mu_{i(\boldsymbol{b}, \boldsymbol{b}')}}_{\text{downward transitions}}, \label{eq_hypercube}
\end{equation}
where $r(\boldsymbol{b}, \boldsymbol{b}')$ is the transition rate from state $\boldsymbol{b}$ to state $\boldsymbol{b}'$ in the spatial hypercube model, $i(\boldsymbol{b}, \boldsymbol{b}')$ is the index at which $\boldsymbol{b}$ and $\boldsymbol{b}'$ differ, and $H(\boldsymbol{b}, \boldsymbol{b}')$ is the Hamming distance from $\boldsymbol{b}$ to $\boldsymbol{b}'$, which is defined as
\begin{equation}
    H(\boldsymbol{b}, \boldsymbol{b}') = w\left(\left[\boldsymbol{b} \cap \neg \boldsymbol{b}'\right] \cup \left[\neg \boldsymbol{b} \cap \boldsymbol{b}'\right]\right), \label{eq_hamming}
\end{equation}
where $\cap$ and $\cup$ are bitwise AND and OR operations, respectively, and $\neg \boldsymbol{b}$ is the logical NOT operation of the binary number representation of state $\boldsymbol{b}$. The function $w(\cdot)$ counts the number of 1s in the binary state representation. Transitions in the spatial hypercube model only take place to adjacent states, i.e., $H(\boldsymbol{b}, \boldsymbol{b}')=1$. \cite{larson74hypercube} introduced upward and downward Hamming distances, which are denoted by $H(\boldsymbol{b}, \boldsymbol{b}')^+$ and $H(\boldsymbol{b}, \boldsymbol{b}')^-$, respectively. An upward distance measures the distance to states reached by call arrivals and a downward distance measures the distance to states reached by service completions. Thus,
\begin{equation}
    H(\boldsymbol{b}, \boldsymbol{b}')^+ = w\left(\left[\boldsymbol{b} \cap \neg \boldsymbol{b}'\right]\right), \qquad H(\boldsymbol{b}, \boldsymbol{b}')^- = w\left(\left[\neg \boldsymbol{b} \cap \boldsymbol{b}'\right]\right).
\end{equation}
Clearly, $H(\boldsymbol{b}, \boldsymbol{b}') = H(\boldsymbol{b}, \boldsymbol{b}')^+ + H(\boldsymbol{b}, \boldsymbol{b}')^-$. We note that in the exact hypercube formulation, service times must be assumed exponential with rate $\mu_i$ for unit $i$ in order to express the balance equations. This restriction contrasts with the approximation algorithm in the main body of the paper, which accommodates general service time distributions and allows the service rate $\mu_{ij}$ to vary by both unit location and subregion. In addition, evaluating the exact objective function is computationally intensive, as it requires solving the transition matrix of a Markov chain with $2^p$ states, where $p$ is the number of units.

\begin{theorem}
\label{thm:NP-hard_exact}
The Exact $p$-Hypercube problem is NP-hard.
\end{theorem}


\proof{Proof.}

We prove by a polynomial-time reduction from the $p$-Median instance. Let $\boldsymbol{b}\in\{0,1\}^p$ denote a server state vector in the exact hypercube notation (where $b_k=1$ means server $k$ is busy). 
By setting service rates $\mu_{ij}\to\infty$, the following hold:
\begin{itemize}
  \item For any state $\boldsymbol{b}$ with at least one busy server, $P(\boldsymbol{b})\to 0$, and the probability mass concentrates on the all-idle state $\boldsymbol{b}=\boldsymbol{0}$ (i.e. $P(\boldsymbol{0})\to 1$).
  \item Consequently, the blocking probabilities satisfy $\Pi \to 0$.
\end{itemize}
Substitute these limits into the exact allocation formula:
\[
q_{ij}(\boldsymbol{x}) \;=\; \frac{\sum_{\boldsymbol{b}\in S_{ij}(\boldsymbol{x})}\lambda_j P(\boldsymbol{b})}
                 {\sum_{j'\in\mathcal J}\lambda_{j'}\bigl(1-\Pi\bigr)}.
\]
In the $\mu_{ij}\to\infty$ limit the numerator is nonzero only for those states in which the (unique) server assigned to serve calls from $j$ is the facility at $i$ and that server is idle; because all busy-state probabilities vanish, this simplifies to the indicator that $i$ is the nearest open facility (in $\boldsymbol{x}$) to $j$. Hence
\[
q_{ij}(\boldsymbol{x}) \;=\; \frac{\lambda_j}{\sum_{j'\in\mathcal J}\lambda_{j'}}\cdot \boldsymbol{1}\{\,i=\arg\min_{k: x_k=1} t_{kj}\,\}.
\]
Intuitively, when servers are never busy, each arriving call is immediately served by the closest open unit according to the preference order.

The exact $p$-Hypercube objective is
\[
\min_{\boldsymbol{x}\in\{0,1\}^N}\;\sum_{i\in\mathcal I}\sum_{j\in\mathcal J}(\tau_i+t_{ij})\,q_{ij}(\boldsymbol{x})
\quad\text{s.t.}\quad \sum_{i\in\mathcal I} x_i = p.
\]
Substituting $\tau_i=0$ the simplified $q_{ij}(\boldsymbol{x})$ above yields
\[
\min_{\boldsymbol{x}}\sum_{i\in\mathcal I}\sum_{j\in\mathcal J} t_{ij}\cdot
\frac{\lambda_j}{\sum_{j'\in\mathcal J}\lambda_{j'}}\,\boldsymbol{1}\{\,i=\arg\min_{k: x_k=1} t_{kj}\,\}
= \frac{1}{\sum_{j'\in\mathcal J}\lambda_{j'}}\min_{\boldsymbol{x}}\sum_{j\in\mathcal J} \lambda_j \min_{i: x_i=1} t_{ij}.
\]
Because $\sum_{j'\in\mathcal J}\lambda_{j'}$ is a positive constant independent of $\boldsymbol{x}$, minimizing the exact $p$-Hypercube objective under this construction is equivalent to solving the original $p$-Median problem.
\Halmos
\endproof

We next show that a special case of $p$-Median provides a lower bound for the optimal value of $p$-Hypercube.

\begin{theorem}
\label{thm_low_bound_exact}
The optimal value to the $p$-Median problem is a lower bound for the optimal value of the exact $p$-Hypercube problem when $w_j = \lambda_j/\sum_{j'\in \mathcal{J}} \lambda_{j'}$, i.e.,
\begin{equation*}
\text{OPT}_H \geq \text{OPT}_M.
\end{equation*}
\end{theorem}

\proof{Proof.}

Let \(\boldsymbol{x}\in\{0,1\}^N\) be any feasible facility placement. For each subregion \(j\) and facility \(i\) denote by \(q_{ij}(\boldsymbol{x})\ge0\) the proportion of calls from \(j\) handled by unit \(i\) under placement \(\boldsymbol{x}\). Assume \(q_{ij}(\boldsymbol{x})=0\) whenever \(x_i=0\) (no unit at \(i\) implies zero allocation to \(i\)). Set
\[
S_j(\boldsymbol{x})\coloneqq\sum_{i\in\mathcal{I}} \frac{1}{w_j}q_{ij}(\boldsymbol{x}),\qquad c_{ij}\coloneqq \tau_i+t_{ij},
\]
and define the objective at \(\boldsymbol{x}\)
\[
\text{OPT}_{H}(\boldsymbol{x})=\sum_{j\in\mathcal{J}}\sum_{i\in\mathcal{I}} c_{ij}q_{ij}(\boldsymbol{x})=\sum_{j\in\mathcal{J}}\sum_{i\in\mathcal{I}} w_j\,c_{ij}\,\frac{1}{w_j}q_{ij}(\boldsymbol{x}).
\]
For the fixed-\(\boldsymbol{x}\) $p$-Median problem write
\[
\text{OPT}_M=\min_{\substack{y_{ij}\ge0\\ \sum_i y_{ij}=1,\ y_{ij}\le x_i}}
\sum_{j\in\mathcal{J}}\sum_{i\in\mathcal{I}} w_j\,c_{ij}\,y_{ij}.
\]

Fix a subregion \(j\). If \(S_j(\boldsymbol{x})=0\) then \(q_{ij}(\boldsymbol{x})=0\) for all \(i\) and the contribution of \(j\) to both sides is zero, so the inequality holds trivially. Suppose \(S_j(\boldsymbol{x})>0\) and define the normalized vector
\[
\hat y_{ij}(\boldsymbol{x}):=\frac{q_{ij}(\boldsymbol{x})}{w_j S_j(\boldsymbol{x})}\quad \forall\ i\in\mathcal{I}.
\]
By construction \(\hat y_{ij}(\boldsymbol{x})\ge0\), \(\sum_i\hat y_{ij}(\boldsymbol{x})=1\), and \(\hat y_{ij}(\boldsymbol{x})=0\) whenever \(x_i=0\); hence \(\hat y_{\cdot j}(\boldsymbol{x})\) is feasible for the $p$-Median subproblem at fixed \(\boldsymbol{x}\). Let \(\text{OPT}_{j,M}\) denote the contribution of subregion \(j\) to \(\text{OPT}_M\). Then
\[
\sum_{i\in\mathcal{I}} w_j c_{ij}\,\hat y_{ij}(\boldsymbol{x}) \ge \text{OPT}_{j,M}.
\]
Multiplying by \(S_j(\boldsymbol{x})\) gives
\begin{equation}\label{eq:lower-bound-1}
    \sum_{i\in\mathcal{I}} w_j c_{ij}\frac{1}{w_j}q_{ij}(\boldsymbol{x})
= S_j(\boldsymbol{x})\sum_{i\in\mathcal{I}} w_j c_{ij}\,\hat y_{ij}(\boldsymbol{x})
\ge S_j(\boldsymbol{x})\,\text{OPT}_{j,M}.
\end{equation}

In the exact $p$-Hypercube formulation, the allocation probabilities are derived from the true joint availability distribution of units, and by construction satisfy, for every feasible \(\boldsymbol{x}\) and every \(j\),
\begin{equation}\label{eq:lower-bound-2}
    S_j(\boldsymbol{x})\coloneqq\sum_{i\in\mathcal{I}} \frac{1}{w_j}q_{ij}(\boldsymbol{x})=1.
\end{equation}
Substituting~\eqref{eq:lower-bound-2} into~\eqref{eq:lower-bound-1} yields the claimed lower bound.
\Halmos
\endproof

The same upper bound also applies following the argument in Lemma~\ref{lemma_upper_bound}.

\subsection{Approximation Result}
\label{ssec:approx_accuracy}
\begin{figure}[hbtp]
    \centering
    \begin{subfigure}[t]{0.44\textwidth}
        \centering
        \includegraphics[width=\linewidth]{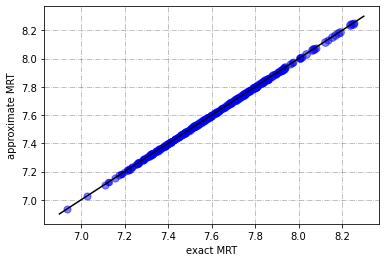}
        \caption{Exact vs. Approximate}
    \end{subfigure}
    \begin{subfigure}[t]{0.48\textwidth}
        \centering
        \includegraphics[width=\linewidth]{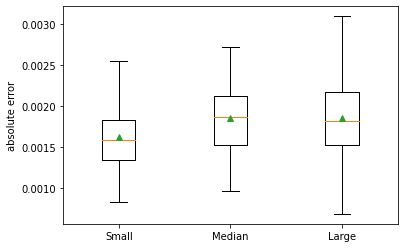}
        \caption{Absolute errors for various system sizes}
    \end{subfigure}%
    \caption{The accuracy of the approximate hypercube model.}
    \label{fig_larson}
\end{figure}

In Figure \ref{fig_larson}, we show that the approximation is very close to the exact value, with a mean absolute error of less than 0.002 minutes, obtained by generating 100 random setups for small (15 units), medium (20 units), and large problems (30 units). For a system with 30 units, the approximation algorithm takes less than 0.01 seconds while the exact solution requires solving $2^{30}$ simultaneous equations which takes hours.

\section{Posterior Distribution of Horseshoe Prior}
\label{ec:hs_posteiror}
To enable efficient sampling under the horseshoe prior, the half-Cauchy distributions on $\beta_k$ and $\tau$ are reparameterized using auxiliary variables $\nu_k$ and $\xi$, which yield inverse-gamma representations~\citep{makalic2016inverseGamma}. This reparameterization transforms the hierarchical model into a form where the \textit{posterior distribution} can be expressed in closed form. Under this reparameterization, the full conditional posterior distributions for the model parameters are given by: 
\begin{equation}\label{eqn:horseshoe-posterior}
    \begin{aligned}
    \boldsymbol{\alpha} &\mid \cdot \sim \mathcal{N}(\boldsymbol{M}^{-1} \boldsymbol{X}^T \boldsymbol{y}, \sigma^2 \boldsymbol{M}^{-1}), \quad\boldsymbol{M} = (\boldsymbol{X}^T \boldsymbol{X} + \Sigma_*^{-1}), \quad \Sigma_* = \tau^2 \text{diag}(\beta_1^2, \ldots, \beta_D^2)\\
    \sigma^2 &\mid \cdot \sim \mathcal{IG} \left( \frac{N+D}{2}, \frac{(\boldsymbol{y} - \boldsymbol{X}\boldsymbol{\alpha})^T (\boldsymbol{y} - \boldsymbol{X}\boldsymbol{\alpha}) + \boldsymbol{\alpha}^T \Sigma_*^{-1} \boldsymbol{\alpha}}{2} \right)\\
    \beta_k^2 & \mid \cdot \sim \mathcal{IG} \left( 1, \frac{1}{\nu_k} + \frac{\alpha_k^2}{2 \tau^2 \sigma^2} \right) \quad k = 1, \ldots, D\\
    \tau^2 &\mid \cdot \sim \mathcal{IG} \left( \frac{D+1}{2}, \frac{1}{\xi} + \frac{1}{2 \sigma^2} \sum_{k=1}^{D} \frac{\alpha_k^2}{\beta_k^2} \right)\\
    \nu_k &\mid \cdot \sim \mathcal{IG} \left( 1, 1 + \frac{1}{\beta_k^2} \right) \quad k = 1, \ldots, D\\
    \xi &\mid \cdot \sim \mathcal{IG} \left( 1, 1 + \frac{1}{\tau^2} \right).
\end{aligned}
\end{equation}

Here, $\mathcal{IG}(a, b)$ denotes the inverse-gamma distribution with shape parameter $a$ and scale parameter $b$. 
The matrix $\boldsymbol{M}$ captures the posterior precision of the regression coefficients $\boldsymbol{\alpha}$, incorporating both the design matrix $\boldsymbol{X}$ and the shrinkage structure imposed by $\Sigma_*$. The auxiliary variables $\nu_k$ and $\xi$ allow the half-Cauchy priors on $\beta_k$ and $\tau$ to be expressed as hierarchical inverse-gamma priors, which are conjugate to the Gaussian likelihood. The result is a computationally tractable yet flexible model that automatically performs shrinkage and variable selection.
\section{Miscellaneous Proofs}

\subsection{Proof of Theorem \ref{thm:NP-hard_approx}}
\label{ssec:proof_NPhard_approx}
\proof{Proof.}
We prove the claim by a polynomial-time reduction from the \emph{$p$-Median problem}, which is known to be NP-hard for a general problem. Set $c_{ij}\coloneqq \tau_i+t_{ij}$. An instance of the $p$-Median problem is given by a set of candidate facility locations $\mathcal{I}$ and a set of demand points $\mathcal{J}$; demand weights $\{w_j : j \in \mathcal{J}\}$, a travel time matrix $\{c_{ij} : i \in \mathcal{I}, j \in \mathcal{J}\}$, and an integer $p$ representing the number of facilities to open.
The objective is to choose $p$ locations to minimize
\[
    \sum_{j \in \mathcal{J}} w_j \min_{i: x_i=1} c_{ij}.
\]

Given any $p$-Median instance $(\mathcal{I}, \mathcal{J}, w_j, c_{ij}, p)$, we construct a $p$-Hypercube instance as follows:
\begin{itemize}
    \item Set the number of candidate sites to $N = |\mathcal{I}|$ and the number of subregions to $M = |\mathcal{J}|$.
    \item Assign the travel time parameters as $c_{ij} = c_{ij}$ and set the turnout times $\tau_i = 0$ for all $i \in \mathcal{I}$.
    \item Choose service rates $\mu_{ij} \to \infty$, which implies that the server utilization $\rho_i = 0$ for all $i$.
    \item Set arrival rates $\lambda_j = w_j$ and define the preference list $(\gamma_{jl})$ for subregion $j$ by sorting candidate sites in non-decreasing order of $c_{ij}$.
    \item Since $\rho_i = 0$, the correction factor $Q(N,\bar{\rho}, \cdot)$ evaluates to $1$ under this construction.
\end{itemize}

Under these parameter settings, the $p$-Hypercube allocation formula simplifies to
\[
    q_{ij}(\boldsymbol{x}) 
    = \mathbf{1}\left\{ i = \arg\min_{k: x_k=1} c_{kj} \right\}.
\]

That is, each request is always served by its closest open facility. The $p$-Hypercube objective becomes:
\[
    \min_{\boldsymbol{x}} \sum_{i \in \mathcal{I}} \sum_{j \in \mathcal{J}} 
        c_{ij}\frac{\lambda_j\cdot q_{i j}(\boldsymbol{x})}{\sum_{j\prime\in \mathcal{J}} \lambda_{j\prime}}
    = \frac{1}{\sum_{j'\in\mathcal J}\lambda_{j'}} \min_{\boldsymbol{x}} \sum_{j \in \mathcal{J}} \lambda_j \min_{i: x_i=1} c_{ij}.
\]
Since $\sum_{j'\in\mathcal J}\lambda_{j'}$ is a constant independent of $\boldsymbol{x}$, this is exactly the $p$-Median objective.
\Halmos
\endproof


\subsection{Proof of Theorem \ref{thm_low_bound_approx}}
\label{ssec_lower_bound_proof_approx}
\proof{Proof.}

Let \(\boldsymbol{x}\in\{0,1\}^N\) be any feasible facility placement. For each subregion \(j\) and facility \(i\) denote by \(q_{ij}(\boldsymbol{x})\ge0\) the proportion of calls from \(j\) handled by unit \(i\) under placement \(\boldsymbol{x}\). Let \(q_{ij}(\boldsymbol{x})=0\) when \(x_i=0\) (no unit at \(i\) implies zero allocation to \(i\)). Set
\[
S_j(\boldsymbol{x})\coloneqq\sum_{i\in\mathcal{I}} \frac{1}{w_j}q_{ij}(\boldsymbol{x}),\qquad c_{ij}\coloneqq \tau_i+t_{ij},\qquad w_j\coloneqq  \frac{\lambda_j}{\sum_{j\in\mathcal J}\lambda_j},
\]
and define the objective at \(\boldsymbol{x}\)
\[
\text{OPT}_{H}(\boldsymbol{x})=\sum_{j\in\mathcal{J}}\sum_{i\in\mathcal{I}} w_j\,c_{ij}\,\frac{q_{ij}(\boldsymbol{x})}{w_j}.
\]
For the fixed-\(\boldsymbol{x}\) $p$-Median problem write
\[
\text{OPT}_M=\min_{\substack{y_{ij}\ge0\\ \sum_i y_{ij}=1,\ y_{ij}\le x_i}}
\sum_{j\in\mathcal{J}}\sum_{i\in\mathcal{I}} w_j\,c_{ij}\,y_{ij}.
\]

Fix a subregion \(j\). If \(S_j(\boldsymbol{x})=0\) then \(q_{ij}(\boldsymbol{x})=0\) for all \(i\) and the contribution of \(j\) to both sides is zero, so the inequality holds trivially. Suppose \(S_j(\boldsymbol{x})>0\) and define the normalized vector
\[
\hat y_{ij}(\boldsymbol{x}):=\frac{q_{ij}(\boldsymbol{x})}{w_j S_j(\boldsymbol{x})}\quad \forall\ i\in\mathcal{I}.
\]
By construction \(\hat y_{ij}(\boldsymbol{x})\ge0\), \(\sum_{i\in\mathcal{I}} \hat y_{ij}(\boldsymbol{x})=1\), and \(\hat y_{ij}(\boldsymbol{x})=0\) whenever \(x_i=0\); hence \(\hat y_{\cdot j}(\boldsymbol{x})\) is feasible for the $p$-Median subproblem at fixed \(\boldsymbol{x}\). Let \(\text{OPT}_{j,M}\) denote the contribution of subregion \(j\) to \(\text{OPT}_M\). Then
\[
S_j\sum_{i\in\mathcal{I}} w_j c_{ij}\,\hat y_{ij}(\boldsymbol{x}) \ge S_j\text{OPT}_{j,M}.
\]
Summing over \(j\in\mathcal{J}\) and by $ S_j(\boldsymbol{x})=1$ yields \(\text{OPT}_{H}(\boldsymbol{x})\ge \text{OPT}_M\). Minimizing both sides over feasible \(\boldsymbol{x}\) gives a general global inequality
\[
\text{OPT}_{H} \ge \mathrm{OPT}_M.
\]

\Halmos
\endproof

\subsection{Proof of Theorem \ref{thm:asym-vanishing-lambda}}
\label{ssec:asym-vanishing-lambda}
\proof{Proof.}

Define
\[
\text{OPT}_H^{0}(\boldsymbol{x}) \coloneqq \frac{1}{\sum_{j'\in\mathcal J}\lambda_{j'}} 
\sum_{j\in\mathcal J}\lambda_j \min_{i: x_i=1}(\tau_i+t_{ij}).
\]
We first show that $\text{OPT}_H(\theta,\cdot)$ converges uniformly to $\text{OPT}_H^{0}(\cdot)$ as $\theta\to 0$. Since the feasible set $\mathcal D$ is finite, so it suffices to show point-wise convergence $\text{OPT}_H(\theta,\boldsymbol{x})\to \text{OPT}_H^0(\boldsymbol{x})$ for each fixed $\boldsymbol{x}$, together with a uniform modulus bound that yields uniform convergence on $\mathcal D$.

From~\eqref{Eq:rhoi}, we have for each unit $i$ and each feasible placement $\boldsymbol{x}$
\[
\rho_i(\theta;\boldsymbol{x})\;=\;1 - \Biggl(1 + \sum_{k=1}^N\sum_{j\in G_{ik}(\boldsymbol{x})}
\frac{\lambda_j(\theta)\,Q\bigl(N,\bar\rho(\theta;\boldsymbol{x}),k-1\bigr)\prod_{l=1}^{k-1}\rho_{\gamma_{jl}}(\theta;\boldsymbol{x})}{\mu_{ij}}\Biggr)^{-1}.
\]
Use the inequality $1-(1+z)^{-1}\le z$ for all $z\ge 0$ to obtain
\[
\rho_i(\theta;\boldsymbol{x})\;\le\; \sum_{k=1}^N\sum_{j\in G_{ik}(\boldsymbol{x})}
\frac{\lambda_j(\theta)\,Q\bigl(N,\bar\rho(\theta;\boldsymbol{x}),k-1\bigr)\prod_{l=1}^{k-1}\rho_{\gamma_{jl}}(\theta;\boldsymbol{x})}{\mu_{ij}},
\]
where $\bar\rho(\theta;\boldsymbol{x})$ is the average server utilization. 

By continuity of $\bar\rho\mapsto Q(N,\bar\rho,r)$, there exists $\delta>0$ and $C_Q<\infty$ such that for all $\bar\rho\in[0,\delta]$ and all relevant $r$ we have $Q(N,\bar\rho,r)\le C_Q$. 
Because $\rho_i(0;\boldsymbol{x})=0$ (when $\theta=0$ no demand exists), continuity of the selected solution $\rho(\theta;\boldsymbol{x})$ in $\theta$ implies that for sufficiently small $\theta$ we have $\bar\rho(\theta;\boldsymbol{x})\le\delta$ and hence the bound $Q(N,\bar\rho,r)\le C_Q$ applies. Substituting $\lambda_j(\theta)=\theta\lambda_j$ gives
\[
\rho_i(\theta;\boldsymbol{x})\le \theta\; \sum_{k=1}^N\sum_{j\in G_{ik}(\boldsymbol{x})} \frac{\lambda_j\,C_Q}{\mu_{ij}} .
\]
The right-hand side is independent of $\theta$ except for the factor $\theta$, so there exists a finite constant $C(\boldsymbol{x})$ (depending on $\boldsymbol{x}$ but not on $\theta$) such that for all sufficiently small $\theta$,
\[
0\le\rho_i(\theta;\boldsymbol{x})\le C(\boldsymbol{x})\,\theta.
\]
Since $\mathcal D$ is finite, we may take $C:=\max_{x\in\mathcal X} C(\boldsymbol{x})$ and conclude the uniform (in $\boldsymbol{x}$ and $i$) bound
\begin{equation}\label{eqn:uniform-bound}
    \max_{x\in\mathcal X}\max_{i\in\mathcal I}\rho_i(\theta;\boldsymbol{x})\le C\theta\qquad\text{for all sufficiently small }\theta.
\end{equation}

In particular $\rho_i(\theta;\boldsymbol{x})\to 0$ uniformly in $\boldsymbol{x}$ as $\theta\rightarrow 0$.

Fix $\boldsymbol{x}\in\mathcal X$ and $i,j$. Recall
\[
q_{ij}(\theta;\boldsymbol{x})=
Q\bigl(N,\bar\rho(\theta;\boldsymbol{x}),\eta_{ij}-1\bigr)\,
\Bigl(\prod_{l=1}^{\eta_{ij}-1}\rho_{\gamma_{jl}}(\theta;\boldsymbol{x})\Bigr)\,
\bigl(1-\rho_i(\theta;\boldsymbol{x})\bigr)
\]
Let $\eta=\eta_{ij}$. There are two cases:

\emph{Case $\eta=1$.} This corresponds to the case where location $i$ is the first preferred facility for subregion $j$. The product over an empty index set equals $1$, i.e., $\prod_{l=1}^{\eta_{ij}-1}\rho_{\gamma_{jl}}(\theta;\boldsymbol{x}) = 1$ . Using continuity of $Q$ at $\bar\rho=0$ and $1-\rho_i(\theta;\boldsymbol{x})\to 1$, we obtain
\[
\lim_{\theta\rightarrow 0}  \frac{\lambda_j}{\sum_{j\in\mathcal J}\lambda_j}q_{ij}(\theta;\boldsymbol{x}) \;=\; \frac{\lambda_j}{\sum_{j\in\mathcal J}\lambda_j}\cdot Q(N,0,0)\cdot 1 \cdot 1
\;=\; \frac{\lambda_j}{\sum_{j\in\mathcal J}\lambda_j},
\]
since $Q(N,0,0)=1$.

\emph{Case $\eta\ge 2$.} This corresponds to the case where at least one strictly preferred facility ahead of $i$. 
Using the uniform bound $\rho_\cdot(\theta;\boldsymbol{x})\le C\theta$ we get
\[
0\le \prod_{l=1}^{\eta-1}\rho_{\gamma_{jl}}(\theta;\boldsymbol{x})\le (C\theta)^{\eta-1}\xrightarrow[\theta\rightarrow 0]{}0.
\]
Hence $q_{ij}(\theta;\boldsymbol{x})\to 0$ in this case.

Combining the two cases, for each fixed $\boldsymbol{x}$ and $j$ the limiting allocation assigns the whole fraction $\lambda_j/\sum_{j'\in\mathcal J}\lambda_j'$ to the first open (preferred) facility in $j$'s preference list. Equivalently, if we denote by $i_j(\boldsymbol{x})$ the highest-preference open facility to $j$ under $\boldsymbol{x}$, then
\[
\lim_{\theta\rightarrow 0}  \frac{\lambda_j}{\sum_{j\in\mathcal J}\lambda_j}q_{ij}(\theta;\boldsymbol{x}) = \frac{\lambda_j}{\sum_{j'\in\mathcal J}\lambda_j'}\boldsymbol{1}\{i=i_j(\boldsymbol{x})\}.
\]

Define $q^0_{ij}(\boldsymbol{x})=\lim_{\theta\rightarrow 0}q_{ij}(\theta;\boldsymbol{x})$
By the above convergence argument, the pointwise convergence $q_{ij}(\theta;\boldsymbol{x})\to q^0_{ij}(\boldsymbol{x})$ implies
\[
\lim_{\theta\rightarrow 0} \text{OPT}_H(\theta,\boldsymbol{x})
= \sum_{i,j}(\tau_i+t_{ij}) \frac{\lambda_j}{\sum_{j\in\mathcal J}\lambda_j}q^0_{ij}(\boldsymbol{x})
= \frac{1}{\sum_{j\in\mathcal J}\lambda_j}\sum_{j\in\mathcal J}\lambda_j\min_{i:x_i=1}(\tau_i+t_{ij})=:\text{OPT}_H^0(\boldsymbol{x}).
\]
Moreover, the convergence is uniform in $\boldsymbol{x}\in\mathcal X$ because we obtained uniform bounds in~\eqref{eqn:uniform-bound} and $\mathcal D$ is finite. Uniform convergence on a finite domain implies that minimizers of $\text{OPT}_H(\theta,\cdot)$ converge (in the sense that any limit point of minimizers is a minimizer of the limit function) to minimizers of $\text{OPT}_H^0(\cdot)$. 


Next, we show that when $\theta$ is small, the minimizers of $\operatorname{OPT}_H(\theta, \cdot)$ are the same as the minimizers of $\mathrm{OPT}_H^0$. Suppose $\boldsymbol{x}^{0}$ is an optimal solution to $\text{OPT}_H^{0}$. If it is unique, there exists a gap $\delta>0$ such that 
\[
\text{OPT}_H^{0}(\boldsymbol{x})-\text{OPT}_H^{0}(\boldsymbol{x}^{0}) \;\ge\;\delta,\qquad \forall \boldsymbol{x}\in \mathcal X,\, \boldsymbol{x}\neq \boldsymbol{x}^{0}.
\]
By uniform convergence, there exists $\theta_0>0$ such that 
\[
\sup_{\boldsymbol{x}\in \mathcal X}\big|\text{OPT}_{H}(\theta,\boldsymbol{x})-\text{OPT}_H^{0}(\boldsymbol{x})\big| < \delta/2,\qquad \forall 0<\theta<\theta_0.
\]
That means, for every $\boldsymbol{x}$, including $\boldsymbol{x}^0$, we have
$$
-\frac{\delta}{2}<\mathrm{OPT}_H(\theta, \boldsymbol{x})-\mathrm{OPT}_H^0(\boldsymbol{x})<\frac{\delta}{2}
$$
Thus, for any $\boldsymbol{x}\neq \boldsymbol{x}^{0}$ and $0<\theta<\theta_0$, we have 
\begin{align*}
    \operatorname{OPT}_H(\theta, x)-\operatorname{OPT}_H(\theta, x^0) =\big(\operatorname{OPT}_H^0(x)-\operatorname{OPT}_H^0(x^0)\big)  &+ \big(\operatorname{OPT}_H(\theta, x)-\operatorname{OPT}_H^0(x)\big)\\
 &- \big(\operatorname{OPT}_H(\theta, x^0)-\operatorname{OPT}_H^0(x^0)\big). 
\end{align*}

Therefore,
\[
\text{OPT}_{H}(\theta,\boldsymbol{x})- \text{OPT}_{H}(\theta,\boldsymbol{x}^{0}) 
\;\ge\; \big(\text{OPT}_H^{0}(\boldsymbol{x})-\text{OPT}_H^{0}(\boldsymbol{x}^{0})\big) - \delta \;>\; 0,
\]
implying that $\boldsymbol{x}^{0}$ is also the unique minimizer of $\text{OPT}_H(\theta,\cdot)$.  

If $\text{OPT}_H^{0}$ admits multiple optimal solutions, then the same argument shows that for sufficiently small $\theta>0$, every minimizer of $\text{OPT}_H(\theta,\cdot)$ must belong to the set of minimizers of $\text{OPT}_H^{0}$.

Recall that in~\eqref{eq:objective}, the $p$-Median objective is given by $\sum_{j \in \mathcal{J}} w_j \min_{i:x_i=1}(\tau_i + t_{i j})$. Choosing $w_j=\frac{\lambda_j}{\sum_{j'\in\mathcal{J}}\lambda_j'}$ yields $\text{OPT}_M=\text{OPT}_H^{0}$.
Therefore, there exists $\theta_0>0$ such that for all $0<\theta<\theta_0$, the optimal solution(s) to $\text{OPT}_H(\theta,\cdot)$ coincide with those of the weighted $p$-Median problem. 
\Halmos
\endproof

\subsection{Proof of Theorem \ref{thm:asym-heavy-lambda}}
\label{ssec:asym-heavy-lambda}
\proof{Proof.}

Consider the case where colocating multiple servers at the same site is allowed. Write $L(\boldsymbol x)$ for the multiset of sites that host the $p$ units (so $|L(\boldsymbol x)|=p$). Let
\[
\pi_j:=\lambda_j/{\sum_{j\in \mathcal{J}}\lambda_j}
\]
be the normalized arrival mix; note that $\{\pi_j\}$ do not depend on $\theta$. For a server located at site $i$, denote by $1/\mu_{ij}$ the mean service time for a job from subregion $j$ and let
\[
m_i:=\sum_{j\in\mathcal J}\pi_j\frac{1}{\mu_{ij}}
\]
be the mean service time under the arrival mix seen by that server.

Define the weighted distance from site $i$ to demand as
$$
V_i:=\sum_{j \in \mathcal{J}} \lambda_j (t_{i j}+\tau_i)
$$
The 1-Median solution is $i^* \in \arg \min _{i \in \mathcal{I}} V_i$.

Consider the long-run behavior as $\theta\to\infty$. Over a long horizon, the fraction of all service completions attributed to a particular server equals its long-run completion rate divided by the total completion rate of all servers in $L(\boldsymbol x)$. By the elementary renewal theorem, the long-run completion rate of a server that sees the job-mix $\{\pi_j\}$ is $1/m_i$. Hence the long-run fraction of completions handled by server $i\in L(\boldsymbol x)$ equals
\[
\psi_i(\boldsymbol x)
:=\frac{1/m_i}{\sum_{k\in L(\boldsymbol x)} 1/m_k},
\]
and the vector $\{\psi_i(\boldsymbol x)\}_{i\in L(\boldsymbol x)}$ is a probability vector, nonnegative and summing to one.

For an arrival from subregion $j$, in the heavy-load regime, the identity of the server that actually responds is distributed according to the completion-frequency weights $\psi_i(\boldsymbol x)$ 
. Therefore, the heavy-load expected travel time for a job from subregion $j$ is
\[
\mathbb{E}[t\mid j] \;=\; \sum_{i\in L(\boldsymbol x)} \psi_i(\boldsymbol x)\, t_{ij}.
\]
Averaging over the arrival mix $\{\pi_j\}$, the location-dependent contribution to the mean response time per arrival in the heavy-load limit is
\begin{align*}
\sum_{j\in\mathcal J}\pi_j\mathbb{E}[t\mid j]
&= \sum_{j\in\mathcal J}\pi_j\sum_{i\in L(\boldsymbol x)} \psi_i(\boldsymbol x)\, t_{ij}
= \sum_{i\in L(\boldsymbol x)} \psi_i(\boldsymbol x)\Bigl(\sum_{j\in\mathcal J}\pi_j t_{ij}\Bigr)\\
&= \frac{1}{{\sum_{j\in \mathcal{J}}\lambda_j}}\sum_{i\in L(\boldsymbol x)} \psi_i(\boldsymbol x)\Bigl(\sum_{j\in\mathcal J}\lambda_j t_{ij}\Bigr)
= \frac{1}{{\sum_{j\in \mathcal{J}}\lambda_j}}\sum_{i\in L(\boldsymbol x)} \psi_i(\boldsymbol x) V_i - \sum_{i\in L(\boldsymbol x)}\psi_i(\boldsymbol x)\tau_i,
\end{align*}
so that, apart from the additive convex-combination term involving the turnout times $\left\{\tau_i\right\}$, the location-dependent component reduces to the convex combination $\tfrac{1}{{\sum_{j\in \mathcal{J}}\lambda_j}}\sum_{i\in L(\boldsymbol x)}\psi_i(\boldsymbol x)V_i$ of the $\{V_i\}$.

Because any convex combination of values is at least their minimum, for any multiset $L$ and probability weights $\{\psi_i\}_{i\in L}$, we have
\[
\sum_{i\in L}\psi_i V_i \;\ge\; \min_{i\in L} V_i,
\]
with equality only if all weights are on an index attaining the minimum. Minimizing over all choices of $L$ yields
\[
\inf_{L:\,|L|=p}\sum_{i\in L}\psi_i V_i \;\ge\; \inf_{L:\,|L|=p}\min_{i\in L} V_i
\;=\; \min_{i\in\mathcal I} V_i =: V_{i^*}.
\]

Finally, if all $p$ units are colocated at some site $m$, then every colocated server faces the same service-time mixture and hence has the same mean service time $m_i$. In this case, the shares $\psi_i$ are equal across the $p$ colocated servers, and since each such server also has the same cost $V_i=V_m$, the convex combination reduces to $V_m$. Choosing $m=i^*$ (a 1-Median) attains the lower bound $V_{i^*}$. Therefore, in the heavy-load limit, the location-dependent part of the mean response time is minimized by colocating all $p$ units at any 1-Median site $i^*$. 
\Halmos
\endproof

\subsection{Proof of Theorem~\ref{thm:linear-regret}}\label{ssec_linear_regret}

The key technique used in our regret analysis is the \textit{posterior contraction} property in high-dimensional Bayesian linear models. In particular, we control the $\ell_1$-distance between the posterior sample $\hat{\alpha}_t$ and the true parameter $\alpha$ at each time step $t \in [T]$. 

The argument proceeds in three stages. First, we establish sparsity recovery by bounding the posterior mass, leveraging the horseshoe's exponential tail decay properties. Second, we control prediction error through the compatibility condition and the prior's concentration around sparse vectors. Finally, we convert these results into $\ell_1$-norm estimation error bounds using the sparsity induced by the prior's shrinkage behavior.

Let $\boldsymbol{\alpha}$ denote the true parameter vector. The bandit environment generates context-reward pairs ($\boldsymbol{x}_t, y_t$), where the contexts are drawn from arm-specific distributions and may depend on the history of actions. The joint distribution $\mathcal{Q}_t$ over $\left(\boldsymbol{\alpha}, \boldsymbol{x}_t, y_t\right)$ is defined by combining the horseshoe prior with the data-generating process 
$$\alpha_k \mid \beta_k, \tau, \sigma \sim \mathcal{N}(0,\beta_k^2\tau^2\sigma^2), \qquad (\boldsymbol{X}_t, y_t) \mid \boldsymbol{\alpha} \sim \text{HS}_t(\boldsymbol{\alpha}, \tau, \mathcal{P}_\epsilon),$$
where $\mathrm{HS}_t$ denotes the horseshoe prior of contexts and rewards at time $t$ given $\boldsymbol{\alpha}$, $\beta_k$ and $\tau$ are local and global hypervariances, and the noise distribution $\mathcal{P}_\epsilon$, as given in \eqref{eqn:horseshoe-prior}.


Let $\mathcal{D}=\{x\in\{0,1\}^N:\|x\|_0=p\}$ and let the surrogate feature map
$\phi:\mathcal{D}\to\{0,1\}^d$ collect the intercept, the $N$ linear terms and all
$\binom{N}{2}$ pairwise interactions, so $D=1+N+\tfrac{N(N-1)}{2}$.
For $t\ge 1$, after choosing $x_t\in\mathcal{D}$ the observation is
\[
y_t=\phi(x_t)^\top \alpha^* + \varepsilon_t.
\]

At round $t$, Thompson sampling (TS) draws $\tilde\alpha_t$ from the current posterior and
selects
\begin{equation}
x_t \in \argmin_{x\in\mathcal{D},\,\|x\|_0=p}\ \phi(x)^\top \tilde\alpha_t.
\label{eq:TS-argmin}
\end{equation}
The (surrogate) regret is $r_T=\sum_{t=1}^T \big(\phi(x_t)^\top\alpha^* - \phi(x^*)^\top\alpha^*\big)$
where $x^*\in\argmin_{x\in\mathcal{D}}\phi(x)^\top\alpha^*$.

\begin{lemma}[Per-round regret controlled by posterior $\ell_1$-radius]\label{lem:l1-reg}
For the TS action \eqref{eq:TS-argmin} and any realization of $\tilde\alpha_t$,
\[
r_t:=\phi(x_t)^\top\alpha^*-\phi(x^*)^\top\alpha^* \ \le\ 2\,\|\tilde\alpha_t-\alpha^*\|_1.
\]
Consequently, $\E[r_t\mid\mathcal{H}_{t-1}] \le 2\,\E_{\Pi(\cdot\mid\mathcal{H}_{t-1})}\|\alpha-\alpha^*\|_1$.
\end{lemma}

\proof{Proof.}
By optimality of $x_t$ for $\tilde\alpha_t$,
$\phi(x_t)^\top\tilde\alpha_t \le \phi(x^*)^\top\tilde\alpha_t$.
Hence
\[
\begin{aligned}
r_t
&= \underbrace{\phi(x_t)^\top\alpha^*-\phi(x_t)^\top\tilde\alpha_t}_{\le \|\alpha^*-\tilde\alpha_t\|_1}
 + \underbrace{\phi(x_t)^\top\tilde\alpha_t-\phi(x^*)^\top\tilde\alpha_t}_{\le 0}
 + \underbrace{\phi(x^*)^\top\tilde\alpha_t-\phi(x^*)^\top\alpha^*}_{\le \|\tilde\alpha_t-\alpha^*\|_1}\\
&\le 2\,\|\tilde\alpha_t-\alpha^*\|_1,
\end{aligned}
\]
where we used $|\phi(x)^\top v|\le \|v\|_1$ since $\|\phi(x)\|_\infty\le 1$ for all $x\in\mathcal{D}$.
Taking conditional expectation over the posterior of $\alpha$ given $\mathcal{H}_{t-1}$ yields the claim.
\Halmos \endproof

\begin{lemma}[Posterior contraction under fixed design and GL shrinkage priors]\label{lem:fixed-design}
Let $(X_t,y_t)_{t=1}^n$ be a fixed-design regression with $y_t = x_t^\top \alpha^* + \varepsilon_t$,
$\varepsilon_t$ i.i.d.\ Gaussian, and design matrix $X=[x_1,\ldots,x_n]^\top$ satisfying the RE condition
$\|X v\|_2/\sqrt{n}\ge \kappa \|v\|_2$ for all $v\in\mathcal{C}(s_0)$.
Consider a \emph{global--local heavy-tailed shrinkage prior} on $\alpha$ that (i) allocate polynomially heavy tails and (ii) put sufficiently large mass in any neighborhood of $0$.
Then there exists a universal constant $C>0$ (depending only on the prior class, $\kappa$ and the noise level)
such that for $n$ large enough,
\begin{equation}
\E\Big[\ \E_{\Pi(\cdot\mid X,y)} \|\alpha-\alpha^*\|_2\ \Big]\ \le\ C\,\sigma \sqrt{\frac{s_0 \log D}{n}},
\qquad
\E\Big[\ \E_{\Pi(\cdot\mid X,y)} \|\alpha-\alpha^*\|_1\ \Big]\ \le\ C\,\sigma\, s_0 \sqrt{\frac{\log D}{n}}.
\label{eq:fixed-design-rate}
\end{equation}
The outer expectation is w.r.t.\ the data. In particular, the horseshoe prior belongs to this class and enjoys the same rates.
\end{lemma}

\proof{Proof.}
We invoke the general posterior contraction theory for heavy-tailed global--local shrinkage priors in
high-dimensional linear regression: under Assumption~\ref{assump:compatibility} and Gaussian noise,
the posterior contracts around $\alpha^*$ at the (nearly) minimax rates
$\|\alpha-\alpha^*\|_2 = O_\Pi\!\big(\sqrt{(s_0\log D)/n}\big)$ and
prediction error $n^{-1}\|X(\alpha-\alpha^*)\|_2^2=O_\Pi\!\big((s_0\log D)/n\big)$, 
uniformly over $\ell_0$-sparse vectors. 
This is established for a broad class of continuous shrinkage priors (including the horseshoe)
by~\cite{song2023nearly} using testing and prior mass arguments adapted to global--local
mixtures with sufficiently heavy tails and spikes at $0$. Concretely, their conditions (flat-heavy tails,
sufficient prior mass near zero, and suitable hyper-prior on the global scale) are satisfied by the
horseshoe prior; see also summaries in~\cite{bhadra2019prediction}.
Applying Theorem~2.2 in~\cite{song2023nearly} yields
\[
\E\Big[\Pi\big(\|\alpha-\alpha^*\|_2 \ge c_0\,\sigma \sqrt{(s_0\log D)/n}\,\big|\,X,y\big)\Big]
\ \le\ e^{-c_1 s_0 \log D},
\]
for constants $c_0,c_1>0$ depending on $(\kappa,\sigma)$ and the prior family.
An integration-by-parts argument gives the bound on the posterior \emph{mean} loss:
for any $r>0$,
\[
\E_{\Pi}\|\alpha-\alpha^*\|_2
\le r + \int_r^\infty \Pi\big(\|\alpha-\alpha^*\|_2 > u \,\big|\, X,y\big)\,du,
\]
and choosing $r\asymp \sigma\sqrt{(s_0\log D)/n}$ plus the exponential tail yields the first inequality
in \eqref{eq:fixed-design-rate}. The $\ell_1$ bound follows from $\|\cdot\|_1\le \sqrt{s_0}\|\cdot\|_2$
on the error cone (standard in sparse regression under RE). This proves \eqref{eq:fixed-design-rate}.
\Halmos \endproof

\begin{lemma}[Horseshoe posterior contraction along the TS trajectory, on the good event]\label{lem:hs-on-policy}
Conditioning on $\mathcal{E}_t$ (the RE event for $\Phi_{t-1}$ with constant $\kappa$ as stated in Assumption~\ref{assump:compatibility}), for some $C_{\mathrm{HS}}>0$ and all $t\ge t_0$, we have
\begin{equation}
\E\Big[\ \E_{\Pi(\cdot\mid\mathcal{H}_{t-1})}\|\alpha-\alpha^*\|_1\ \Big|\ \mathcal{E}_t\Big]
\ \le\ C_{\mathrm{HS}}\,\sigma\, s_0 \sqrt{\frac{\log D}{t}}\ .
\label{eq:HS-rate-conditional}
\end{equation}
\end{lemma}

\proof{Proof.}
Conditioning on $\mathcal{H}_{t-1}$ and $\mathcal E_t$, the regression is fixed-design with $n=t-1$ samples
and design matrix $X=\Phi_{t-1}$ satisfying RE with constant $\kappa$.
Applying Lemma~\ref{lem:fixed-design} with $n=t-1$ gives
\[
\E\Big[\ \E_{\Pi(\cdot\mid\mathcal{H}_{t-1})}\|\alpha-\alpha^*\|_1\ \Big|\ \mathcal{E}_t\Big]
\ \le\ C\,\sigma\, s_0 \sqrt{\frac{\log D}{t-1}}
\ \le\ 2C\,\sigma\, s_0 \sqrt{\frac{\log D}{t}}\ ,
\]
for $t$ large enough; rename $C_{\mathrm{HS}}:=2C$ to obtain the stated Lemma.
\Halmos \endproof

\begin{lemma}[Deterministic bound on regret under the bad event]\label{lem:bad-event}
For every $t\ge 1$ and every history,
\[
r_t\ \le\ 2\,\|\alpha^*\|_1\ \le\ 2 W.
\]
Consequently, $\E\big[r_t\cdot \one_{\mathcal E_t^c}\big]\ \le\ 2 W\,\Prob(\mathcal E_t^c)$.
\end{lemma}

\proof{Proof.}
For any $x\in\mathcal D$, by (A2), $|\phi(x)^\top \alpha^*|\le \|\alpha^*\|_1$. Hence
\[
r_t = \phi(x_t)^\top\alpha^*-\phi(x^*)^\top\alpha^*
\ \le\ \abs{\phi(x_t)^\top\alpha^*}+\abs{\phi(x^*)^\top\alpha^*}
\ \le\ 2\,\|\alpha^*\|_1\ \le\ 2 W.
\]
Multiplying by $\one_{\mathcal E_t^c}$ and taking expectation yields the second claim.
\Halmos \endproof

Decompose
\[
\E[r_t] \ =\ \E\big[r_t \one_{\mathcal E_t}\big] + \E\big[r_t \one_{\mathcal E_t^c}\big].
\]
For the first term, by Lemma~\ref{lem:l1-reg} and then Lemma~\ref{lem:hs-on-policy},
\[
\E\big[r_t \one_{\mathcal E_t}\big]
\ \le\ 2\,\E\Big[\ \E_{\Pi(\cdot\mid\mathcal{H}_{t-1})}\|\alpha-\alpha^*\|_1 \cdot \one_{\mathcal E_t}\ \Big]
\ \le\ 2 C_{\mathrm{HS}}\,\sigma\, s_0 \sqrt{\frac{\log D}{t}}.
\]
For the second term, Lemma~\ref{lem:bad-event} yields
\[
\E\big[r_t \one_{\mathcal E_t^c}\big]\ \le\ 2 W\,\Prob(\mathcal E_t^c)\ \le\ 2 W\,e^{-ct}\cdot 2
\ =\ 4 W\,e^{-ct},
\]
for $t\ge t_0$ (absorbing the constant in $c$ if desired). Summing $t=1$ to $T$,
\[
\E[R(T)] \ \le\ 2 C_{\mathrm{HS}}\,\sigma\, s_0 \sum_{t=1}^T \sqrt{\frac{\log D}{t}}
\ +\ \sum_{t=1}^T 4 W\,e^{-ct}.
\]
Using $\sum_{t=1}^T t^{-1/2}\le 2\sqrt{T}$ gives
\[
\E[R(T)] \ \le\ 4 C_{\mathrm{HS}}\,\sigma\, s_0 \sqrt{T\,\log D}
\ +\ \underbrace{\sum_{t=1}^{\infty} 4 W\,e^{-ct}}_{=:C'_{\mathrm{bad}}<\infty}.
\]
Rename $C_0:=4 C_{\mathrm{HS}}$ and $C_{\mathrm{bad}}:=C'_{\mathrm{bad}}$ to conclude.

\subsection{Proof of Lemma \ref{lemma_rand_swap}}
\label{ssec_swap_lemma}
\proof{Proof.}

We begin by proving the first part of the lemma. Recall that $H(\boldsymbol{x}, \boldsymbol{x}') = H(\boldsymbol{x}, \boldsymbol{x}')^+ + H(\boldsymbol{x}, \boldsymbol{x}')^-$ and
\begin{equation}
     H(\boldsymbol{x}, \boldsymbol{x}')^+ = w\left(\left[\boldsymbol{x} \cap \neg \boldsymbol{x}'\right]\right), \qquad H(\boldsymbol{x}, \boldsymbol{x}')^- = w\left(\left[\neg \boldsymbol{x} \cap \boldsymbol{x}'\right]\right).
\end{equation}
When $m=1$, we have by definition $\xi^m(\boldsymbol{x}) = \xi(\boldsymbol{x}) = \{x_1,\cdots, x_j, \cdots, x_i, \cdots, x_{N}\}$ and $x_i\not=x_j$. Without loss of generality, let $x_i=1$ and $x_j=0$. Then, we have 
\begin{equation}
\begin{aligned}
     \boldsymbol{x} \cap \neg \xi(\boldsymbol{x}) &= \{x_1,\cdots, x_i, \cdots, x_j, \cdots, x_{N}\}\cap \neg\{x_1,\cdots, x_j, \cdots, x_i, \cdots, x_{N}\}\\
     &= \{x_1,\cdots, x_i, \cdots, x_j, \cdots, x_{N}\}\cap\{\neg x_1,\cdots, \neg x_j, \cdots, \neg x_i, \cdots, \neg x_{N}\}\\
     &= \{x_1 \cap \neg x_1,\cdots, x_i \cap \neg x_j, \cdots, x_j \cap \neg x_i, \cdots, x_{N} \cap \neg x_N\}
\end{aligned}
\end{equation}
We note that $x_k \cap \neg x_k = 0$ for all $k$. In addition, $x_i \cap \neg x_j = 1\cap 1=1$ and $x_j \cap \neg x_i = 0\cap 0=0$. Thus, we have 
\begin{equation}
    H(\boldsymbol{x}, \xi(\boldsymbol{x}))^+ = 1.
\end{equation}
Following the same logic, we can determine that $H(\boldsymbol{x}, \xi(\boldsymbol{x}))^- = 1, $ and $H(\boldsymbol{x}, \xi(\boldsymbol{x})) = H(\boldsymbol{x}, \xi(\boldsymbol{x}))^+ + H(\boldsymbol{x}, \xi(\boldsymbol{x}))^- = 2$. 

When $m\geq 2$, we first show that $H(\boldsymbol{x}, \xi^m(\boldsymbol{x})) \leq 2 \min \{p, N-p\}$ using induction. When $p\leq N/2$, we have $\min \{p, N-p\}=p$, and we show that $H(\boldsymbol{x}, \xi^m(\boldsymbol{x})) \leq 2p$. This holds because in order to increase the Hamming distance by two, we need to swap a new unit with a new empty site. In this case, we have more empty sites than units, so the number of units $p$ becomes the bottleneck, and the maximum number of new swaps is $2p$. The same logic applies when $p> N/2$, and in this case $H(\boldsymbol{x}, \xi^m(\boldsymbol{x})) \leq 2(N-p)$. Therefore, we must have $H(\boldsymbol{x}, \xi^m(\boldsymbol{x})) \leq 2 \min \{p, N-p\}$.

Assume that $H(\boldsymbol{x}, \xi^m(\boldsymbol{x})) \in \llbracket 0, 2\min\{m, p, N-p\} \rrbracket \cap \mathbb{N}_{2k}$, where $\llbracket a, b \rrbracket$ indicates the interval of all integers between $a$ and $b$ included and $\mathbb{N}_{2k}$ is the set of all even natural numbers. This holds for $m=2$, because $H(\boldsymbol{x}, \xi^2(\boldsymbol{x})) \in \{0 , 2, 4\}$. If in both steps, the same two units $x_i$ and $x_j$ are swapped, $H(\boldsymbol{x}, \xi^2(\boldsymbol{x}))=0$; if in the second step, $x_i$ is swapped with another unit $x_m$, then $H(\boldsymbol{x}, \xi^2(\boldsymbol{x}))=2$; if in the second step, two other units $x_m$ and $x_n$ are swapped, then $H(\boldsymbol{x}, \xi^2(\boldsymbol{x}))=4$. 

For the $m+1$-st swap, $H(\boldsymbol{x}, \xi^{m+1}(\boldsymbol{x})) \in \llbracket 0, 2\min\{(m+1), p, N-p\} \rrbracket \cap \mathbb{N}_{2k}$, because when two other units are swapped, the Hamming distance increases by two, as shown above. The maximum Hamming distance cannot exceed $\max \{p, N-p\}$.
\Halmos

\subsection{Proof of Theorem \ref{thm:convergence}}
\label{ssec_optimality_proof}

To prove Theorem \ref{thm:convergence}, we first present the following lemma. 
\begin{lemma}
\label{lemma:FTR_solution}
   Let $N$ and $p$ be the number of candidate locations and available units, respectively, in the $p$-Hypercube problem, and let $d$ be the edge-length of a feasible trust region (FTR). Then, the number of feasible solutions in the FTR is given by $\sum_{i=0}^{\lfloor d/2 \rfloor} {p \choose i}{N-p \choose i}$. 
\end{lemma}
\proof{Proof.}
Consider a feasible trust region (FTR) with edge-length $d$. The FTR contains all feasible solutions that are at most a Hamming distance of $\lfloor d/2 \rfloor$ away from the center of the FTR $\boldsymbol{x}^{c}$. To see this, note that we can generate a new feasible solution from the center solution $\boldsymbol{x}^{c}$ by moving $i$ units from $p$ occupied locations to $i$ empty locations from the $N-p$ set. This results in a Hamming distance of $2i$ between the new solution and the center solution. So, the maximum possible Hamming distance between any feasible solution and the center solution is $2\lfloor d/2 \rfloor$, the nearest even number that is less than $d$. 

Thus, the FTR contains all feasible solutions that are at most a Hamming distance of $\lfloor d/2 \rfloor$ away from the center solution. To count the number of such feasible solutions, we count the number of ways to choose $i$ units to move from $p$ occupied locations to $i$ empty locations from the $N-p$ locations, for $i=0,1,\ldots,\lfloor d/2 \rfloor$. There are ${p \choose i}$ ways to choose the occupied locations to move units from, and ${N-p \choose i}$ ways to choose the vacant locations to move units to. Therefore, the number of feasible solutions within the FTR is given by $\sum_{i=0}^{\lfloor d/2 \rfloor} {p \choose i}{N-p \choose i}$.
\Halmos \endproof

\proof{Proof of Theorem \ref{thm:convergence}.}
\color{black}
We prove the theorem in two parts. First, we show that ${\boldsymbol{x}_t}$ generated by our approach converges to the global optimum of $p$-Hypercube. Then, we prove that the algorithm converges in a finite number of iterations. 

Recall that $\text{OPT}_H$ is the global optimum of $p$-Hypercube, i.e., the minimum value of $f(\boldsymbol{x})$ over all feasible solutions $\boldsymbol{x} \in \mathcal{S}$. By definition, we have $\text{OPT}_H\leq g_t = \min_{k\leq t} f(\boldsymbol{x}_{k}) \leq f(\boldsymbol{x}_t)$ for all $t$. Thus, $\lim_{t\rightarrow\infty}g_t \geq \text{OPT}_H$. On the other hand, since our approach generates a new feasible solution in each iteration, we have $g_{t+1} \leq g_t$ for all $t$. This means that the sequence $\left\{g_t\right\}_{t=1}^{\infty}$ is monotonically non-increasing, and the actual objective function $f(\boldsymbol{x_t})$ and the constructed function $g_t$ are bounded below by $\text{OPT}_H$.  Therefore, by the monotone convergence theorem \citep{bibby1974axiomatisations}, $\left\{g_t\right\}_{t=1}^{\infty}$  converges to the global optimum of the objective function, i.e., $\lim_{t \rightarrow \infty} g_t = \text{OPT}_H$. This confirms the algorithm's convergence to the global optimum.

Next, we prove that the algorithm always converges in a finite number of iterations. The only possible cause for the algorithm's failure to converge in finite steps is if the local search method becomes trapped in some feasible trust region (FTR) and is unable to initiate a restart and move to a new FTR with new solutions. We prove this situation will not arise.

Assume, by contradiction, that there exists a feasible trust region $\mathrm{FTR}{d}\left(\boldsymbol{x}^{c}\right)$, characterized by a center $\boldsymbol{x}^c$ and edge-length $d$, defined as $\mathrm{FTR}_{d}\left(\boldsymbol{x}^{c}\right)=\left\{\boldsymbol{x}\in \{0,1\}^N \mid H(\boldsymbol{x}, \boldsymbol{x}^{c}) \leq d \text{  and  } \sum_{i=1}^N x_i = p \right\}$, such that the restart mechanism of our method fails to trigger. This implies that in this FTR, the edge-length $d$ fails to shrink below 2 according to the design of our algorithm. However, by the design of our algorithm, the feasible trust region will shrink after a consecutive number of $n_f$ failures, which means that the new selected solution does not strictly improve, and each time is reduced by a factor of $\alpha_f<1$. The adaptive swapping search method will finish exploring a current FTR in $n_f \left\lceil\log {1/\alpha_{f}}d\right\rceil$ evaluations, a finite number, if no further improvements can be found. According to Lemma \ref{lemma:FTR_solution}, there are a finite number of solutions within each FTR, so the adaptive swapping method will eventually trigger a restart and explore a new FTR. Given a finite number of FTRs with a fixed edge-length $d$, our algorithm will always terminate, and thus, it converges within a finite number of iterations.
\Halmos
\endproof

\subsection{Proof of Theorem \ref{thm:bound}}
\label{ssec_bound}

To prove Theorem \ref{thm:bound}, we first introduce some preliminary definitions and lemmas. Following \cite{cover1999elements}, we define the \textit{information gain} ($\mathrm{IG}$) of a set of sampling points $\boldsymbol{Y}_t=\{y_1, y_2, \cdots, y_t\}$ with respect to a function $f$ as
\begin{equation}
    \mathrm{IG}\left(\boldsymbol{Y}_t, f\right)=\mathrm{H}\left(\boldsymbol{Y}_t\right)-\mathrm{H}\left(\boldsymbol{Y}_t \mid f\right).
\end{equation}
The first term, $\mathrm{H}\left(\boldsymbol{Y}_t\right)$, represents the Shannon entropy of the set of points $\boldsymbol{Y}_t$. This measures the level of uncertainty or randomness in the set. If the points are close to each other, then the entropy is low, indicating that there is little uncertainty. On the other hand, if the points are spread out, the entropy is high, indicating a high level of uncertainty. The second term, $\mathrm{H}\left(\boldsymbol{Y}_t \mid f\right)$, represents the conditional entropy of $\boldsymbol{Y}_t$ given the function $f$. This measures the remaining level of uncertainty in $\boldsymbol{Y}_t$ after $f$ has been observed. If $f$ is informative about $\boldsymbol{Y}_t$, then the conditional entropy is low, indicating that there is little uncertainty left. Conversely, if $f$ is not very informative about $\boldsymbol{Y}_t$, then the conditional entropy is high, indicating that there remains a lot of uncertainty.

Therefore, IG represents the reduction in uncertainty of the function $f$ obtained by observing the set of points $\boldsymbol{Y}_t$. 
To further derive the bound on regret, we require the expression for the entropy of a multivariate normal variable. The following lemma presents this expression. 
\begin{lemma}
    The entropy of a multivariate normal variable $\boldsymbol{Z}\sim \mathcal{N}(\boldsymbol{\mu}, \boldsymbol{\Sigma})$ is given by 
\begin{equation}
    \mathrm{H}(\boldsymbol{Z}) = \frac{1}{2} \log \det(2 \pi e \boldsymbol{\Sigma}),
\end{equation}
where $\det(\cdot)$ denotes the determinant of a matrix and $\boldsymbol{\Sigma}$ is the covariance matrix of $\boldsymbol{Z}$.
\end{lemma}

The covariance matrix of a multivariate normal distribution is symmetric and positive-definite, so its determinant is always positive. Therefore, the entropy of a set of multivariate normally distributed variables is always non-negative.

\begin{lemma}
The information gain of a set of sampling points $\boldsymbol{Y}_t$ with respect to $f$ follows
\begin{equation}
    \mathrm{IG}\left(\boldsymbol{Y}_t, f\right)=\frac{1}{2} \log \det(\boldsymbol{I}_t+\sigma^{-2} \boldsymbol{K}_{t}),
\end{equation}
where $\boldsymbol{K}_{t} = [k(\boldsymbol{x},\boldsymbol{x}^{\prime})]_{\boldsymbol{x},\boldsymbol{x}^{\prime}\in \boldsymbol{X}_t}$ is the covariance matrix of the posterior distribution at the $t$-th step, and $\sigma$ is the standard deviation of the evaluation noise.
\end{lemma}
\proof{Proof:}
Since both $\boldsymbol{Y}_t$ and $\boldsymbol{Y}_t \mid f$ are multivariate normal distributions, according to the property of Gaussian processes, we have $\mathrm{H}(\boldsymbol{Y}_t) = \frac{1}{2} \log \det(2 \pi e \boldsymbol{\Sigma}_1)$ and $\mathrm{H}(\boldsymbol{Y}_t \mid f) = \frac{1}{2} \log \det(2 \pi e \boldsymbol{\Sigma}_2)$, where $\boldsymbol{\Sigma}_1$ and $\boldsymbol{\Sigma}_2$ are the covariance matrices of $\boldsymbol{Y}_t$ and $\boldsymbol{Y}_t \mid f$, respectively. We note that $\boldsymbol{\Sigma}_1 = \boldsymbol{K}_t+\sigma^2 \boldsymbol{I}_t$ and $\boldsymbol{\Sigma}_2 = \sigma^2 \boldsymbol{I}_t$ according to the Gaussian process. We thus have 
\begin{equation}
\begin{aligned}
\mathrm{IG}\left(\boldsymbol{Y}_t, f\right) & = \mathrm{H}\left(\boldsymbol{Y}_t\right)-\mathrm{H}\left(\boldsymbol{Y}_t \mid f\right) \\
& = \frac{1}{2} \log \det(2 \pi e (\boldsymbol{K}_t+\sigma^2 \boldsymbol{I}_t)) - \frac{1}{2} \log \det(2 \pi e \sigma^2 \boldsymbol{I}_t)\\& = \frac{1}{2} \log \left((2 \pi e)^t \det(\boldsymbol{K}_t+\sigma^2 \boldsymbol{I}_t)\right) - \frac{1}{2} \log \left((2 \pi e)^t \sigma^{2t} \right)\\
& = \frac{1}{2} \log \det(2 \pi e (\boldsymbol{K}_t+\sigma^2 \boldsymbol{I}_t)) - \frac{1}{2} \log \det(2 \pi e \sigma^2 \boldsymbol{I}_t)\\& = \frac{1}{2} \log \left( \frac{\det(\boldsymbol{K}_t+\sigma^2 \boldsymbol{I}_t)}{\sigma^{2t}}\right) \\
& =\frac{1}{2} \log \det(\boldsymbol{I}_t+\sigma^{-2} \boldsymbol{K}_{t}).
\end{aligned}
\end{equation}

\Halmos
\endproof

We next follow \cite{srinivas2009gaussian} in defining the \textit{maximum information gain} $\kappa_t$ at step $t$ as 
\begin{equation}
    \kappa_t = \max_{\boldsymbol{X}_t}\mathrm{IG}\left(\boldsymbol{Y}_t, f\right),
\end{equation}
which is the maximum reduction in uncertainty of the function $f$ obtained by observing the set of points $\boldsymbol{Y}_t$ for $\boldsymbol{X}_t$. We note that finding $\kappa_t$ is NP-hard, but a simple greedy algorithm guarantees a near-optimal solution due to the submodularity of IG \citep{srinivas2009gaussian}. Specifically, 
\begin{equation}
\boldsymbol{x}_t=\underset{\boldsymbol{x} \in [0,1]^n}{\operatorname{argmax}} \  k_{t-1}(\boldsymbol{x}, \boldsymbol{x}),
\end{equation}
and we have 
\begin{equation}
    \mathrm{IG}\left(\boldsymbol{Y}_t, f\right) \geq (1-1/e)\max_{t^{\prime}\leq t} \kappa_{t^{\prime}}.
\end{equation}

The maximum information gain $\kappa_t$ is used to bound the average cumulative regret of our algorithm. Next, we show the bound of the maximum information gain of our kernel.

\begin{lemma}
\label{lemma:max_info_gain}
    For the kernel $k(\cdot, \cdot)$ that we defined in \eqref{eq:kernel}, the maximum information gain $\kappa_T$ is bounded by $\mathcal{O}(2^N \log T)$, i.e., 
    \begin{equation}
        \kappa_T = \mathcal{O}(2^N \log T).
    \end{equation}
\end{lemma}

\proof{Proof.}
We first show the maximum information gain of our kernel in the 1-dimensional case and then extend it to the $N$-dimensional case. Recall that our kernel is a combination of two kernels
\begin{equation}
    k(\boldsymbol{x},\boldsymbol{x}^{\prime}) = k_c(\boldsymbol{x},\boldsymbol{x}^{\prime}) + k_d(\boldsymbol{x},\boldsymbol{x}^{\prime}),
\end{equation}
where $k_c(\boldsymbol{x},\boldsymbol{x}^{\prime})=\exp ( \sum_{i=1}^{N} \ell_{i}\delta\left(x_{i}, x_{i}^{\prime}\right)/N )$, and  $k_d(\boldsymbol{x},\boldsymbol{x}^{\prime}) = (\tanh{\gamma})^{\frac{H(\boldsymbol{x}, \boldsymbol{x}^{\prime})}{2}}$. 

In the 1-dimensional case, the Hamming distance is $H(x, x^{\prime}) = 1 - \delta(x, x^{\prime})$. Hence, we can write $k_d(x, x^{\prime}) = (\tanh{\gamma})^{\frac{1}{2}(1 - \delta(x, x^{\prime}))}$. If we set the parameter $l = -\frac{1}{2}\ln{(\tanh{\gamma}})$ in the categorical kernel $k_c$, then $k_c(x, x^{\prime}) = (\tanh{\gamma})^{-\frac{1}{2} \delta(x, x^{\prime})}$. We can arbitrarily set this parameter $l$ in our proof because its value does not influence the bound of the maximum information gain. As a result, we can represent $k_d$ by $k_c$ as
\begin{equation}
k_d({x},{x}^{\prime}) = (\tanh{\gamma})^{\frac{1}{2}} k_c({x},{x}^{\prime}).
\end{equation}
Then, we can further represent our kernel $k$ by the kernel $k_c$ as
\begin{equation}
    k({x},{x}^{\prime}) = k_c({x},{x}^{\prime}) + k_d({x},{x}^{\prime}) = (1 + (\tanh{\gamma})^{\frac{1}{2}} )k_c({x},{x}^{\prime}).
\end{equation}
Following the similar logic in \cite{wan2021think}, we bound the maximum information gain of our kernel $k$ by bounding the following
\begin{equation}
\kappa_T(k) \leq 
\log \det(I_T+\sigma^{-2} K_T) = 
 \log \det(I_T+\sigma^{-2}(1+(\tanh{\gamma})^{\frac{1}{2}}) K_T^c),
\end{equation}
where $K_T = \left[k\left(x_i, x_j\right)\right]_{i, j=1}^T$ and $K_T^c = \left[k_c\left(x_i, x_j\right)\right]_{i, j=1}^T$ are the transposes of the kernel matrices of our kernel $k$ and $k_c$, respectively, and $\kappa_T( k)$ represents the maximum information gain of our kernel $k$ in the 1-dimensional case at iteration $T$. Following \cite{wan2021think}, we have $\kappa_T( k) = \mathcal{O}(\zeta \log{(1+\sigma^{-2}(1+(\tanh{\gamma})^{\frac{1}{2}}) T(\exp (l)+\zeta-1))}) = \mathcal{O}(\zeta \log T)$, where $\zeta$ is the number of values the variable $x_i$ is allowed to take.

Next we consider the $N$-dimensional case. We can construct the $N$-dimensional kernel $k^{(N)}$ by taking the product of $N$ identical $1$-dimensional kernels, each of which depends on a single input dimension \citep{duvenaud2014automatic}. We use Theorem 2 from \cite{krause2011contextual} to derive the maximum information gain for the 2-dimensional kernel $k^{(2)}$ as follows:
\begin{equation}
    \kappa_T(k^{(2)} ) \leq \zeta \kappa_T(k)+\zeta \log T \leq \mathcal{O}(\zeta^2 \log T).
\end{equation}
We extend this result to the $N$-dimensional case using induction as follows:
\begin{equation}
    \kappa_T(k^{(N)} ) \leq \zeta \kappa_T( k^{(N-1)})+\zeta \log T \leq \mathcal{O}(\zeta^N \log T).
\end{equation}

In our model, there are only two distinct values $x \in\{0,1\}$, so we have $\zeta = 2$. Therefore, the maximum information gain $\kappa_T$ is bounded by $\mathcal{O}(2^N \log T)$.
\Halmos\endproof



We have now introduced all the necessary background information and tools to prove Theorem \ref{thm:bound}. Our proof follows a similar logic as in \cite{srinivas2009gaussian} for the GP-UCB algorithm. 

\proof{Proof of Theorem \ref{thm:bound}:}
\color{black}
In the Gaussian process model $GP^*$, given the collected data up to restart $v$, we denote $\boldsymbol{X}_{v} = \{\boldsymbol{x}_{1}, \cdots, \boldsymbol{x}_{v}\}$, $\boldsymbol{Y}_{v} = \{y_{1}, \cdots, y_{v}\}$. For a new solution point $\boldsymbol{x}$, the Gaussian process surrogate model learns the posterior mean and variance, denoted by $\mu(\boldsymbol{x};\boldsymbol{X}_{v-1}, \boldsymbol{Y}_{v-1})$ and $\sigma(\boldsymbol{x};\boldsymbol{X}_{v-1}, \boldsymbol{Y}_{v-1})$, respectively, from the observed data $\boldsymbol{X}_{v-1}$ and $\boldsymbol{Y}_{v-1}$. According to the Gaussian process property, we have $f(\boldsymbol{x}) \sim N(\mu(\boldsymbol{x};\boldsymbol{X}_{v-1}, \boldsymbol{Y}_{v-1}), \sigma(\boldsymbol{x};\boldsymbol{X}_{v-1}, \boldsymbol{Y}_{v-1}))$. 

From the properties of the normal distribution, for any $c$, we can write:
\begin{equation}
\begin{aligned}
\operatorname{Pr}\{ \frac{f(\boldsymbol{x}) - \mu(\boldsymbol{x};\boldsymbol{X}_{v-1}, \boldsymbol{Y}_{v-1})}{\sigma(\boldsymbol{x};\boldsymbol{X}_{v-1}, \boldsymbol{Y}_{v-1})} >c\} & = \frac{1}{\sqrt{2\pi}} \int_c^{\infty} e^{-u^2 / 2} d u = \frac{1}{\sqrt{2\pi}} \int_c^{\infty} e^{-(u-c)^2 / 2 - uc + c^2 / 2} d u \\
& =e^{-c^2 / 2}\frac{1}{\sqrt{2\pi}} \int_c^{\infty} e^{-(u-c)^2 / 2-c(u-c)} d u.
\end{aligned}
\end{equation}

Let $r=u-c$ and let $z$ be a standard normal variable. Then, we have:
\begin{equation}
\begin{aligned}
\operatorname{Pr}\{ \frac{f(\boldsymbol{x}) - \mu(\boldsymbol{x};\boldsymbol{X}_{v-1}, \boldsymbol{Y}_{v-1})}{\sigma(\boldsymbol{x};\boldsymbol{X}_{v-1}, \boldsymbol{Y}_{v-1})} >c\} 
& =e^{-c^2 / 2}\frac{1}{\sqrt{2\pi}} \int_0^{\infty} e^{-r^2 / 2-cr} d r\\
& \leq e^{-c^2 / 2} \operatorname{Pr}\{z>0\}= \frac{e^{-c^2 / 2}}{2}.
\end{aligned}
\end{equation}

Similarly, we apply the same argument to the left tail of the distribution. We have
\begin{equation}
\operatorname{Pr}\{ \frac{f(\boldsymbol{x}) - \mu(\boldsymbol{x};\boldsymbol{X}_{v-1}, \boldsymbol{Y}_{v-1})}{\sigma(\boldsymbol{x};\boldsymbol{X}_{v-1}, \boldsymbol{Y}_{v-1})} <-c\}  \leq  \frac{e^{-c^2 / 2}}{2}.
\end{equation}

Therefore, let $c=\beta_v^{1 / 2}$, and we obtain 
\begin{equation}
    \operatorname{Pr}\{|f(\boldsymbol{x}) - \mu(\boldsymbol{x};\boldsymbol{X}_{v-1}, \boldsymbol{Y}_{v-1})|> \beta_v^{1 / 2} \sigma(\boldsymbol{x};\boldsymbol{X}_{v-1}, \boldsymbol{Y}_{v-1})\} \leq e^{-\beta_v / 2},
\end{equation}
and this holds for all feasible solutions $\boldsymbol{x}$. We apply the union bound and obtain 
\begin{equation}
    \operatorname{Pr}\{|f(\boldsymbol{x}) - \mu(\boldsymbol{x};\boldsymbol{X}_{v-1}, \boldsymbol{Y}_{v-1})|\leq \beta_v^{1 / 2} \sigma(\boldsymbol{x};\boldsymbol{X}_{v-1}, \boldsymbol{Y}_{v-1}), \forall \boldsymbol{x}\} \geq 1 - {N \choose p} e^{-\beta_v / 2},
\end{equation}
where $N$ is the number of available locations and $p$ is the number of available units. Let $\delta = {N \choose p} e^{-\beta_v / 2} \frac{\pi^2 v^2}{6}$ and we have  $\sum_{v=1}^{\infty} \frac{6}{\pi^2 v^2} = 1$. Then, using the union bound for all $v$, we have for any $\delta \in (0,1]$,
\begin{equation}
    \operatorname{Pr}\{|f(\boldsymbol{x}) - \mu(\boldsymbol{x};\boldsymbol{X}_{v-1}, \boldsymbol{Y}_{v-1})|\leq \beta_v^{1 / 2} \sigma(\boldsymbol{x};\boldsymbol{X}_{v-1}, \boldsymbol{Y}_{v-1}), \forall \boldsymbol{x}, \forall v\} \geq 1 - \delta,
\end{equation}
where $\beta_v = 2\log(\frac{{N \choose p}\pi^2 v^2}{6\delta})$. 


Recall that our algorithm obtains the center of the feasible trust region at the $v$-th restart $\boldsymbol{x}_{v}$ using the lower confidence bound (LCB) function with the Gaussian process $GP^*$, which we fit from the selected observations $(\boldsymbol{X}_{v-1}, \boldsymbol{Y}_{v-1})$. We obtain
{
\setlength{\abovedisplayskip}{3pt}
\setlength{\belowdisplayskip}{3pt}
\begin{equation}
    \boldsymbol{x}_{v}=\argmin_{\boldsymbol{x}} \mu(\boldsymbol{x};\boldsymbol{X}_{v-1}, \boldsymbol{Y}_{v-1})-\beta_{v}^{1 / 2} \sigma(\boldsymbol{x};\boldsymbol{X}_{v-1}, \boldsymbol{Y}_{v-1}).
\end{equation}
}
By definition of the LCB function, we have
{\setlength{\abovedisplayskip}{3pt}
\setlength{\belowdisplayskip}{3pt}
\begin{align*}
    &\mu\left(\boldsymbol{x}_v;\boldsymbol{X}_{v-1}, \boldsymbol{Y}_{v-1}\right)-\beta_{v}^{1 / 2} \sigma\left(\boldsymbol{x}_v;\boldsymbol{X}_{v-1}, \boldsymbol{Y}_{v-1}\right) \\&
    \leq \mu\left(\boldsymbol{x}^{*};\boldsymbol{X}_{v-1}, \boldsymbol{Y}_{v-1}\right)-\beta_{v}^{1 / 2} \sigma\left(\boldsymbol{x}^{*};\boldsymbol{X}_{v-1}, \boldsymbol{Y}_{v-1}\right) \leq f\left(\boldsymbol{x}^{*}\right),
\end{align*}
}
where $\boldsymbol{x}^{*}$ is the optimal solution and $ f\left(\boldsymbol{x}^{*}\right) = \text{OPT}_H$ is the optimal objective function value. We have 
{
\setlength{\abovedisplayskip}{3pt}
\setlength{\belowdisplayskip}{3pt}
\begin{equation}
f\left(\boldsymbol{x}_{v}\right) - \text{OPT}_H \leq \beta_{v}^{1 / 2} \sigma\left(\boldsymbol{x}_v;\boldsymbol{X}_{v-1}, \boldsymbol{Y}_{v-1}\right)-\mu\left(\boldsymbol{x}_v;\boldsymbol{X}_{v-1}, \boldsymbol{Y}_{v-1}\right)+f\left(\boldsymbol{x}_{v}\right). 
\end{equation}
}

Then, with probability greater than or equal to $1-\delta$, we bound the regret of the best observation $\boldsymbol{x}_{v}^{*}$ at the $v$-th restart as
\begin{equation}
    \begin{aligned}
    r_{v} & = f(\boldsymbol{x}_{v}^{*}) - \text{OPT}_H \leq f\left(\boldsymbol{x}_{v}\right) - \text{OPT}_H \\ 
    & \leq \beta_{v}^{1 / 2} \sigma\left(\boldsymbol{x}_{v};\boldsymbol{X}_{v-1}, \boldsymbol{Y}_{v-1}\right)-\mu\left(\boldsymbol{x}_{v};\boldsymbol{X}_{v-1}, \boldsymbol{Y}_{v-1}\right)+f\left(\boldsymbol{x}_{v}\right) \\
     & \leq \beta_{v}^{1 / 2} \sigma\left(\boldsymbol{x}_{v};\boldsymbol{X}_{v-1}, \boldsymbol{Y}_{v-1}\right) + \beta_{v}^{1 / 2} \sigma\left(\boldsymbol{x}_v;\boldsymbol{X}_{v-1}, \boldsymbol{Y}_{v-1}\right)\\
     & = 2\beta_{v}^{1 / 2} \sigma\left(\boldsymbol{x}_v;\boldsymbol{X}_{v-1}, \boldsymbol{Y}_{v-1}\right).
    \end{aligned}
\end{equation}

Then, taking the square of both sides, we obtain, with probability greater than or equal to $1-\delta$,
\begin{equation}
r_{v}^2 = (f(\boldsymbol{x}_{v}) - \text{OPT}_H)^2  \leq 4 \beta_{v}^{} \sigma^2\left(\boldsymbol{x}_v;\boldsymbol{X}_{v-1}, \boldsymbol{Y}_{v-1}\right).
\end{equation}

Since $\beta_v$ is non-decreasing, setting $C=1 / \log \left(1+\sigma^{-2}\right)$, we have
\begin{equation}
    \begin{aligned}
\beta_{v} \sigma^2\left(\boldsymbol{x}_v;\boldsymbol{X}_{v-1}, \boldsymbol{Y}_{v-1}\right) &\leq \beta_{V} \sigma^2\left(\boldsymbol{x}_v;\boldsymbol{X}_{v-1}, \boldsymbol{Y}_{v-1}\right)=\beta_{V} \sigma^2 (\sigma^{-2} \sigma^2\left(\boldsymbol{x}_v;\boldsymbol{X}_{v-1}, \boldsymbol{Y}_{v-1}\right))\\
&\leq \beta_{V}  C\log \left(1+\sigma^{-2} \sigma^2\left(\boldsymbol{x}_v;\boldsymbol{X}_{v-1}, \boldsymbol{Y}_{v-1}\right)\right).
    \end{aligned}
\end{equation}

Using Lemma 5.3 in \cite{srinivas2009gaussian}, we show that the information gain before the $V$-th restart is
\begin{equation}
    \mathrm{IG}\left(\boldsymbol{Y}_{V}, f\right)=\frac{1}{2} \sum_{v=1}^{V} \log \left(1+\sigma^{-2} \sigma^{2}\left(\boldsymbol{x}_{v} ;\boldsymbol{X}_{v-1}, \boldsymbol{Y}_{v-1}\right)\right).
\end{equation}

Then, with probability greater than or equal to $1-\delta$, we have
\begin{equation}
    \begin{aligned}
    \sum_{v=1}^{V} r_{v}^{2} &\leq \sum_{v=1}^{V} 4 \beta_{v} \sigma^2\left(\boldsymbol{x}_v;\boldsymbol{X}_{v-1}, \boldsymbol{Y}_{v-1}\right) \leq \sum_{v=1}^{V} 4 C \beta_{V} \log(1 + \sigma^{-2} \sigma^2\left(\boldsymbol{x}_v;\boldsymbol{X}_{v-1}, \boldsymbol{Y}_{v-1}\right)) \\
     &=  4 C \beta_{V} \sum_{v=1}^{V} \log(1 + \sigma^{-2} \sigma^2\left(\boldsymbol{x}_v;\boldsymbol{X}_{v-1}, \boldsymbol{Y}_{v-1}\right))\\
     &=  8 C \beta_{V} \mathrm{IG}\left(\boldsymbol{Y}_V, f\right)\leq  8 C \beta_{V} \kappa_V.
    \end{aligned}
\end{equation}

We now bound the average regret before the $V$-th restart, denoted by $\bar{r}_{V} = \frac{1}{V} \sum_{v=1}^{V} r_v$, using the Cauchy-Schwarz inequality:
\begin{equation*}
    \bar{r}_{V} = \frac{1}{V}\sum_{v=1}^{V}r_v \leq \sqrt{\frac{1}{V} \sum_{v=1}^{V} r_{v}^{2}} \leq \sqrt{8 C\beta_{V}  \kappa_{V} / V}, 
\end{equation*}
where $\kappa_{V} = \mathcal{O}(2^N \log V)$ as shown in Lemma \ref{lemma:max_info_gain}.
\Halmos\endproof

\color{black}

\end{appendices}

\end{document}